\documentclass{elsarticle}

\usepackage{bm}
\usepackage{amsmath}
\usepackage{amsfonts}
\usepackage{amsthm}
\usepackage{amssymb}
\usepackage{graphicx}
\usepackage{paralist}
\usepackage{algorithm}
\usepackage{pstool}
\usepackage{wrapfig}
\usepackage{algpseudocode}
\usepackage{mathtools}
\usepackage{fullpage}
\usepackage{etoolbox}

\usepackage{tikz}
\usepackage{pgfplots}
\usepackage{pgfplotstable, filecontents}
\pgfplotsset{compat=1.9}
\usetikzlibrary{pgfplots.groupplots}
\usepgfplotslibrary{fillbetween}
\usetikzlibrary{spy,calc,fit,matrix,arrows,automata,positioning,shapes}
\usepackage{pgfplotstable, booktabs}
\pgfplotsset{select coords between index/.style 2 args={
    x filter/.code={
        \ifnum\coordindex<#1\fi
        \ifnum\coordindex>#2\fi
    }
}}

\newtheorem*{Rem}{Remark}

\newcommand{\func}[3]{\ensuremath{#1 : #2 \rightarrow #3}}

\newcommand{\norm}[1]{\ensuremath{\left\| #1 \right\|}}
\newcommand{\optunc}[2]{\underset{#1}{\text{minimize}} ~~ #2}
\newcommand{\optcona}[3]{
\begin{aligned}
& \underset{#1}{\text{minimize}}
& & #2 \\
& \text{subject to} & & #3
\end{aligned}}

\newcommand{\pder}[2]{\ensuremath{\frac{\partial #1}{\partial #2}}} 

\newcommand{\Acal}{\ensuremath{\mathcal{A}}}

\newcommand{\Ecal}{\ensuremath{\mathcal{E}}}
\newcommand{\Fcal}{\ensuremath{\mathcal{F}}}
\newcommand{\Gcal}{\ensuremath{\mathcal{G}}}

\newcommand{\Lcal}{\ensuremath{\mathcal{L}}}

\newcommand{\Ocal}{\ensuremath{\mathcal{O}}}
\newcommand{\Pcal}{\ensuremath{\mathcal{P}}}

\newcommand{\Tcal}{\ensuremath{\mathcal{T}}}

\newcommand{\Vcal}{\ensuremath{\mathcal{V}}}

\newcommand{\Rbb}{\ensuremath{\mathbb{R} }}

\newcommand\Ibm{{\ensuremath{\bm{I}}}}

\newcommand\Mbm{{\ensuremath{\bm{M}}}}

\newcommand\Wbm{{\ensuremath{\bm{W}}}}
\newcommand\Xbm{{\ensuremath{\bm{X}}}}

\newcommand\ebm{{\ensuremath{\bm{e}}}}

\newcommand\rbm{{\ensuremath{\bm{r}}}}

\newcommand\ubm{{\ensuremath{\bm{u}}}}

\newcommand\xbm{{\ensuremath{\bm{x}}}}

\newcommand\lambdabold{{\ensuremath{\boldsymbol{\lambda}}}}

\newcommand\etabold{{\ensuremath{\boldsymbol{\eta}}}}

\newcommand\phibold{{\ensuremath{\boldsymbol{\phi}}}}

\newcommand\zerobold{\ensuremath{\mathbf{0}}}

\newbool{fastcompile}
\setbool{fastcompile}{false}

\begin{document}

\title{An optimization-based approach for high-order accurate
       discretization of conservation laws with discontinuous solutions}

\author[rvt1,rvt2]{M.~J.~Zahr\fnref{fn1}\fnref{fn2}\corref{cor1}}
\ead{mjzahr@lbl.gov}

\author[rvt1,rvt3]{P.-O.~Persson\fnref{fn3}}
\ead{persson@berkeley.edu}

\address[rvt1]{Mathematics Group, Lawrence Berkeley National Laboratory,
               1 Cyclotron Road, Berkeley, CA 94720, United States}
\address[rvt2]{Department of Aerospace and Mechanical Engineering,
               University of Notre Dame, Notre Dame, IN, 46556, United States}
\address[rvt3]{Department of Mathematics, University of California, Berkeley,
               Berkeley, CA 94720, United States}
\cortext[cor1]{Corresponding author}
\fntext[fn1]{Luis W. Alvarez Postdoctoral Fellow,
             Computational Research Division,
             Lawrence Berkeley National Laboratory}
\fntext[fn2]{Assistant Professor, Department of Aerospace and Mechanical
             Engineering, University of Notre Dame}
\fntext[fn3]{Associate Professor, Department of Mathematics, University of
             California, Berkeley}

\begin{keyword}
 $r$-adaptivity, %
 shock tracking, %
 high-order methods, %
 discontinuous Galerkin, %
 full space PDE-constrained optimization, %
 transonic and supersonic flow
\end{keyword}

\begin{abstract}
This work introduces a novel discontinuity-tracking framework for
resolving discontinuous solutions of conservation laws with high-order
numerical discretizations that support inter-element solution
discontinuities, such as discontinuous Galerkin or finite volume methods.
The proposed method aims to align inter-element boundaries with
discontinuities in the solution by deforming the computational mesh.
A discontinuity-aligned mesh ensures the discontinuity is represented through
inter-element jumps while smooth basis functions interior to elements are only
used to approximate smooth regions of the solution, thereby avoiding Gibbs'
phenomena that create well-known stability issues. Therefore, very coarse
high-order discretizations accurately resolve the piecewise smooth solution
throughout the domain, provided the discontinuity is tracked. Central to the
proposed discontinuity-tracking framework is a discrete PDE-constrained
optimization formulation that simultaneously aligns the computational
mesh with discontinuities in the solution and solves the discretized
conservation law on this mesh. The optimization objective is taken as a
combination of the the deviation of the finite-dimensional solution from
its element-wise average and a mesh distortion metric to simultaneously
penalize Gibbs' phenomena and distorted meshes. It will be shown that our
objective function satisfies two critical properties that are required for this
discontinuity-tracking framework to be practical: (1) possesses a local
minima at a discontinuity-aligned mesh and (2) decreases monotonically to
this minimum in a neighborhood of approximately $\Ocal(h)$, whereas other
popular discontinuity indicators fail to satisfy the latter. Another
important contribution of this work is the observation that traditional
reduced space PDE-constrained optimization solvers that repeatedly solve
the conservation law at various mesh configurations are not viable in this
context since severe overshoot and undershoot in the solution,
i.e., Gibbs' phenomena, may make it impossible to solve the discrete
conservation law on non-aligned meshes.
Therefore, we advocate a gradient-based, full space solver where the mesh
and conservation law solution converge to their optimal values simultaneously
and therefore never require the solution of the discrete conservation law on
a non-aligned mesh. The merit of the proposed method is demonstrated on a
number of one- and two-dimensional model problems including the $L^2$
projection of discontinuous functions, Burgers' equation with a discontinuous
source term, transonic flow through a nozzle, and supersonic flow around a
bluff body. We demonstrate optimal $\Ocal(h^{p+1})$ convergence rates in
the $L^1$ norm for up to polynomial order $p=6$ and show that accurate
solutions can be obtained on extremely coarse meshes.
\end{abstract}

\maketitle

\section{Introduction}\label{sec:intro}

Even with the continued advancement of computational hardware and algorithms,
it is apparent that important developments are still required to improve the
predictive capability of the computational simulation tools that are used for
design. For example, it is widely believed that higher fidelity is required
for problems with propagating waves, turbulent fluid flow, nonlinear
interactions, and multiple scales \cite{wang2013high}. This has resulted in a
significant interest in high-order accurate methods, such as discontinuous
Galerkin (DG) methods \cite{cockburn01rkdg,hesthaven08dgbook}, which have the
potential to produce more accurate solutions on coarser meshes than
traditional compatible discretizations, such as mimetic finite differences
and finite volumes.

One of the remaining problems with the DG methods is their sensitivity to
under-resolved features, in particular for non-linear problems where the
spurious oscillations often cause a break-down of the numerical solvers.  This
is exacerbated for problems with shocks, where the natural dissipative
mechanisms introduced by DG methods through jump terms is
insufficient to stabilize the solution for high-order approximations.  Since
shocks are present in many important problems in fields such as aerospace,
astrophysics, and combustion, the lack of efficient ways to handle them is a
fundamental challenge for the wide adoption of these new numerical methods.

One of the most straight-forward approaches for stabilizing high-order methods
in the presence of shocks is to identify the elements close to the shocks
using a sensor, and reduce the corresponding approximation orders
\cite{BauOden,burbeau01limiter}. This increases the numerical dissipation
which helps stabilize the solution, and in the extreme case of piecewise
constant solutions will completely avoid any oscillations. An obvious
drawback with this strategy is that it will typically result in a first-order
scheme in the affected elements. This can partially be addressed by using
$h$-adaptivity and remeshing \cite{dervieux03adaptation,casoni12shock}; however,
generating efficient anisotropic meshes in 3D is still a difficult task and
the resulting meshes typically have a very large number of degrees of freedom
around the shocks. Further issues with these methods include moving shocks,
which require constant remeshing and solution transfer, and a difficulty for
implicit solvers with large timesteps since order reduction strategies
typically cannot be differentiated for a nonlinear Newton solver. A related
class of methods use limiters based on the ENO/WENO schemes
\cite{eno1,weno1,weno2} to limit the high-order representation of the solution
inside the elements \cite{cockburnshuV,DGWENO1,DGWENO2,Krivodonova,shu2010}.
While these schemes can be remarkably robust and stable for very strong shocks,
they suffer from spurious oscillations, have difficulties with fully
unstructured meshes, particularly in 3D, and are typically only used
for explicit timestepping.

Another popular way to stabilize problems with shocks is to explicitly add
artificial viscosity to the governing equations \cite{neumann50shocks}.  In
\cite{bassirebay95shocks,BauOden,hartmann02shock}, this was applied to the DG
method using element-wise viscosity based on the residual of the Euler
equations.  To address issues with the consistency of the resulting scheme, a
proper discretization of the artificial viscosity was used in
\cite{persson06shock} together with a resolution-based indicator based on the
decay-rate of the high-order terms. The resulting scheme could smoothly
resolve shocks using subcell resolution, and the continuous nature of both the
sensor and the viscosity made it possible to get fully converged steady-state
solutions and accurate time-stepping using implicit methods.  Further
developments of this approach include some refinements of the parameter
selection \cite{klockner2011shock}, the solution of an auxiliary PDE for the
viscosity field \cite{barter08}, and the extension to implicit time-stepping
for transient problems with smooth viscosity fields \cite{persson13shock}.  An
alternative approach for the shock sensor is to somehow incorporate the
physics of the problem
\cite{jameson81jst,lele2009viscosity,moro14thesis}. Overall, these artificial
viscosity approaches are widely used, but like limiting they require
anisotropic $h$-adaptivity to be competitive
\cite{alauzet16anisotropic,yano12anisotropic}, which again is difficult in 3D,
particularly for moving shocks, and produces large meshes.

In this work, we propose a fundamentally different approach based on the idea
of shock tracking or shock fitting
\cite{shubin1981steady, shubin1982steady, bell1982fully, harten1983self,
      rosendale1994floating,
      trepanier1996conservative, zhong98tracking,
      taghaddosi1999adaptive, baines2002multidimensional,
      roe2002adaptive, glimm2003conservative,
      palaniappan2008sub}
using $r$-adaptivity. For generality,
we will use the term \emph{discontinuity tracking} in this work, although our
motivating applications involve shock waves. Since the solution basis in DG
methods (as well as the finite volume methods) naturally support
discontinuities between each element, it is in principle possible to keep an
unmodified high-order DG discretization for problems with strong
discontinuities, provided the element faces can be aligned with
the discontinuities of the problem. However, this is
a very challenging problem since the unknown discontinuity surfaces can have
complex geometries and if the elements are even just slightly misaligned,
the discrete equations cannot be solved due to spurious oscillations. Instead,
we employ a new optimization formulation that aims to align faces of the
computational mesh with discontinuities and solve the discrete conservation
law. The optimization formulation is classified as fully discrete
$r$-adaptive PDE-constrained optimization since the discretized conservation
law defines nonlinear equality constraints and the state vector and nodal
positions of the computational mesh are optimization variables. In this
framework, the objective function must possess a local minima at a
discontinuity-aligned mesh and must monotonically approach this minima
in a neighborhood of radius approximately $h/2$ for a solver to
successfully locate such a minima. Here, we propose a simple robustness-based
discontinuity indicator that will be experimentally shown to satisfy both
criteria, while a number of intuitive physics- and error-based discontinuity
indicators fail to satisfy the monotonicity property. Finally, we
use a full space approach to solve the PDE-constrained optimization problem
whereby the discrete solution of the conservation law and the nodal positions
of the mesh \emph{simultaneously} converge to their optimal values without
ever requiring the solution of the discrete conservation law on a non-aligned
mesh. The popular alternative, the reduced space approach, is not a viable
solver as it requires the solution of the discrete conservation law on
non-aligned meshes, which may be impossible due to Gibbs' instabilities.
We demonstrate our method can produce high-order accurate solutions on
remarkably coarse meshes, for a range of problems in 1D and 2D.

The remainder of the paper is organized as follows. Section~\ref{sec:disc}
introduces the governing system of steady, inviscid conservation laws on a
parametrized domain and its discretization using a discontinuous Galerkin
method. Since deformation of the computational domain is a key ingredient in
this work, a parametrized domain deformation is introduced at the continuous
level and the conservation laws are recast on a fixed reference domain. In
the end, Section~\ref{sec:disc} reduces the conservation law to a discrete
nonlinear system of equations that depend on the discrete state vector and
computational mesh, which is the point of departure for the proposed
discontinuity-tracking framework, introduced in Section~\ref{sec:track}.
The proposed PDE-constrained, $r$-adaptive formulation is presented in
Section~\ref{sec:track-optform}. The objective function used in this work,
introduced in Section~\ref{sec:track-obj}, is a combination of new
discontinuity indicator and a standard mesh distortion indicator. In this
section, we show the indicator possesses a minima at a discontinuity-aligned
mesh and decreases monotonically to this minima, whereas the latter condition
fails for some existing indicators. Section~\ref{sec:track-solver}
discusses the full space solver required for the proposed
discontinuity-tracking and a detailed discussion on why a traditional
reduced space PDE-constrained optimization approach is not sufficient.
Section~\ref{sec:track-practical} discusses a number of practical details
including implementation, initialization of the nonlinear optimization solver,
and robustness and efficiency of the method. Finally,
Section~\ref{sec:num-exp} presents a number of one- and two-dimensional
test problems that demonstrate optimal $\Ocal(h^{p+1})$ convergence up to
$p = 6$ and accurate solutions are obtained on extremely coarse meshes.

\section{Governing equations and high-order numerical discretization}
\label{sec:disc}
Consider a general system of $N_c$ conservation laws, defined on the physical
domain $\Omega \subset \Rbb^d$,
\begin{equation} \label{eqn:claw-phys}
  \nabla\cdot \Fcal(U) = 0 \quad \text{in}~~\Omega,
\end{equation}
where $U(x) \in \Rbb^{N_c}$ is the solution of the system of conservation laws
at $x \in \Omega$ and $\Fcal(U) \in \Rbb^{N_c\times d}$ is the physical flux. While we
solely consider steady conservation laws, i.e., $\Omega$ is a $d$-dimensional
spatial domain and the solution is independent of time, the conservation law in
(\ref{eqn:claw-phys}) and the proposed numerical method encapsulate the
unsteady case where $\Omega$ is the $d$-dimensional space-time domain and
the spatial domain is $(d-1)$-dimensional. Furthermore, we assume the solution
$U$ contains discontinuities, in which case the conservation law
(\ref{eqn:claw-phys}) holds away from these discontinuities.

The proposed optimization-based method for tracking discontinuities,
introduced in Section~\ref{sec:track}, is built upon existing numerical
discretizations that possess the following properties:
\begin{inparaenum}[1)]
 \item represents a stable and convergent discretization of the conservation
       law in (\ref{eqn:claw-phys}),
 \item allows for deformation of the computational domain, and
 \item employs a solution basis that supports discontinuities between
       computational cells or elements.
\end{inparaenum}
While the proposed method will be valid for any numerical discretization that
satisfies these requirements, such as finite volume or discontinuous Galerkin
methods, we focus on high-order DG methods given their potential to deliver
accurate solutions on very coarse discretizations, provided the discontinuities
are tracked.

The remainder of this section will detail the discretization of the
conservation law (\ref{eqn:claw-phys}) using DG such that it reduces
to the discrete form
\begin{equation} \label{eqn:claw-disc}
 \rbm(\ubm,\,\xbm) = \zerobold
\end{equation}
where $\ubm \in \Rbb^{N_\ubm}$ is the discrete representation of the
conservation law state $U$, $\xbm \in \Rbb^{N_\xbm}$ is the discrete
representation of the conservation law domain $\Omega$, and $\rbm$ is
the discretized conservation law. A standard nodal, isoparametric DG method
is used for the discretization, with special attention given to treatment of
the domain deformation and the numerical fluxes, which must be selected
carefully in a tracking method since the inter-element jumps do not tend to
zero on the discontinuity surface. We will discuss the formulation
of the proposed optimization-based tracking method in the context of the general
discrete conservation law (\ref{eqn:claw-disc}). As a result, our developments
will be applicable to a broad range of admissible discretizations with only
minor modifications.
Section~\ref{sec:track} will mention the modifications required to use a
finite volume method as the underlying discretization.

\subsection{Transformed conservation law from deformation
            of physical domain}
\label{sec:disc:transf}
Before introducing a discretization of (\ref{eqn:claw-phys}) it is convenient
to explicitly treat deformations to the domain of the conservation law $\Omega$.
Given the divergence-form of the conservation law in (\ref{eqn:claw-phys}),
i.e., no explicit time derivative term, the domain deformation is handled
directly. That is, the physical domain can be taken as the result of a
parametrized diffeomorphism applied to a reference domain
(Figure~\ref{fig:dom-map})
\begin{equation} \label{eqn:dom-map}
 \Omega = \Gcal(\Omega_0,\,\mu),
\end{equation}
where $\Omega_0 \subset \Rbb^d$ is a fixed reference domain,
$\mu \in \Rbb^{N_\mu}$ is a vector of parameters, and
$\func{\Gcal}{\Rbb^d \times \Rbb^{N_\mu}}{\Rbb^d}$ is the parametrized
diffeomorphism defining the domain mapping. Under the domain mapping
(\ref{eqn:dom-map}), the conservation law becomes
\begin{equation} \label{eqn:claw-phys-map}
 \nabla \cdot \Fcal(U) = 0 \quad \text{in}~~\Gcal(\Omega_0,\,\mu).
\end{equation}
\ifbool{fastcompile}{}{
\begin{figure}
  \centering
  \includegraphics[width=2.5in]{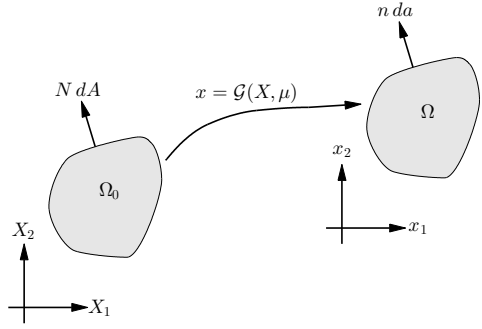}
  \caption{Mapping between reference and physical domains.}
  \label{fig:dom-map}
\end{figure}
}
For convenience, the conservation law on the physical domain $\Omega$ is
transformed to a conservation law on the reference domain $\Omega_0$ using
the procedure in \cite{persson2009dgdeform} to yield
\begin{equation} \label{eqn:claw-ref}
 \nabla_X \cdot F(u,\,\mu) = 0 \quad \text{in}~~\Omega_0
\end{equation}
where $\nabla_X$ denotes spatial derivatives with respect to the reference
domain $\Omega_0$ with coordinates $X$. The transformed state vector, $u$,
and flux, $F$, take the form
\begin{equation} \label{eqn:transf-stvc-flux}
u = g_\mu U, \qquad F(u,\,\mu) = g_\mu \Fcal(g_\mu^{-1}u)G_\mu^{-T}
\end{equation}
where $G_\mu(X) = \pder{}{X}\Gcal(X,\,\mu)$ is the deformation gradient
of the domain mapping and $g_\mu(X) = \det G_\mu(X)$ is the
Jacobian. The transformed conservation law in
(\ref{eqn:claw-ref})-(\ref{eqn:transf-stvc-flux}) is, by construction,
identical to the Arbitrary Lagrangian-Eulerian formulation of the steady
conservation law in (\ref{eqn:claw-phys-map}). For details regarding the
derivation of the transformed equations, the reader is referred to
\cite{persson2009dgdeform}.

While a numerical scheme and implementation can be designed around
(\ref{eqn:claw-phys-map}) or (\ref{eqn:claw-ref}), we chose the formulation
on the reference domain in (\ref{eqn:claw-ref}) as it provides a convenient
framework to handle domain deformation. In this case, the governing equation
depends on the domain deformation \emph{solely} through the transformed
flux function and definition of the state vector. This detail leads to a
clean implementation, particularly when computing required derivative
quantities, as discussed in Section~\ref{sec:track-practical}.

\subsection{Discontinuous Galerkin discretization of
            transformed conservation law}
\label{sec:disc:dg}
In this section, a standard nodal discontinuous Galerkin method
\cite{cockburn01rkdg,hesthaven08dgbook} is used to
discretize the transformed conservation law (\ref{eqn:claw-ref}) with special
attention given to the numerical fluxes. Let $\Ecal_{h,p}$ represent
a discretization of the reference domain $\Omega_0$ into non-overlapping,
potentially curved, computational elements, where $h$ is a mesh element
size parameter and $p$ is the polynomial order associated with
the curved elements.
The DG construction begins by considering the element-wise weak form of the
governing equation in (\ref{eqn:claw-ref}) that results from multiplication
by a test function $\psi$, integration over a single element
$K \in \Ecal_{h,p}$, and application of the divergence theorem
\begin{equation} \label{eqn:claw-weak-elem}
 \int_K \nabla_X \cdot F(u,\,\mu) \cdot \psi \,dV =
 \int_{\partial K} \psi \cdot F(u,\,\mu) N\,dA - \int_K F(u,\,\mu) :
                   \nabla_X \psi\,dV = 0,
\end{equation}
where $N$ is the outward normal to the surface $\partial K$. The global
weak form, upon which DG methods are built, arises from the summation of
the local residuals over all elements in $\Ecal_{h,p}$
\begin{equation} \label{eqn:claw-weak-dg}
 \sum_{K\in\Ecal_{h,p}}\int_{\partial K} \psi \cdot F(u,\,\mu) N\,dA -
 \int_{\Omega_0} F(u,\,\mu) : \nabla_X \psi\,dV = 0.
\end{equation}
At this point, a numerical flux, $F^*$, is introduced in the first term of
(\ref{eqn:claw-weak-dg})
\begin{equation} \label{eqn:claw-weak-dg-numflux}
 \sum_{K\in\Ecal_{h,p}}\int_{\partial K} \psi\cdot F^*(u,\,\mu,\,N)\,dA -
 \int_{\Omega_0} F(u,\,\mu) : \nabla_X \psi\,dV = 0
\end{equation}
to ensure the flux is single-valued along $\partial K$ where $u$ is
multi-valued. In DG methods, the numerical flux must be chosen such that it
is consistent to ensure local conservation \cite{cockburn2001rkdg}, i.e.,
$F^*(\bar{u},\,\mu,\,N) = F(\bar{u},\,\mu) N$ where $\bar{u}(X)$ is single-valued, and
leads to a stable discretization. Given that inter-element jumps tend to zero
under refinement, the particular choice of numerical flux is less important,
provided it is stable and consistent \cite{cockburn01rkdg, wheatley2010riemann}.
However, recall that the goal of this
work is to \emph{track} discontinuities with element faces, which implies
inter-element jumps at the discontinuity will not tend to zero under
refinement and care must be taken when choosing a numerical flux such that
it is consistent with the governing equation at discontinuities. In this
work, Roe's method \cite{roe1981approximate} with the Harten-Hyman entropy
fix \cite{harten1983self} is used for the numerical fluxes at interior
faces and the appropriate boundary conditions determine the numerical
fluxes on faces that intersect $\partial\Omega$. The Harten-Hyman entropy
fix is used to ensure non-physical rarefaction shocks do not appear.

To establish the finite-dimensional form of (\ref{eqn:claw-weak-dg-numflux}),
we introduce the isoparametric finite element space of piecewise polynomial
functions associated with the mesh $\Ecal_{h,p}$:
\begin{equation*}
 \begin{aligned}
  \Vcal_{h,p} &= \left\{v \in [L^2(\Omega_0)]^{N_c} \mid
          \left.v\right|_K \circ \Tcal_K \in [\Pcal(K_0)]^{N_c}
          ~\forall K \in \Ecal_{h,p}\right\} \\
  \hat\Vcal_{h,p}(\mu) &= \left\{v \in [L^2(\Omega)]^{N_c} \mid
          \left.v\right|_{\Gcal(K,\,\mu)} \circ \Tcal_{\Gcal(K,\,\mu)}
          \in [\Pcal(K_0)]^{N_c}
          ~\forall K \in \Ecal_{h,p}\right\}
 \end{aligned}
\end{equation*}
where $\Pcal_p(K_0)$ is the space of polynomial functions of degree at most
$p \geq 1$ on the parent element $K_0$ and $K = \Tcal_K(K_0)$ defines a
mapping from the parent element to element $K \in \Ecal_{h,p}$. For notational
brevity, we assume all elements map from a single parent element. The
parametrized space $\hat\Vcal_{h,p}(\mu)$ will be used in later sections to
define an $L^2$ projection onto a piecewise polynomial space attached to the
physical domain $\Omega = \Gcal(\Omega_0,\,\mu)$. The
Petrov-Galerkin weak form in (\ref{eqn:claw-weak-dg-numflux}) becomes: find
$u_{h,p} \in \Vcal_{h,p}$ such that for all $\psi_{h',p'} \in \Vcal_{h',p'}$
\begin{equation} \label{eqn:claw-findim-dg-pg}
 r_{h,p}^{h',p'}(u_{h,p},\,\mu) \coloneqq
                       \sum_{K\in\Ecal_{h,p}}\int_{\partial K}
                       \psi_{h',p'}\cdot F^*(u_{h,p},\,\mu,\,N) \,dA -
 \int_{\Omega_0} F(u_{h,p},\,\mu) : \nabla_X \psi_{h',p'}\,dV = 0,
\end{equation}
where $r_{h,p}^{h',p'}$ is the finite-dimensional residual of
(\ref{eqn:claw-weak-dg-numflux}) corresponding to the trial space
$\Vcal_{h,p}$ and test space $\Vcal_{h',p'}$. In this work, these spaces
are taken to be the same and we have the Galerkin weak form: find
$u_{h,p} \in \Vcal_{h,p}$ such that for all $\psi_{h,p} \in \Vcal_{h,p}$
\begin{equation} \label{eqn:claw-findim-dg}
 \sum_{K\in\Ecal_{h,p}}\int_{\partial K}
 \psi_{h,p}\cdot F^*(u_{h,p},\,\mu,\,N) \,dA -
 \int_{\Omega_0} F(u_{h,p},\,\mu) : \nabla_X \psi_{h,p}\,dV = 0.
\end{equation}
The Petrov-Galerkin residual in (\ref{eqn:claw-findim-dg-pg}) will be used
in Section~\ref{sec:track-obj-compare} to define an accuracy-based
discontinuity indicator. To obtain the equivalent algebraic
representation of (\ref{eqn:claw-findim-dg}), we introduce a basis
$\{\varphi_i\}_{i=1}^{N_p}$ for $\Pcal_p(K_0)$, where $N_p$ is the number of
basis functions. While this
can be any valid basis, a flexibility afforded by the DG framework, it is
convenient for the proposed shock tracking method to work with a
\emph{nodal basis}. Define a set of nodes $\{\xi_j\}_{j=1}^{N_p}$ within the
parent element $K_0$ and let $\varphi_i$ be the Lagrange polynomial associated
with nodes $\xi_i$ over element $K_0$ with the property
$\varphi_i(\xi_j) = \delta_{ij}$. For brevity, we introduce notation for the
Lagrange basis defined on the reference element $K \in \Ecal_{h,p}$:
$\varphi_i^K(X) = \varphi_i(\Tcal_K^{-1}(X))$. The finite-dimensional
solution, $u_{h,p}$, in each element $K \in \Ecal_{h,p}$ is written in terms of
its discrete expansion coefficients as
\begin{equation} \label{eqn:trial-fcn-findim}
 \left.u_{h,p}(X)\right|_K = \sum_{i=1}^{N_p}\ubm_i^K\varphi_i^K(X),
\end{equation}
where $\ubm_i^K \in \Rbb^{N_c}$ is the solution at node $i$ of
element $K$. 
Under the isoparametric assumption, the physical coordinates are expanded in
the nodal basis as
\begin{equation}
 \left.x_{h,p}(X)\right|_K = \sum_{i=1}^{N_p}\xbm_i^K\varphi_i^K(X),
\end{equation}
where $\xbm_i^K \in \Rbb^d$ is the coordinate of node $i$ of
element $K$. It will also be convenient in later sections to introduce the
corresponding nodal coordinates in the reference domain, $\Xbm_i^K \in \Rbb^d$.
From this expansion, the deformation gradient, $G$, and Jacobian, $g$,
required to define the transformed state vector and fluxes are
\begin{equation}
 \begin{aligned}
  \left.G_{h,p}(X)\right|_K &= \sum_{i=1}^{N_p}
                               \xbm_i^K\pder{\varphi_i^K}{X}(X) \\
  \left.g_{h,p}(X)\right|_K &= \left.\det G_{h,p}(X)\right|_K.
 \end{aligned}
\end{equation}
Therefore, in the isoparametric setting, the mapping to the physical domain is
completely determined from the nodal positions of each element $\xbm_i^K$.
To ensure the domain mapping is continuous, nodal positions
co-located in the reference domain are required to be co-located in the
physical domain. This is accomplished by using a single parameter to define
the position of all nodes co-located in the reference domain, which is
equivalent to using the nodal positions of the \emph{continuous} high-order
mesh corresponding to $\Ecal_{h,p}$ as the parameter set. Let
$\xbm \in \Rbb^{N_\xbm}$ denote these positions, then $\xbm_i^K$ is
reconstructed from $\xbm$ by selecting the appropriate entry based on the
position of the corresponding node in the reference domain. The corresponding
nodal coordinates, denoted $\Xbm \in \Rbb^{N_\xbm}$, in the reference domain
are defined similarly. With this notation, the continuous domain mapping
is parametrized by $\mu = \xbm$, i.e., $\Omega = \Gcal(\Omega_0,\,\xbm)$.
The remainder of the document will use $\xbm$ to parametrize the domain
deformation instead of $\mu$. Finally, the integrals in
(\ref{eqn:claw-findim-dg}) are evaluated using high-order Gaussian quadrature
rules to yield the discrete form of the governing equations
\begin{equation} \label{eqn:claw-disc}
 \rbm(\ubm,\,\xbm) = 0,
\end{equation}
where $\ubm \in \Rbb^{N_\ubm}$ ($N_\ubm = N_cN_p|\Ecal_{h,p}|$) is the solution vector comprised of the coefficients $\ubm_i^K$ for all elements.
The discrete form of the governing equations in (\ref{eqn:claw-disc}) will
be the point of departure for the remainder of the document.
The proposed shock tracking method presented in the next section
will be applicable to any discretization of the governing equations that
assumes the general form (\ref{eqn:claw-disc}) and satisfies the requirements
outlined in the beginning of this section, such as finite volume methods.

\section{High-order discontinuity tracking via optimization-based
         $r$-adaptivity}
\label{sec:track}
In this section, we present a novel method designed to align
discontinuous features in a finite-dimensional solution basis with
features in the solution itself, with a driving application being
transonic and supersonic compressible flow where shock waves are present.
In the discretization setting outlined in Section~\ref{sec:disc}, e.g.,
discontinuous Galerkin or finite volumes, this amounts to aligning element
faces, where discontinuities are supported, with discontinuities in the
solution. With the discontinuous features tracked with element faces,
very coarse high-order discretizations effectively resolve the solution
throughout the domain. However, in general, the location of the
discontinuities is not known \emph{a-priori} and therefore the face alignment
cannot be done explicitly. To circumvent this difficulty that has plagued
previous attempts at discontinuity tracking, the proposed method is built
around three critical contributions:
\begin{enumerate}[1)]
 \item a PDE-constrained optimization formulation with optimization variables
       taken as the nodal positions of the \emph{high-order} mesh, i.e., an
       optimization-based $r$-adaptive framework, that allows elements
       to curve to accurately track the location of the discontinuity,
 \item an objective function that attains a (local) minimum at a
       discontinuity-aligned mesh and monotonically decreases to this
       minimum in a neighborhood of radius approximately $h/2$,
 \item a PDE-constrained optimization solver that does not require solution
       of the discrete PDE (\ref{eqn:claw-disc}) on meshes not aligned with
       the discontinuity as this cannot be done reliably in a numerical
       scheme due to Gibbs' phenomena around the discontinuity (see
       Section~\ref{sec:intro} for the complete discussion).
\end{enumerate}
The remainder of this section will detail these contributions and summarize
the proposed high-order discontinuity-tracking algorithm in its entirety.

\subsection{Optimization formulation for $r$-adaptivity}
\label{sec:track-optform}
The proposed method for high-order resolution of discontinuous solutions of
conservation laws reformulates the discrete nonlinear system in
(\ref{eqn:claw-disc}) as a PDE-constrained optimization problem. The
original nonlinear system formulation searches for a discrete solution $\ubm$
on a given and fixed mesh $\xbm$, whereas the proposed optimization
formulation searches for the discrete solution and mesh that minimize some
objective function, $\func{f}{\Rbb^{N_\ubm}\times\Rbb^{N_\xbm}}{\Rbb}$, and
satisfy the discrete PDE
\begin{equation} \label{eqn:claw-disc-opt}
 \optcona{\ubm\in\Rbb^{N_\ubm},\,\xbm\in\Rbb^{N_\xbm}}
         {f(\ubm,\,\xbm)}{\rbm(\ubm,\,\xbm) = \zerobold.}
\end{equation}
The objective function, discussed in detail in Section~\ref{sec:track-obj},
is constructed such that it is locally minimized at a discontinuity-aligned
mesh. The proposed method falls into the category of $r$-adaptive methods
given that it moves nodes to improve approximation quality; however, it has
stronger requirements than existing $r$-adaptive methods as we aim to align
mesh \emph{faces} with the discontinuity rather than just nodes. See
Figure~\ref{fig:mshparam-1dof-2d} (top left) for a case where several mesh
nodes are aligned with a discontinuity without any faces being aligned. Such
situations are avoided through the construction of the objective function in
Section~\ref{sec:track-obj} that is minimized when faces are aligned with
discontinuities.

While the general formulation of the optimization problem in
(\ref{eqn:claw-disc-opt}) is desirable in that it does not require any
\emph{a-priori} knowledge of the location of the discontinuity, it can be
cumbersome due to the large number of optimization variables and
difficult to avoid ill-conditioning and non-uniqueness in the
domain mapping, i.e., severe mesh distortion or inversion.
Ill-conditioning and non-uniqueness of the domain mapping will be addressed
in the construction of the objective function in Section~\ref{sec:track-obj}.
To reduce the number of optimization variables, we re-parametrize the
optimization problem in (\ref{eqn:claw-disc-opt}). Let
$\phibold \in \Rbb^{N_\phibold}$ be a vector of parameters such that
$N_\phibold \leq N_\xbm$ and $\func{\Acal}{\Rbb^{N_\phibold}}{\Rbb^{N_\xbm}}$ be
an injective mapping that maps $\phibold \mapsto \xbm$, i.e.,
$\xbm = \Acal(\phibold)$. For clarity, we use $\xbm(\phibold)$ in place of
$\Acal(\phibold)$. Then the optimization problem in (\ref{eqn:claw-disc-opt})
becomes
\begin{equation} \label{eqn:claw-disc-opt1}
 \optcona{\ubm\in\Rbb^{N_\ubm},\,\phibold\in\Rbb^{N_\phibold}}
         {f(\ubm,\,\xbm(\phibold))}
         {\rbm(\ubm,\,\xbm(\phibold)) = \zerobold.}
\end{equation}
In addition to reducing the number of optimization parameters, the mapping
$\Acal$ can incorporate mesh smoothing operations to improve mesh quality.
However, despite the advantages of re-parametrization, the fact that
$\Acal$ is not necessarily bijective limits the reachable mesh configurations
and care must be taken to ensure the parametrization is capable of aligning with
discontinuities. In this document, three different mesh parametrizations will
be considered: the entire mesh, a single node parametrization with piecewise
linear mesh smoothing, and a $n$-node parametrization with linear elastic
smoothing in $d$ dimensions. In the most general case where
all mesh nodes are used as optimization variables, the
re-parametrization map is the identity
map, $\Acal(\phibold)=\phibold$, or any other bijection, and
(\ref{eqn:claw-disc-opt1}) is equivalent to the original optimization problem.
The other mesh parametrizations used in this work will be discussed in
Sections~\ref{sec:track-obj}~and~\ref{sec:num-exp:invisc-2d}.

\subsection{Discontinuity tracking objective function}
\label{sec:track-obj}
A critical contribution of this work is the specific form of the objective
function in the PDE-constrained optimization problem in
(\ref{eqn:claw-disc-opt1}). For the proposed framework to successfully align
faces of the computational mesh with discontinuities in the solution the
objective function must
\begin{inparaenum}[1)]
 \item attain a local minimum at some discontinuity-aligned mesh and
 \item monotonically decrease to such a minima is a neighborhood of
       radius approximately $h/2$.
\end{inparaenum}
The first requirement ensures a local minima of (\ref{eqn:claw-disc-opt1})
leads to a desired mesh and the second requirement gives gradient-based
numerical optimization methods hope to locate it even in the worst-case
scenario where element faces are initially $h/2$ away from
discontinuities. Thus, we seek an objective function that is an
effective discontinuity indicator that monotonically decreases from
large values at non-aligned meshes to smaller values at aligned meshes.

A host of discontinuity indicators exist in the literature
\cite{baines2002multidimensional,persson06shock} as sensing
discontinuous features is a backbone of many numerical methods for either
capturing or tracking discontinuities. Two common classes of discontinuity
indicators are physics-based and error-based. Physics-based indicators
typically involve the residual associated with the integral form of the
conservation law, i.e., the Rankine-Hugoniot conditions
\cite{baines2002multidimensional}.
Error-based indicators use the magnitude of an approximate error measure
such as adjoint-based error estimators \cite{fidkowski2011review} or the
residual of the solution in a richer test space \cite{ainsworth1997posteriori}.
In this work, we propose a departure from these traditional discontinuity
indicators and consider an indicator driven purely by the \emph{robustness}
of the numerics \cite{persson06shock}. Recall from Section~\ref{sec:intro}
that high-order discretizations of conservation laws fail when attempting
to resolve discontinuous solutions due to severe spurious oscillations.
Therefore we propose an objective function that penalizes these oscillations
by considering the \emph{element-wise} distance of the numerical solution
from its mean value, i.e.,
\begin{equation} \label{eqn:obj-shk}
 f_{shk}(\ubm,\,\xbm) = h_0^{-2}
                        \sum_{K \in \Ecal_{h,p}}
                        \int_{\Gcal(K,\,\xbm)}
                        \norm{u_{h,p} - \bar{u}_{h,p}^K}_\Wbm^2 \, dV
\end{equation}
where the dependence of the finite dimensional solution $u_{h,p}$ on the
discrete representation is implied, $\bar{u}_{h,p}^K$ is the mean value
of $u_{h,p}$ over element $K$, $\Wbm \in \Rbb^{N_c\times N_c}$ is the
symmetric positive semi-definite matrix that defines the local semi-norm,
and $h_0$ is the length scale of the reference mesh $\Ecal_{h,p}$
\begin{equation} \label{eqn:elem-mean}
  \bar{u}_{h,p}^K = \frac{1}{|\Gcal(K,\,\xbm)|}
                \int_{\Gcal(K,\,\xbm)} u_{h,p}\,dV,\qquad
  |\Gcal(K,\,\xbm)| = \int_{\Gcal(K,\,\xbm)} dV,\qquad
  h_0 = \left(\frac{1}{|\Ecal_{h,p}|}\int_{\Omega_0} dV\right)^{1/d}.
\end{equation}
The $h_0^{-2}$ factor ensures the indicator scales as $\Ocal(1)$ instead of
$\Ocal(h^2)$.
This objective function does not explicitly depend on the governing equation,
only the approximation space $\Vcal_{h,p}$, nor does it take into account the
accuracy of the finite-dimensional solution, which marks a departure
from physics-based and error-based discontinuity indicators.

\begin{Rem}
The discontinuity indicator in (\ref{eqn:obj-shk}) will be identically zero,
and therefore a useless objective function, if a piecewise constant
polynomial approximation space is used ($p=0$), which is equivalent to a
finite volume scheme. In this case, the pointwise values of the finite
dimensional solution $u_{h,p}$ in (\ref{eqn:obj-shk}) can be replaced by its
non-limited reconstruction using the solution in neighboring elements
\cite{levequeFVM}.
\end{Rem}

The discontinuity indicator in (\ref{eqn:obj-shk}) is not well-suited as
the objective function in the discontinuity-tracking optimization setting
(\ref{eqn:claw-disc-opt1}) in its current form because it is agnostic to a poor
quality or inverted mesh that may arise from certain choices of $\xbm$.
Therefore we construct the objective function as a weighted combination
of (\ref{eqn:obj-shk}) and a function $\func{f_{msh}}{\Rbb^{N_\xbm}}{\Rbb}$
that penalizes mesh distortion
\begin{equation} \label{eqn:obj}
 f(\ubm,\,\xbm;\,\alpha) = f_{shk}(\ubm,\,\xbm) + \alpha f_{msh}(\xbm).
\end{equation}
In the remainder of the document the dependence on $\alpha$ will be dropped
unless explicitly required. In this work, we use the following mesh distortion
measures
\begin{equation} \label{eqn:obj-msh}
 f_{msh}(\xbm) =
\begin{cases}
 \displaystyle{h_0\sum_{K\in\Ecal_{h,p}} \left|\frac{h_0}{|\Gcal(K,\,\xbm)|} - 1\right|} & \quad d = 1 \\
 \displaystyle{h_0^d\sum_{K\in\Ecal_{h,p}} \frac{1}{|\Gcal(K,\,\xbm)} \int_{\Gcal(K,\,\xbm)} \left(\frac{\norm{G_{h,p}}_F^2}{\left(\det G_{h,p}\right)_+^{2/d}}\right)^r} & \quad \text{otherwise},
\end{cases}
\end{equation}
where $r = 2$ is used in this work. For the general case of $d > 1$, the
quantity in (\ref{eqn:obj-msh}) is widely accepted in the high-order meshing
community as reasonable metric to penalize mesh distortion and is well-suited
for optimization \cite{knupp01quality,roca16distortion}. However, this mesh distortion measure is not useful
for the case $d = 1$. Therefore we treat the case $d=1$ separately and define
mesh distortion as the deviation from a uniform mesh, modified to heavily
penalize zero-volume elements. The $h_0^d$ factors are again included to
ensure the $f_{msh}$ scales independently of $h$. Since both $f_{shk}$ and
$f_{msh}$ are independent of $h$ and $p$, the value of $\alpha$ that
appropriately balances discontinuity tracking and mesh regularization should
be mesh independent, i.e., independent of $h$ and $p$.

\subsubsection{Behavior of objective function: $L^2$ projection of discontinuous
               function}
\label{sec:track-obj-behavior}
In this section, we use a simple one- and two-dimensional model problem to
demonstrate the behavior of the proposed discontinuity indicator
(\ref{eqn:obj-shk}) near a discontinuity-aligned mesh and verify it possesses
the desired properties suggested at the beginning of
Section~\ref{sec:track-obj}, namely attains a local minimum at a
discontinuity-aligned mesh and monotonically approaches such a minima.
Additionally, we show the impact of the mesh regularization term
(\ref{eqn:obj-msh}) on the combined objective function (\ref{eqn:obj}) as a
function of the weighting parameter $\alpha$. Unlike the physics- and
error-based indicator, the proposed indicator is independent of
the conservation law, which provides the opportunity to study the indicator
\emph{without} the complication of a partial differential equation. Instead,
we consider the $L^2$ projection of a discontinuous function onto the piecewise
polynomial space $\hat\Vcal_{h,p}(\xbm)$.

A complete description of the $L^2$ projection in the finite-dimensional
setting of Section~\ref{sec:disc:dg} is provided in
Section~\ref{sec:num-exp:l2proj}. Here, we let $\ubm_H(\xbm)$ denote the
discrete representation of the $L^2$ projection of the discontinuous function
$\func{H}{\Omega}{\Rbb}$ onto the space $\hat\Vcal_{h,p}(\xbm)$ and
investigate the behavior of
\begin{equation} \label{eqn:obj-redsp-l2proj}
 j_\alpha(\phibold) = f(\ubm_H(\xbm(\phibold)),\,\xbm(\phibold);\,\alpha),
\end{equation}
where $\xbm(\phibold)$ denotes a parametrization of the computational mesh.
In order to graphically study the objective function, it is important to
use a single degree of freedom to parametrize the mesh. This is accomplished
in one spatial dimension ($d=1$) by using $\phibold$ to define the mapping
of a single point $\bar{x}$ in the reference domain $\Omega_0$ and mapping
the remaining points piecewise linearly between the fixed boundaries and
$\phibold$
\begin{equation} \label{eqn:mshparam1d}
x_{h,p}(X;\,\phibold) =
\begin{cases}
\displaystyle{x_l + (X-x_l)\frac{\phibold-x_l}{\bar{x}-x_l}} & X \leq\bar{x} \\
\displaystyle{x_r - (x_r-X)\frac{x_r-\phibold}{x_r-\bar{x}}} & X > \bar{x}
\end{cases}
\end{equation}
where $\Omega_0 = (x_l,\,x_r) \subset \Rbb$. In the discrete setting, the
mapping in (\ref{eqn:mshparam1d}) is applied to the nodal positions of the
reference domain $\Xbm$ to define the nodal positions of the physical
domain $\xbm$ as
\begin{equation} \label{eqn:mshparam1d-disc}
\xbm_i(\phibold) =
 \begin{cases}
   \displaystyle{x_l + (\Xbm_i-x_l)\frac{\phibold-x_l}{\bar{x}-x_l}}
                                                  \quad & \Xbm_i \leq \bar{x} \\
   \displaystyle{x_r - (x_r-\Xbm_i)\frac{x_r-\phibold}{x_r-\bar{x}}}
                                                  \quad & \Xbm_i > \bar{x}
 \end{cases}
\end{equation}
for $j \in \{1,\,\dots,\,N_\xbm\}$. This single degree of freedom mesh
parametrization possesses the obvious property
$\xbm_i(\bar{x}) = \Xbm_i$, which implies that the mesh quality will be
perfect in terms of the distortion metric $f_{msh}$ in (\ref{eqn:obj-msh})
at $\phibold=\bar{x}$ if the reference nodes are uniformly distributed.
Additionally, if the reference nodes are uniformly distributed, the
parametrization in (\ref{eqn:mshparam1d-disc}) inherently incorporates ideal
smoothing since the resulting nodes will be uniformly spaced in the intervals
$[x_l,\,\bar{x}]$ and $(\bar{x},\,x_r]$. A depiction of the single degree
of freedom mesh parametrization is provided in
Figure~\ref{fig:mshparam-1dof-1d}. This figure also shows the range that
$\phibold$ will sweep to study the behavior of $j_\alpha(\phibold)$
next.

\ifbool{fastcompile}{}{
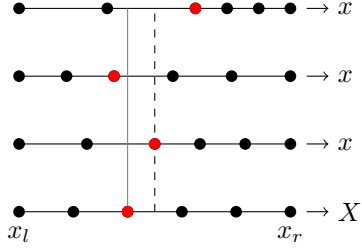
\begin{figure}
 \centering
 \begin{tikzpicture}
\begin{axis}[
axis equal image,
axis lines=none,
xmin=-0.2,
ymin=-0.2,
ymax=0.95,
xmax=1.7]
\addplot [solid, forget plot]
coordinates {
(  0.00000000,   0.00000000)
(  0.20000000,   0.00000000)
(  0.40000000,   0.00000000)
(  0.60000000,   0.00000000)
(  0.80000000,   0.00000000)
(  1.00000000,   0.00000000)};

\addplot [red, forget plot]
coordinates {
(  0.40000000,   0.00000000)};

\addplot [solid, mark options={solid}, mark=*, mark size=3, only marks, mark size=2]
coordinates {
(  0.00000000,   0.00000000)
(  0.20000000,   0.00000000)
(  0.40000000,   0.00000000)
(  0.60000000,   0.00000000)
(  0.80000000,   0.00000000)
(  1.00000000,   0.00000000)};\label{line:mshparam_1dof_1d:other_nodes}

\addplot [mark options={solid}, mark=*, mark size=3, red, only marks, mark size=2]
coordinates {
(  0.40000000,   0.00000000)};\label{line:mshparam_1dof_1d:phi_node}

\addplot [solid, forget plot]
coordinates {
(  0.00000000,   0.25000000)
(  0.25000000,   0.25000000)
(  0.50000000,   0.25000000)
(  0.66666667,   0.25000000)
(  0.83333333,   0.25000000)
(  1.00000000,   0.25000000)};

\addplot [red, forget plot]
coordinates {
(  0.50000000,   0.25000000)};

\addplot [solid, mark options={solid}, mark=*, mark size=3, only marks, mark size=2]
coordinates {
(  0.00000000,   0.25000000)
(  0.25000000,   0.25000000)
(  0.50000000,   0.25000000)
(  0.66666667,   0.25000000)
(  0.83333333,   0.25000000)
(  1.00000000,   0.25000000)};\label{line:mshparam_1dof_1d:other_nodes}

\addplot [mark options={solid}, mark=*, mark size=3, red, only marks, mark size=2]
coordinates {
(  0.50000000,   0.25000000)};\label{line:mshparam_1dof_1d:phi_node}

\addplot [solid, forget plot]
coordinates {
(  0.00000000,   0.50000000)
(  0.17500000,   0.50000000)
(  0.35000000,   0.50000000)
(  0.56666667,   0.50000000)
(  0.78333333,   0.50000000)
(  1.00000000,   0.50000000)};

\addplot [red, forget plot]
coordinates {
(  0.35000000,   0.50000000)};

\addplot [solid, mark options={solid}, mark=*, mark size=3, only marks, mark size=2]
coordinates {
(  0.00000000,   0.50000000)
(  0.17500000,   0.50000000)
(  0.35000000,   0.50000000)
(  0.56666667,   0.50000000)
(  0.78333333,   0.50000000)
(  1.00000000,   0.50000000)};\label{line:mshparam_1dof_1d:other_nodes}

\addplot [mark options={solid}, mark=*, mark size=3, red, only marks, mark size=2]
coordinates {
(  0.35000000,   0.50000000)};\label{line:mshparam_1dof_1d:phi_node}

\addplot [solid, forget plot]
coordinates {
(  0.00000000,   0.75000000)
(  0.32500000,   0.75000000)
(  0.65000000,   0.75000000)
(  0.76666667,   0.75000000)
(  0.88333333,   0.75000000)
(  1.00000000,   0.75000000)};

\addplot [red, forget plot]
coordinates {
(  0.65000000,   0.75000000)};

\addplot [solid, mark options={solid}, mark=*, mark size=3, only marks, mark size=2]
coordinates {
(  0.00000000,   0.75000000)
(  0.32500000,   0.75000000)
(  0.65000000,   0.75000000)
(  0.76666667,   0.75000000)
(  0.88333333,   0.75000000)
(  1.00000000,   0.75000000)};\label{line:mshparam_1dof_1d:other_nodes}

\addplot [mark options={solid}, mark=*, mark size=3, red, only marks, mark size=2]
coordinates {
(  0.65000000,   0.75000000)};\label{line:mshparam_1dof_1d:phi_node}

\addplot [dashed, black]
coordinates {
(  0.50000000,   0.00000000)
(  0.50000000,   0.75000000)};\label{line:mshparam_1dof_1d:xshk}

\addplot [solid, gray]
coordinates {
(  0.40000000,   0.00000000)
(  0.40000000,   0.75000000)};\label{line:mshparam_1dof_1d:xloc}

\node[below]    at    (axis cs:0.0, -0.02) {$x_l$};
\node[below]    at    (axis cs:1.0, -0.02) {$x_r$};
\draw [->] plot coordinates {(axis cs:1.06, 0.0) (axis cs:1.14, 0.0)
};

\node[right]    at    (axis cs:1.14, 0.0) {$X$};
\draw [->] plot coordinates {(axis cs:1.06, 0.25) (axis cs:1.14, 0.25)
};

\node[right]    at    (axis cs:1.14, 0.25) {$x$};
\draw [->] plot coordinates {(axis cs:1.06, 0.5) (axis cs:1.14, 0.5)
};

\node[right]    at    (axis cs:1.14, 0.5) {$x$};
\draw [->] plot coordinates {(axis cs:1.06, 0.75) (axis cs:1.14, 0.75)
};

\node[right]    at    (axis cs:1.14, 0.75) {$x$};
\end{axis}
\end{tikzpicture}
 \caption{Single degree of freedom $1d$ mesh parametrization in
          (\ref{eqn:mshparam1d}) corresponding to a function with a
          discontinuity at (\ref{line:mshparam_1dof_1d:xshk}). The choice of
          $\bar{x}$ (\ref{line:mshparam_1dof_1d:xloc}) is chosen to align with
          a mesh node (\ref{line:mshparam_1dof_1d:phi_node}) near the position
          of the discontinuity and the parameter $\phibold$ reduces to the
          position of this node.
          The remaining nodes (\ref{line:mshparam_1dof_1d:other_nodes}) are
          determined by uniformly distributing them in the intervals
          $(x_l,\,\bar{x})$ and $(\bar{x},\,x_r)$. The first configuration
          (\emph{bottom}) shows the reference domain with a uniform nodal
          distribution and the remaining configurations corresponds to the
          deformation with
          $\phibold = x^*,\,x^*-0.75h,\,x^*+0.75h$, respectively.
          This range of $\phibold$ is the same range used to perform the
          sweep of the objective function (\ref{eqn:obj-redsp-l2proj}) in
          Section~\ref{sec:track-obj-behavior} and
          Figure~\ref{fig:l2proj1d_midsweep}.}
 \label{fig:mshparam-1dof-1d}
\end{figure}
}

With this single degree of freedom parametrization, we consider the graph of
$j_\alpha(\phibold)$ for various values of $\alpha$ in
Figure~\ref{fig:l2proj1d_midsweep}, where $\ubm_H$ corresponds to the
$L^2$ projection of the piecewise sinusoidal function in (\ref{eqn:l2proj1d})
with $k = 2$ onto a piecewise polynomial basis with $17$ elements and
polynomial orders $p \in \{1,\,2,\,3,\,4\}$. The reference nodes are
uniformly distributed and therefore $f_{msh}$ possesses a minimum at
$\phibold = \bar{x}$. For the case with no mesh
regularization, the proposed objective function satisfies the two desired
properties, i.e., possesses a local minima when the mesh is aligned with
the discontinuity and monotonically approaches this minima in a neighborhood
of radius $h/2$. As the penalty parameter is increased to
$\alpha = 1$, the shape of the graph changes in that values near
$\bar{x}$ decrease, but monotonicity is not broken. As $\alpha$ increases
to $10$, the monotonicity property breaks down and a local minima is
introduced at $\bar{x}$, which becomes stronger as $\alpha$ is increased
to $20$. Once $\alpha$ is sufficiently large ($\alpha = 1000$), the mesh
regularization term dominates causing the local minima at the shock location
to vanish. Similar trends were observed when using coarser and finer
discretization into piecewise polynomials. This simple test shows that the
discontinuity indicator proposed in (\ref{eqn:obj}) is well-suited as the
objective function in our discontinuity-tracking framework, provided the
penalty parameter is chosen to \emph{balance} the scale of the discontinuity
indicator and mesh regularization terms.

\ifbool{fastcompile}{}{
\begin{figure}
 \begin{tikzpicture}
\begin{groupplot}[
    group style={
        group size=2 by 3,
        horizontal sep=2.4cm
    },
    xmin=0.45, xmax=0.55,
    width=8.0cm,
    height=5.0cm,
]

\nextgroupplot[ymin=0, ymax=3.0,
               ylabel={$j_\alpha(\phibold),\,\alpha=0.0$}]
\addplot [black, solid, thick, mark=*, mark size=2, mark options={solid}, mark repeat={8}]  table[x index=0, y index=1] {dat/l2proj1d_midnode_nel0017_p1.dev_from_mean.dev_from_unifmsh.sweep.dat}; \label{line:l2proj1d_midsweep:p1}
\addplot [blue, solid, thick, mark=square*, mark size=2, mark options={solid}, mark repeat={8}]  table[x index=0, y index=1] {dat/l2proj1d_midnode_nel0017_p2.dev_from_mean.dev_from_unifmsh.sweep.dat}; \label{line:l2proj1d_midsweep:p2}
\addplot [red, solid, thick, mark=triangle*, mark size=2, mark options={solid}, mark repeat={8}]  table[x index=0, y index=1] {dat/l2proj1d_midnode_nel0017_p3.dev_from_mean.dev_from_unifmsh.sweep.dat}; \label{line:l2proj1d_midsweep:p3}
\addplot [green, solid, thick, mark=diamond*, mark size=2, mark options={solid}, mark repeat={8}]  table[x index=0, y index=1] {dat/l2proj1d_midnode_nel0017_p4.dev_from_mean.dev_from_unifmsh.sweep.dat}; \label{line:l2proj1d_midsweep:p4}
\addplot [draw=none, fill=red, fill opacity=0.1]
coordinates {
(  0.470588235,  0)
(  0.470588235,  3.0)
(  0.529411765,  3.0)
(  0.529411765,  0)
(  0.470588235,  0)}; \label{line:l2proj1d_midsweep:h_over_2}
\addplot [gray, solid, thin]
coordinates {
(0.47058823529411764, -2e-3)
(0.47058823529411764, 3.0)}; \label{line:l2proj1d_midsweep:xloc}
\addplot [black, dashed, thin]
coordinates {
(  0.5,  -2e-3)
(  0.5,  3.0)}; \label{line:l2proj1d_midsweep:xshk}

\nextgroupplot[ymin=0, ymax=3.0,
               ylabel={$j_\alpha(\phibold),\,\alpha=1$}]
\addplot [black, solid, thick, mark=*, mark size=2, mark options={solid}, mark repeat={8}]  table[x index=0, y index=7] {dat/l2proj1d_midnode_nel0017_p1.dev_from_mean.dev_from_unifmsh.sweep.dat};
\addplot [blue, solid, thick, mark=square*, mark size=2, mark options={solid}, mark repeat={8}]  table[x index=0, y index=7] {dat/l2proj1d_midnode_nel0017_p2.dev_from_mean.dev_from_unifmsh.sweep.dat};
\addplot [red, solid, thick, mark=triangle*, mark size=2, mark options={solid}, mark repeat={8}]  table[x index=0, y index=7] {dat/l2proj1d_midnode_nel0017_p3.dev_from_mean.dev_from_unifmsh.sweep.dat};
\addplot [green, solid, thick, mark=diamond*, mark size=2, mark options={solid}, mark repeat={8}]  table[x index=0, y index=7] {dat/l2proj1d_midnode_nel0017_p4.dev_from_mean.dev_from_unifmsh.sweep.dat};
\addplot [draw=none, fill=red, fill opacity=0.1]
coordinates {
(  0.470588235,  0)
(  0.470588235,  3.0)
(  0.529411765,  3.0)
(  0.529411765,  0)
(  0.470588235,  0)};
\addplot [gray, solid, thin]
coordinates {
(0.47058823529411764, -2e-3)
(0.47058823529411764, 3.0)};
\addplot [black, dashed, thin]
coordinates {
(  0.5,  -2e-3)
(  0.5,  3.0)};

\nextgroupplot[ymin=0, ymax=4.1,
               ylabel={$j_\alpha(\phibold),\,\alpha=10$}]
\addplot [black, solid, thick, mark=*, mark size=2, mark options={solid}, mark repeat={8}]  table[x index=0, y index=8] {dat/l2proj1d_midnode_nel0017_p1.dev_from_mean.dev_from_unifmsh.sweep.dat};
\addplot [blue, solid, thick, mark=square*, mark size=2, mark options={solid}, mark repeat={8}]  table[x index=0, y index=8] {dat/l2proj1d_midnode_nel0017_p2.dev_from_mean.dev_from_unifmsh.sweep.dat};
\addplot [red, solid, thick, mark=triangle*, mark size=2, mark options={solid}, mark repeat={8}]  table[x index=0, y index=8] {dat/l2proj1d_midnode_nel0017_p3.dev_from_mean.dev_from_unifmsh.sweep.dat};
\addplot [green, solid, thick, mark=diamond*, mark size=2, mark options={solid}, mark repeat={8}]  table[x index=0, y index=8] {dat/l2proj1d_midnode_nel0017_p4.dev_from_mean.dev_from_unifmsh.sweep.dat};
\addplot [draw=none, fill=red, fill opacity=0.1]
coordinates {
(  0.470588235,  0)
(  0.470588235,  4.1)
(  0.529411765,  4.1)
(  0.529411765,  0)
(  0.470588235,  0)};
\addplot [gray, solid, thin]
coordinates {
(0.47058823529411764, -2e-3)
(0.47058823529411764, 4.1)};
\addplot [black, dashed, thin]
coordinates {
(  0.5,  -2e-3)
(  0.5,  4.1)};

\nextgroupplot[ymin=0, ymax=6.0,
               ylabel={$j_\alpha(\phibold),\,\alpha=20$}]
\addplot [black, solid, thick, mark=*, mark size=2, mark options={solid}, mark repeat={8}]  table[x index=0, y index=9] {dat/l2proj1d_midnode_nel0017_p1.dev_from_mean.dev_from_unifmsh.sweep.dat};
\addplot [blue, solid, thick, mark=square*, mark size=2, mark options={solid}, mark repeat={8}]  table[x index=0, y index=9] {dat/l2proj1d_midnode_nel0017_p2.dev_from_mean.dev_from_unifmsh.sweep.dat};
\addplot [red, solid, thick, mark=triangle*, mark size=2, mark options={solid}, mark repeat={8}]  table[x index=0, y index=9] {dat/l2proj1d_midnode_nel0017_p3.dev_from_mean.dev_from_unifmsh.sweep.dat};
\addplot [green, solid, thick, mark=diamond*, mark size=2, mark options={solid}, mark repeat={8}]  table[x index=0, y index=9] {dat/l2proj1d_midnode_nel0017_p4.dev_from_mean.dev_from_unifmsh.sweep.dat};
\addplot [draw=none, fill=red, fill opacity=0.1]
coordinates {
(  0.470588235,  0)
(  0.470588235,  6.0)
(  0.529411765,  6.0)
(  0.529411765,  0)
(  0.470588235,  0)};
\addplot [gray, solid, thin]
coordinates {
(0.47058823529411764, -2e-3)
(0.47058823529411764, 6.0)};
\addplot [black, dashed, thin]
coordinates {
(  0.5,  -2e-3)
(  0.5,  6.0)};

\nextgroupplot[ymin=0, ymax=9.0,
               ylabel={$j_\alpha(\phibold),\,\alpha=40$},
               xlabel={$\phibold$ (position of node closest to shock)}]
\addplot [black, solid, thick, mark=*, mark size=2, mark options={solid}, mark repeat={8}]  table[x index=0, y index=10] {dat/l2proj1d_midnode_nel0017_p1.dev_from_mean.dev_from_unifmsh.sweep.dat};
\addplot [blue, solid, thick, mark=square*, mark size=2, mark options={solid}, mark repeat={8}]  table[x index=0, y index=10] {dat/l2proj1d_midnode_nel0017_p2.dev_from_mean.dev_from_unifmsh.sweep.dat};
\addplot [red, solid, thick, mark=triangle*, mark size=2, mark options={solid}, mark repeat={8}]  table[x index=0, y index=10] {dat/l2proj1d_midnode_nel0017_p3.dev_from_mean.dev_from_unifmsh.sweep.dat};
\addplot [green, solid, thick, mark=diamond*, mark size=2, mark options={solid}, mark repeat={8}]  table[x index=0, y index=10] {dat/l2proj1d_midnode_nel0017_p4.dev_from_mean.dev_from_unifmsh.sweep.dat};
\addplot [draw=none, fill=red, fill opacity=0.1]
coordinates {
(  0.470588235,  0)
(  0.470588235,  9.0)
(  0.529411765,  9.0)
(  0.529411765,  0)
(  0.470588235,  0)};
\addplot [gray, solid, thin]
coordinates {
(0.47058823529411764, -2e-3)
(0.47058823529411764, 9.0)};
\addplot [black, dashed, thin]
coordinates {
(  0.5,  -2e-3)
(  0.5,  9.0)};

\nextgroupplot[ymin=0, ymax=150.0,
               ylabel={$j_\alpha(\phibold),\,\alpha=1000$},
               xlabel={$\phibold$ (position of node closest to shock)}]
\addplot [black, solid, thick, mark=*, mark size=2, mark options={solid}, mark repeat={8}]  table[x index=0, y index=11] {dat/l2proj1d_midnode_nel0017_p1.dev_from_mean.dev_from_unifmsh.sweep.dat};
\addplot [blue, solid, thick, mark=square*, mark size=2, mark options={solid}, mark repeat={8}]  table[x index=0, y index=11] {dat/l2proj1d_midnode_nel0017_p2.dev_from_mean.dev_from_unifmsh.sweep.dat};
\addplot [red, solid, thick, mark=triangle*, mark size=2, mark options={solid}, mark repeat={8}]  table[x index=0, y index=11] {dat/l2proj1d_midnode_nel0017_p3.dev_from_mean.dev_from_unifmsh.sweep.dat};
\addplot [green, solid, thick, mark=diamond*, mark size=2, mark options={solid}, mark repeat={8}]  table[x index=0, y index=11] {dat/l2proj1d_midnode_nel0017_p4.dev_from_mean.dev_from_unifmsh.sweep.dat};
\addplot [draw=none, fill=red, fill opacity=0.1]
coordinates {
(  0.470588235,  0)
(  0.470588235,  150)
(  0.529411765,  150)
(  0.529411765,  0)
(  0.470588235,  0)};
\addplot [gray, solid, thin]
coordinates {
(0.47058823529411764, -2e-3)
(0.47058823529411764, 150)};
\addplot [black, dashed, thin]
coordinates {
(  0.5,  -2e-3)
(  0.5,  150)};

\end{groupplot}
\end{tikzpicture}
 \caption{The objective function $j_\alpha(\phibold)$ corresponding to the
          $L^2$ projection of the discontinuous function in
          (\ref{eqn:l2proj1d}) with $k = 2$ onto a piecewise polynomial
          basis of degree
          $p = 1$ (\ref{line:l2proj1d_midsweep:p1}),
          $p = 2$ (\ref{line:l2proj1d_midsweep:p2}),
          $p = 3$ (\ref{line:l2proj1d_midsweep:p3}),
          $p = 4$ (\ref{line:l2proj1d_midsweep:p4})
          with $|\Ecal_{h,p}| = 17$ elements for various values of the mesh
          regularization parameter $\alpha$ under the mesh parameterization
          in (\ref{eqn:mshparam1d}) and Figure~\ref{fig:mshparam-1dof-1d}.
          The parameter $\phibold$ is swept over the range
          $[x^*-0.75h,\,x^*+0.75h]$, where $x^*$ is the location
          of the discontinuity (\ref{line:l2proj1d_midsweep:xshk}). The position
          of $\bar{x}$ in (\ref{eqn:mshparam1d}) is indicated with
          (\ref{line:l2proj1d_midsweep:xloc}); when $\phibold = \bar{x}$
          nodes are uniformly distributed and $f_{msh}$ takes its minimum
          values. The red shaded region identifies the neighborhood of radius
          $h/2$ about $x^*$. The objective function with no mesh regularization
          ($\alpha = 0$) in the \emph{top left} plot shows the two desired
          properties discussed in Section~\ref{sec:track-obj}, namely, a
          minimum at $x^*$ and monotonically decreases to this minimum in a
          neighborhood centered at $x^*$ with radius $h/2$. The
          remaining plots show increasing values of the regularization
          parameter $\alpha$. The minimum does not noticeably move until
          $\alpha$ is sufficiently large ($\alpha = 40$) at which point is
          moves from $x^*$ to $\bar{x}$, i.e., the objective function
          transitions from prioritizing discontinuity tracking to mesh
          regularization.}
 \label{fig:l2proj1d_midsweep}
\end{figure}
}

Given the promising performance of the indicator, we aim to study
it in the situation where the solution oscillates within an
element and the basis may be underresolved. To this end, we
repeat the previous numerical experiment without mesh regularization
on a piecewise polynomial basis with two elements and polynomial
orders $p \in \{1,\,2,\,3,\,4\}$ for increasingly oscillatory
discontinuous functions: (\ref{eqn:l2proj1d}) with $k = 1$,
$k = 3$, $k = 5$, $k = 7$; see Figure~\ref{fig:l2proj1d_midsweep_nperiod}.
As expected, for all functions and polynomial orders, the
objective function obtains its minima when the mesh is
aligned with the shock. As the frequency of the oscillations
in the discontinuous function increases, the graphs of
the objective function become more oscillatory \emph{yet
remain mostly monotonic within a radius of $h/2$ of the shock
location}. Monotonicity does break down for the most extreme case
($k = 7$) for the lowest polynomial orders $p = 1$ and $p = 2$, which
is an artifact of an extremely underresolved basis.

\ifbool{fastcompile}{}{
\begin{figure}
 \begin{tikzpicture}
\begin{groupplot}[
    group style={
        group size=2 by 3,
        horizontal sep=1.8cm
    },
    xmin=0.0, xmax=1.0,
    width=8.0cm,
    height=5.0cm,
    tick label style={/pgf/number format/fixed}
]

\nextgroupplot[ymin=0, ymax=0.15,
               ylabel={$j_\alpha(\phibold),\,\alpha=0.0$}]
\addplot [black, solid, thick, mark=*, mark size=2, mark options={solid}, mark repeat={8}]  table[x index=0, y index=1] {dat/l2proj1d_midnode_nel0002_p1.dev_from_mean.dev_from_unifmsh.sweep_nperiod1.dat};
\addplot [blue, solid, thick, mark=square*, mark size=2, mark options={solid}, mark repeat={8}]  table[x index=0, y index=1] {dat/l2proj1d_midnode_nel0002_p2.dev_from_mean.dev_from_unifmsh.sweep_nperiod1.dat};
\addplot [red, solid, thick, mark=triangle*, mark size=2, mark options={solid}, mark repeat={8}]  table[x index=0, y index=1] {dat/l2proj1d_midnode_nel0002_p3.dev_from_mean.dev_from_unifmsh.sweep_nperiod1.dat};
\addplot [green, solid, thick, mark=diamond*, mark size=2, mark options={solid}, mark repeat={8}]  table[x index=0, y index=1] {dat/l2proj1d_midnode_nel0002_p4.dev_from_mean.dev_from_unifmsh.sweep_nperiod1.dat};
\addplot [draw=none, fill=red, fill opacity=0.1]
coordinates {
(  0.25,  0)
(  0.25,  0.15)
(  0.75,  0.15)
(  0.75,  0)
(  0.25,  0)};
\addplot [black, dashed, thin]
coordinates {
(  0.5,  0.0)
(  0.5,  0.15)};

\nextgroupplot[ymin=0, ymax=0.4]
\addplot [black, solid, thick, mark=*, mark size=2, mark options={solid}, mark repeat={8}]  table[x index=0, y index=1] {dat/l2proj1d_midnode_nel0002_p1.dev_from_mean.dev_from_unifmsh.sweep_nperiod3.dat};
\addplot [blue, solid, thick, mark=square*, mark size=2, mark options={solid}, mark repeat={8}]  table[x index=0, y index=1] {dat/l2proj1d_midnode_nel0002_p2.dev_from_mean.dev_from_unifmsh.sweep_nperiod3.dat};
\addplot [red, solid, thick, mark=triangle*, mark size=2, mark options={solid}, mark repeat={8}]  table[x index=0, y index=1] {dat/l2proj1d_midnode_nel0002_p3.dev_from_mean.dev_from_unifmsh.sweep_nperiod3.dat};
\addplot [green, solid, thick, mark=diamond*, mark size=2, mark options={solid}, mark repeat={8}]  table[x index=0, y index=1] {dat/l2proj1d_midnode_nel0002_p4.dev_from_mean.dev_from_unifmsh.sweep_nperiod3.dat};
\addplot [draw=none, fill=red, fill opacity=0.1]
coordinates {
(  0.25,  0)
(  0.25,  0.4)
(  0.75,  0.4)
(  0.75,  0)
(  0.25,  0)}; 
\addplot [black, dashed, thin]
coordinates {
(  0.5,  0.0)
(  0.5,  0.4)}; 

\nextgroupplot[ymin=0, ymax=0.4,
               xlabel={$\phibold$ (position of node closest to shock)},
               ylabel={$j_\alpha(\phibold),\,\alpha=0.0$}]
\addplot [black, solid, thick, mark=*, mark size=2, mark options={solid}, mark repeat={8}]  table[x index=0, y index=1] {dat/l2proj1d_midnode_nel0002_p1.dev_from_mean.dev_from_unifmsh.sweep_nperiod5.dat};
\addplot [blue, solid, thick, mark=square*, mark size=2, mark options={solid}, mark repeat={8}]  table[x index=0, y index=1] {dat/l2proj1d_midnode_nel0002_p2.dev_from_mean.dev_from_unifmsh.sweep_nperiod5.dat};
\addplot [red, solid, thick, mark=triangle*, mark size=2, mark options={solid}, mark repeat={8}]  table[x index=0, y index=1] {dat/l2proj1d_midnode_nel0002_p3.dev_from_mean.dev_from_unifmsh.sweep_nperiod5.dat};
\addplot [green, solid, thick, mark=diamond*, mark size=2, mark options={solid}, mark repeat={8}]  table[x index=0, y index=1] {dat/l2proj1d_midnode_nel0002_p4.dev_from_mean.dev_from_unifmsh.sweep_nperiod5.dat};
\addplot [draw=none, fill=red, fill opacity=0.1]
coordinates {
(  0.25,  0)
(  0.25,  0.4)
(  0.75,  0.4)
(  0.75,  0)
(  0.25,  0)}; 
\addplot [black, dashed, thin]
coordinates {
(  0.5,  0.0)
(  0.5,  0.4)}; 

\nextgroupplot[ymin=0, ymax=0.6,
               xlabel={$\phibold$ (position of node closest to shock)}]
\addplot [black, solid, thick, mark=*, mark size=2, mark options={solid}, mark repeat={8}]  table[x index=0, y index=1] {dat/l2proj1d_midnode_nel0002_p1.dev_from_mean.dev_from_unifmsh.sweep_nperiod7.dat};
\addplot [blue, solid, thick, mark=square*, mark size=2, mark options={solid}, mark repeat={8}]  table[x index=0, y index=1] {dat/l2proj1d_midnode_nel0002_p2.dev_from_mean.dev_from_unifmsh.sweep_nperiod7.dat};
\addplot [red, solid, thick, mark=triangle*, mark size=2, mark options={solid}, mark repeat={8}]  table[x index=0, y index=1] {dat/l2proj1d_midnode_nel0002_p3.dev_from_mean.dev_from_unifmsh.sweep_nperiod7.dat};
\addplot [green, solid, thick, mark=diamond*, mark size=2, mark options={solid}, mark repeat={8}]  table[x index=0, y index=1] {dat/l2proj1d_midnode_nel0002_p4.dev_from_mean.dev_from_unifmsh.sweep_nperiod7.dat};
\addplot [draw=none, fill=red, fill opacity=0.1]
coordinates {
(  0.25,  0)
(  0.25,  0.6)
(  0.75,  0.6)
(  0.75,  0)
(  0.25,  0)}; 
\addplot [black, dashed, thin]
coordinates {
(  0.5,  0.0)
(  0.5,  0.6)}; 

\end{groupplot}
\end{tikzpicture}
 \caption{The objective function $j_\alpha(\phibold)$ with no mesh
          regularization ($\alpha = 0.0$) corresponding to the
          $L^2$ projection of the periodic discontinuous function in
          (\ref{eqn:l2proj1d}) of varying frequency
          (\emph{top left}: $k = 1$, \emph{top right}: $k = 3$,
           \emph{bottom left}: $k = 5$, \emph{bottom right}: $k = 7$)
          onto a piecewise polynomial basis of degree
          $p = 1$ (\ref{line:l2proj1d_midsweep:p1}),
          $p = 2$ (\ref{line:l2proj1d_midsweep:p2}),
          $p = 3$ (\ref{line:l2proj1d_midsweep:p3}),
          $p = 4$ (\ref{line:l2proj1d_midsweep:p4})
          with $|\Ecal_{h,p}| = 2$ elements under the mesh parameterization
          in (\ref{eqn:mshparam1d}) and Figure~\ref{fig:mshparam-1dof-1d}.
          The parameter $\phibold$ is swept over the range
          $[x^*-0.75h,\,x^*+0.75h]$, where $x^*$ is the location
          of the discontinuity (\ref{line:l2proj1d_midsweep:xshk}).
          The red shaded region identifies the neighborhood of radius
          $h/2$ about $x^*$.}
 \label{fig:l2proj1d_midsweep_nperiod}
\end{figure}
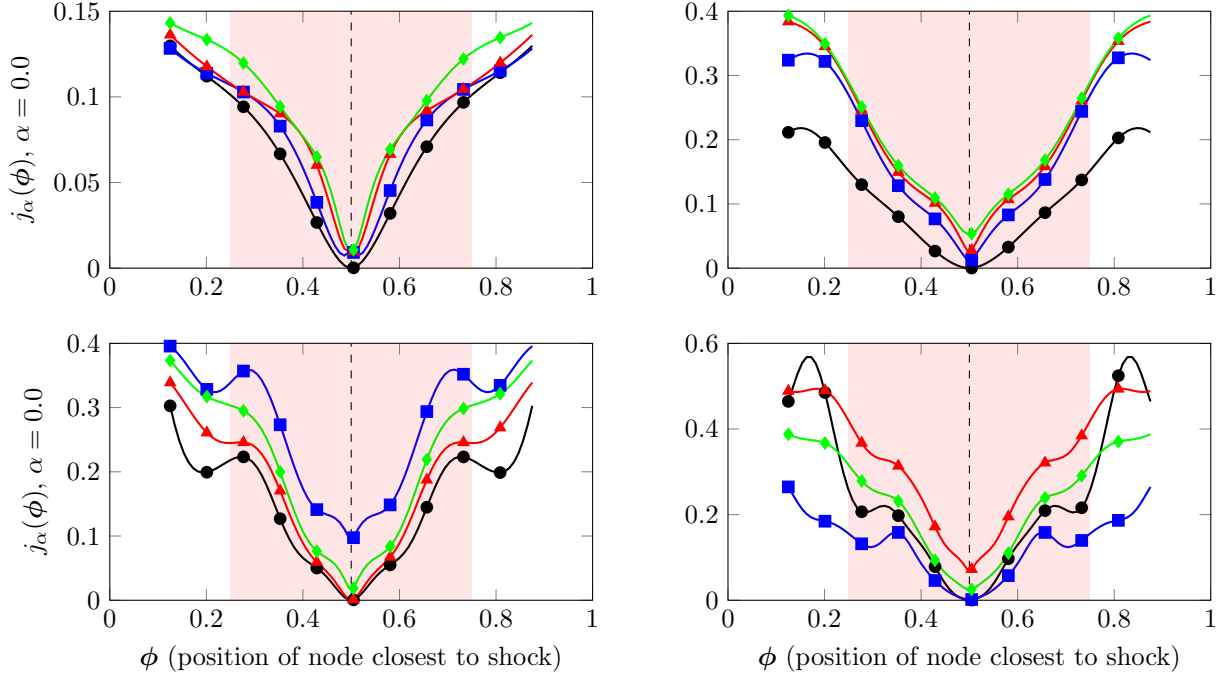
}

Since the mesh distortion metric is defined separately for $d = 1$ and $d > 1$,
it is important to repeat the investigation into $j_\alpha(\phibold)$ in
$d = 2$ dimensions. The same mesh parametrization is applied to the
\emph{first} coordinate in the two-dimensional setting with the second
coordinate held fixed
\begin{equation} \label{eqn:mshparam2d-disc}
\begin{aligned}
\xbm_i(\phibold) \cdot \ebm_1 &=
 \begin{cases}
   \displaystyle{x_l + (\Xbm_i\cdot\ebm_1-x_l)\frac{\phibold-x_l}{\bar{x}-x_l}}
                                                  \quad & \Xbm_i \leq \bar{x} \\
   \displaystyle{x_r - (x_r-\Xbm_i\cdot\ebm_1)\frac{x_r-\phibold}{x_r-\bar{x}}}
                                                  \quad & \Xbm_i > \bar{x}
 \end{cases} \\
 \xbm_i(\phibold) \cdot \ebm_2 &= \Xbm_i \cdot \ebm_2
\end{aligned}
\end{equation}
A depiction of the single degree of freedom mesh parametrization in $d=2$
dimensions is provided in Figure~\ref{fig:mshparam-1dof-2d}. This figure also
shows the range that $\phibold$ will be swept to study the behavior of
$j_\alpha(\phibold)$ next.

\ifbool{fastcompile}{}{
\begin{figure}
 \centering
 \begin{tikzpicture}

\begin{groupplot}[
    group style={ 
        group size=2 by 2,
        horizontal sep=1.3cm
    },
    enlargelimits=false,
    axis on top, axis equal image,
    xmin=0, xmax=1, ymin=0, ymax=1,
    width=0.4\textwidth,
]

\nextgroupplot
\addplot graphics [xmin=0, xmax=1, ymin=0, ymax=1] {{img/l2proj2d_midnode.param1d.config0}.png};
\addplot [solid, red, mark options={solid}, mark=*, mark size=0.75, only marks]
coordinates {
(  0.4,   0.0)
(  0.4,   0.1)
(  0.4,   0.2)
(  0.4,   0.3)
(  0.4,   0.4)
(  0.4,   0.5)
(  0.4,   0.6)
(  0.4,   0.7)
(  0.4,   0.8)
(  0.4,   0.9)
(  0.4,   1.0)}; \label{line:mshparam_1dof_2d:phi_node}
\addplot [dashed, black, thin]
coordinates {
(0.5, 0.0)
(0.5, 1.0)}; \label{line:mshparam_1dof_2d:xshk}
\addplot [solid, gray, thin]
coordinates {
(0.4, 0.0)
(0.4, 1.0)}; \label{line:mshparam_1dof_2d:xloc}

\nextgroupplot
\addplot graphics [xmin=0, xmax=1, ymin=0, ymax=1] {{img/l2proj2d_midnode.param1d.config1}.png};
\addplot [solid, red, mark options={solid}, mark=*, mark size=0.75, only marks]
coordinates {
(  0.5,   0.0)
(  0.5,   0.1)
(  0.5,   0.2)
(  0.5,   0.3)
(  0.5,   0.4)
(  0.5,   0.5)
(  0.5,   0.6)
(  0.5,   0.7)
(  0.5,   0.8)
(  0.5,   0.9)
(  0.5,   1.0)};
\addplot [dashed, black, thin]
coordinates {
(0.5, 0.0)
(0.5, 1.0)};
\addplot [solid, gray, thin]
coordinates {
(0.4, 0.0)
(0.4, 1.0)};

\nextgroupplot
\addplot graphics [xmin=0, xmax=1, ymin=0, ymax=1] {{img/l2proj2d_midnode.param1d.config2}.png};
\addplot [solid, red, mark options={solid}, mark=*, mark size=0.75, only marks]
coordinates {
(  0.35,   0.0)
(  0.35,   0.1)
(  0.35,   0.2)
(  0.35,   0.3)
(  0.35,   0.4)
(  0.35,   0.5)
(  0.35,   0.6)
(  0.35,   0.7)
(  0.35,   0.8)
(  0.35,   0.9)
(  0.35,   1.0)};
\addplot [dashed, black, thin]
coordinates {
(0.5, 0.0)
(0.5, 1.0)};
\addplot [solid, gray, thin]
coordinates {
(0.4, 0.0)
(0.4, 1.0)};

\nextgroupplot
\addplot graphics [xmin=0, xmax=1, ymin=0, ymax=1] {{img/l2proj2d_midnode.param1d.config3}.png};
\addplot [solid, red, mark options={solid}, mark=*, mark size=0.75, only marks]
coordinates {
(  0.65,   0.0)
(  0.65,   0.1)
(  0.65,   0.2)
(  0.65,   0.3)
(  0.65,   0.4)
(  0.65,   0.5)
(  0.65,   0.6)
(  0.65,   0.7)
(  0.65,   0.8)
(  0.65,   0.9)
(  0.65,   1.0)};
\addplot [dashed, black, thin]
coordinates {
(0.5, 0.0)
(0.5, 1.0)};
\addplot [solid, gray, thin]
coordinates {
(0.4, 0.0)
(0.4, 1.0)};
\end{groupplot}

\end{tikzpicture}
 \caption{Single degree of freedom $2d$ mesh parametrization in
          (\ref{eqn:mshparam1d}) and problem setup for investigation into
          the proposed objective function in (\ref{eqn:obj-redsp-l2proj})
          as applied to the $L^2$ projection of the discontinuous function
          in (\ref{eqn:l2proj2d0}). The choice of $\bar{x}$
          (\ref{line:mshparam_1dof_2d:xloc}) is chosen to align with a line of
          mesh nodes (\ref{line:mshparam_1dof_2d:phi_node}) near the position
          of the discontinuity (\ref{line:mshparam_1dof_2d:xshk}) and the
          parameter $\phibold$ reduces to the position of this node.
          The remaining nodes ($x$-coordinates only; $y$-coordinates are frozen)
          are determined by uniformly distributing them in the intervals
          $(x_l,\,\bar{x})$ and $(\bar{x},\,x_r)$. The first configuration
          (\emph{top left}) shows the reference domain with a uniform nodal
          distribution and the remaining configurations corresponds to the
          deformation with
          $\phibold = x^*,\,x^*-0.75h,\,x^*+0.75h$.
          This range of $\phibold$ is the same range used to perform the
          sweep of the objective function (\ref{eqn:obj-redsp-l2proj}) in
          Section~\ref{sec:track-obj-behavior} and
          Figure~\ref{fig:l2proj2d_midsweep}.}
 \label{fig:mshparam-1dof-2d}
\end{figure}
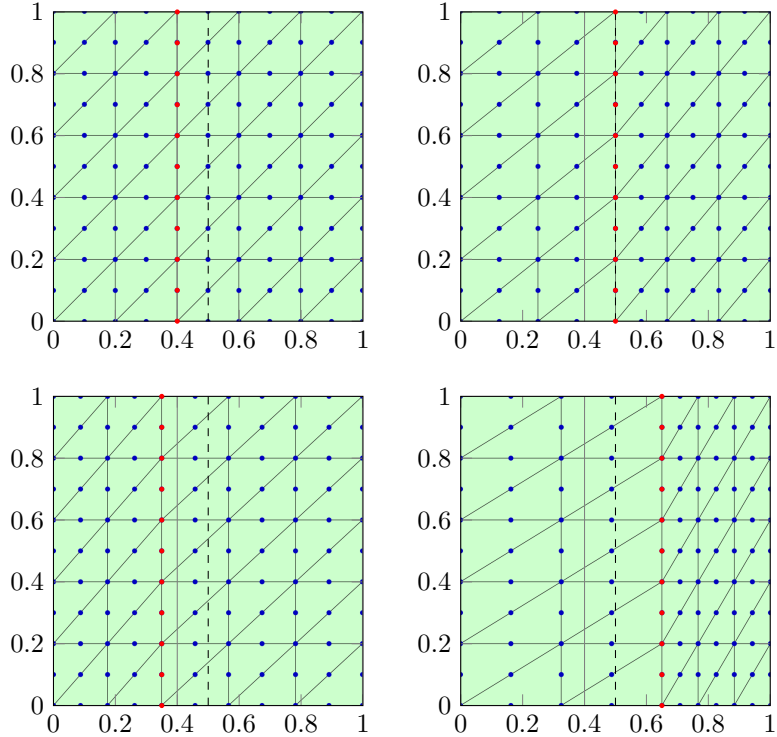
}

Figure~\ref{fig:l2proj2d_midsweep} shows the graph of $j_\alpha(\phibold)$
for various values of $\alpha$ corresponding to this single degree of
freedom parametrization,  where $\ubm_H$ corresponds to the
$L^2$ projection of the piecewise sinusoidal function
\begin{equation} \label{eqn:l2proj2d0}
 \eta(x,\,y) =
 \begin{cases}
  0.25\sin(4\pi x) - 0.5, \qquad &x<0 \\
  0.25\sin(4\pi x) + 0.5, \qquad &x>0,
 \end{cases}
\end{equation}
$(x,\,y) \in \Omega = (-1,\,1) \times (-1,\,1)$, onto a piecewise
polynomial basis with $50$ elements and polynomial orders
$p \in \{1,\,2,\,3\}$. Similar to the one-dimensional case, numerical
studies indicate that the proposed objective satisfies the desired properties
for a broad range of problem configurations,
i.e., possesses a local minima when the mesh is aligned with
the discontinuity and monotonically approaches this minima in a neighborhood
of radius about $h/2$. Also, as the penalty parameter is increased the mesh
regularization is preferred and the local minima at the discontinuity
disappears. Similar trends were observed when using coarser and finer
discretization into piecewise polynomials.

\ifbool{fastcompile}{}{
\begin{figure}
 \begin{tikzpicture}
\begin{groupplot}[
    group style={
        group name=conv plots,
        group size=2 by 3,
        horizontal sep=2.4cm
    },
    xmin=0.32, xmax=0.68,
    width=8.0cm,
    height=5.0cm,
]

\nextgroupplot[ymin=0, ymax=12,
               ylabel={$j_\alpha(\phibold),\,\alpha=0.0$}]
\addplot [black, solid, thick, mark=*, mark size=2, mark options={solid}, mark repeat={8}]  table[x index=0, y index=1] {dat/l2proj2d_midnode_nel0050_p1.dev_from_mean.distortion_nd_intg_phys.sweep.dat}; \label{line:l2proj2d_midsweep:p1}
\addplot [blue, solid, thick, mark=square*, mark size=2, mark options={solid}, mark repeat={8}]  table[x index=0, y index=1] {dat/l2proj2d_midnode_nel0050_p2.dev_from_mean.distortion_nd_intg_phys.sweep.dat}; \label{line:l2proj2d_midsweep:p2}
\addplot [red, solid, thick, mark=triangle*, mark size=2, mark options={solid}, mark repeat={8}]  table[x index=0, y index=1] {dat/l2proj2d_midnode_nel0050_p3.dev_from_mean.distortion_nd_intg_phys.sweep.dat}; \label{line:l2proj2d_midsweep:p3}
\addplot [draw=none, fill=red, fill opacity=0.1]
coordinates {
(  0.4,  0)
(  0.4,  12)
(  0.6,  12)
(  0.6,  0)
(  0.4,  0)}; \label{line:l2proj2d_midsweep:h_over_2}
\addplot [gray, solid, thin]
coordinates {
(0.4, -2e-3)
(0.4, 12)}; \label{line:l2proj2d_midsweep:xloc}
\addplot [black, dashed, thin]
coordinates {
(  0.5,  -2e-3)
(  0.5,  12)}; \label{line:l2proj2d_midsweep:xshk}

\nextgroupplot[ymin=0, ymax=15,
               ylabel={$j_\alpha(\phibold),\,\alpha=1$}]
\addplot [black, solid, thick, mark=*, mark size=2, mark options={solid}, mark repeat={8}]  table[x index=0, y index=8] {dat/l2proj2d_midnode_nel0050_p1.dev_from_mean.distortion_nd_intg_phys.sweep.dat};
\addplot [blue, solid, thick, mark=square*, mark size=2, mark options={solid}, mark repeat={8}]  table[x index=0, y index=8] {dat/l2proj2d_midnode_nel0050_p2.dev_from_mean.distortion_nd_intg_phys.sweep.dat};
\addplot [red, solid, thick, mark=triangle*, mark size=2, mark options={solid}, mark repeat={8}]  table[x index=0, y index=8] {dat/l2proj2d_midnode_nel0050_p3.dev_from_mean.distortion_nd_intg_phys.sweep.dat};
\addplot [draw=none, fill=red, fill opacity=0.1]
coordinates {
(  0.4,  0)
(  0.4,  15)
(  0.6,  15)
(  0.6,  0)
(  0.4,  0)};
\addplot [gray, solid, thin]
coordinates {
(0.4, -2e-3)
(0.4, 15)};
\addplot [black, dashed, thin]
coordinates {
(  0.5,  -2e-3)
(  0.5,  15)};

\nextgroupplot[ymin=20, ymax=35,
               ylabel={$j_\alpha(\phibold),\,\alpha=10$}]
\addplot [black, solid, thick, mark=*, mark size=2, mark options={solid}, mark repeat={8}]  table[x index=0, y index=9] {dat/l2proj2d_midnode_nel0050_p1.dev_from_mean.distortion_nd_intg_phys.sweep.dat};
\addplot [blue, solid, thick, mark=square*, mark size=2, mark options={solid}, mark repeat={8}]  table[x index=0, y index=9] {dat/l2proj2d_midnode_nel0050_p2.dev_from_mean.distortion_nd_intg_phys.sweep.dat};
\addplot [red, solid, thick, mark=triangle*, mark size=2, mark options={solid}, mark repeat={8}]  table[x index=0, y index=9] {dat/l2proj2d_midnode_nel0050_p3.dev_from_mean.distortion_nd_intg_phys.sweep.dat};
\addplot [draw=none, fill=red, fill opacity=0.1]
coordinates {
(  0.4,  20)
(  0.4,  35)
(  0.6,  35)
(  0.6,  20)
(  0.4,  20)};
\addplot [gray, solid, thin]
coordinates {
(0.4, -2e-3)
(0.4, 35)};
\addplot [black, dashed, thin]
coordinates {
(  0.5,  -2e-3)
(  0.5,  35)};

\nextgroupplot[ymin=2000, ymax=2400,
               ylabel={$j_\alpha(\phibold),\,\alpha=1000$}]
\addplot [black, solid, thick, mark=*, mark size=2, mark options={solid}, mark repeat={8}]  table[x index=0, y index=13] {dat/l2proj2d_midnode_nel0050_p1.dev_from_mean.distortion_nd_intg_phys.sweep.dat};
\addplot [blue, solid, thick, mark=square*, mark size=2, mark options={solid}, mark repeat={8}]  table[x index=0, y index=13] {dat/l2proj2d_midnode_nel0050_p2.dev_from_mean.distortion_nd_intg_phys.sweep.dat};
\addplot [red, solid, thick, mark=triangle*, mark size=2, mark options={solid}, mark repeat={8}]  table[x index=0, y index=13] {dat/l2proj2d_midnode_nel0050_p3.dev_from_mean.distortion_nd_intg_phys.sweep.dat};
\addplot [draw=none, fill=red, fill opacity=0.1]
coordinates {
(  0.4,  2000)
(  0.4,  2400)
(  0.6,  2400)
(  0.6,  2000)
(  0.4,  2000)};
\addplot [gray, solid, thin]
coordinates {
(0.4, 2000)
(0.4, 2400)};
\addplot [black, dashed, thin]
coordinates {
(  0.5,  2000)
(  0.5,  2400)};

\end{groupplot}
\end{tikzpicture}
 \caption{The objective function $j_\alpha(\phibold)$ corresponding to the
          $L^2$ projection of the discontinuous function in
          (\ref{eqn:l2proj2d0}) onto a piecewise polynomial basis of degree
          $p = 1$ (\ref{line:l2proj2d_midsweep:p1}),
          $p = 2$ (\ref{line:l2proj2d_midsweep:p2}),
          $p = 3$ (\ref{line:l2proj2d_midsweep:p3})
          with $|\Ecal_{h,p}| = 50$ elements for various values of the mesh
          regularization parameter $\alpha$ under the mesh parameterization
          in (\ref{eqn:mshparam1d}) and Figure~\ref{fig:mshparam-1dof-2d}.
          The parameter $\phibold$ is swept over the range
          $[x^*-0.75h,\,x^*+0.75h]$, where $x^*$ is the location
          of the discontinuity (\ref{line:l2proj2d_midsweep:xshk}). The position
          of $\bar{x}$ in (\ref{eqn:mshparam1d}) is indicated with
          (\ref{line:l2proj2d_midsweep:xloc}); when $\phibold = \bar{x}$
          nodes are uniformly distributed and $f_{msh}$ takes its minimum
          values. The red shaded region identifies the neighborhood of radius
          $h/2$ about $x^*$. The objective function with no mesh regularization
          ($\alpha = 0$) in the \emph{top left} plot shows the two desired
          properties discussed in Section~\ref{sec:track-obj}, namely, a
          minimum at $x^*$ and monotonically decreases to this minimum in a
          neighborhood centered at $x^*$ with radius $h/2$. The
          remaining plots show increasing values of the regularization
          parameter $\alpha$. The minimum does not noticeably move until
          $\alpha$ is sufficiently large at which point is moves from
          $x^*$ to $\bar{x}$, i.e., the objective function transitions from
          prioritizing discontinuity tracking to mesh regularization.}
 \label{fig:l2proj2d_midsweep}
\end{figure}
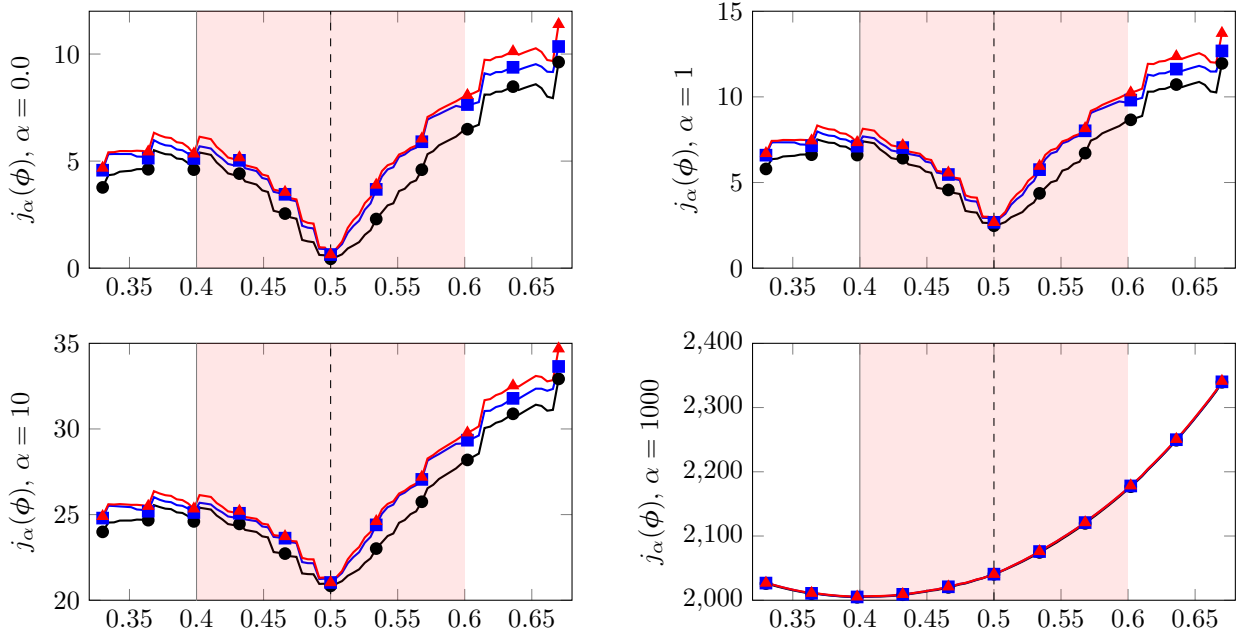
}

\subsubsection{Comparison with physics-based and error-based indicators:
               inviscid Burgers' equation}
\label{sec:track-obj-compare}
With the merit of the proposed discontinuity indicator verified in the
previous section, we compare it with other metrics that have been used
as discontinuity or refinement indicators. The first is an \emph{error-based}
indicator that measures the extent that a Galerkin solution corresponding to
a test and trial space $\Vcal_{h,p}$ satisfies the residual in an enriched
test space $\Vcal_{h'}^{p'}$
\begin{equation} \label{eqn:obj-res}
 f_{h,p}^{h',p'}(\ubm,\,\xbm) = r_{h,p}^{h',p'}(u_{h,p},\,\xbm)
\end{equation}
where $r_{h,p}^{h',p'}$ is the finite-dimensional Petrov-Galerkin residual
corresponding to the trial space $\Vcal_{h,p}$ and test space
$\Vcal_{h'}^{p'}$ defined in (\ref{eqn:claw-findim-dg-pg}). The rationale
behind this type of indicator is there will be no contribution from
regions of the domain away from the discontinuity that are well-resolved and
near the discontinuity the oscillations will be penalized by the additional
test functions. We will consider two instances of this class of
indicators: $h$-refinement where $h' = h/2$, $p' = p$ and $p$-refinement
where $h' = h$, $p' = p+1$. We also consider a \emph{physics-based} indicator
based on the Rankine-Hugoniot conditions \cite{baines2002multidimensional}.
Recall the Rankine-Hugoniot conditions for a steady state conservation law
reduce to
\begin{equation} \label{eqn:rank-hug}
 F(u^+) = F(u^-) \quad \text{on}~\Gamma_s,
\end{equation}
where $\Gamma_s$ is any surface in the domain and $u^+$, $u^-$ are the values
of the solution $u$ on either side of $\Gamma_s$. This will hold in
smooth regions of the domain or in the case $\Gamma_s$ is the shock surface;
however, it will be violated in regions in the domain with non-physical
jumps that arise from non-aligned meshes.

The condition in (\ref{eqn:rank-hug}) is converted to a discontinuity
indicator by integrating the residual over all inter-element faces
\begin{equation} \label{eqn:obj-rh}
 f_{rh}(\ubm,\,\xbm) = \sum_{\partial K \in \partial\Ecal_{h,p}} \int_{\partial K}
                 \norm{F(u_{h,p}^-)-F(u_{h,p}^+)}_\Wbm^2\,ds,
\end{equation}
where $\Wbm \in \Rbb^{N_c \times N_c}$ is a symmetric positive semi-definite
matrix that defines the local semi-norm. In this work, we take
$\Wbm = \ebm_1\ebm_1^T$, where $\ebm_1 \in \Rbb^d$ is the first canonical
unit vector.

To study the indicators in (\ref{eqn:obj-res})-(\ref{eqn:obj-rh}) we consider
the modified inviscid Burgers' equation with a discontinuous source term
(\ref{eqn:burg1d-eqn}) whose solution is given by (\ref{eqn:burg1d-soln}).
A complete description of this problem is deferred to
Section~\ref{sec:num-exp:burg1d}. Here, we let $\ubm$ denote the
discrete representation of the numerical solution of the inviscid Burgers'
equation in (\ref{eqn:burg1d-eqn}) using the test and trial space
$\Vcal_{h,p}$ and investigate the behavior of various indicators:
\begin{inparaenum}[1)]
 \item the proposed discontinuity tracking indicator in (\ref{eqn:obj-shk}),
 \item a residual-based indicator under $p$-refinement ($p' = p+1$),
 \item a residual-based indicator under $h$-refinement ($h = h/2$),
 \item and the Rankine-Hugoniot indicator
\end{inparaenum}
\begin{align}
  j_{shk}(\phibold) &= f_{shk}(\ubm(\xbm(\phibold)),\,\xbm(\phibold)) \label{eqn:obj-shk-redsp}\\
  j_{h,p}^{h,p+1}(\phibold) &= f_{h,p}^{h,p+1}(\ubm(\xbm(\phibold)),\,\xbm(\phibold)) \label{eqn:obj-resp-redsp}\\
  j_{h,p}^{h/2,p}(\phibold) &= f_{h,p}^{h/2,p}(\ubm(\xbm(\phibold)),\,\xbm(\phibold)) \label{eqn:obj-resh-redsp}\\
  j_{rh}(\phibold) &= f_{rh}(\ubm(\xbm(\phibold)),\,\xbm(\phibold)) \label{eqn:obj-rh-redsp}.
\end{align}
To study these indicators in isolation, the mesh regularization term is
excluded. We consider the same single degree of freedom parametrization from
the previous section, i.e., $\xbm(\phibold)$ is given by
(\ref{eqn:mshparam1d-disc}). The graphs of these functions over the range
$\phibold \in [x^*-0.75h,\,x^*+0.75h]$
where $x^* = 0.0$ is the discontinuity location, are shown in
Figure~\ref{fig:burg1d_midsweep} for a discretization with $40$ elements of
orders $p = 1,\,2,\,3$. It is clear from this figure that \emph{all} of the
functions considered possess a local minimum when the mesh is perfectly
aligned with the discontinuity ($\phibold = x^*$). However, the proposed
objective function is the only one that \emph{monotonically} approaches this
minimum in a neighborhood of radius $h/2$. The other indicators are extremely
oscillatory and in many cases are only monotonic in a neighborhood of radius
$h/10$ about the minimum. It is interesting to note that as the polynomial
order of the discretization increases, the physics- and error-based indicators
becomes more oscillatory, while the proposed indicator appears to be
agnostic to the polynomial order. This implies the physics- and error-based
indicators are not suitable candidates for the objective function in the
discontinuity-tracking framework since the mesh faces would have to be
initialized very close to the unknown discontinuity for any hope of convergence.
On the other hand, the monotonicity of the proposed indicator around
the discontinuity suggests that a gradient-based optimization solver will
be capable of aligning the mesh with discontinuities even in the worst
case where a face is $\Ocal(h/2)$ away from the discontinuity. In the
remainder of this document, we will solely consider the proposed
discontinuity indicator with mesh regularization term (\ref{eqn:obj}).

\ifbool{fastcompile}{}{
\begin{figure}
 \begin{tikzpicture}
\begin{groupplot}[
    group style={
        group size=2 by 2,
        vertical sep=1.3cm,
        horizontal sep=2.4cm
    },
    width=8.0cm,
    height=5.0cm,
    unbounded coords=jump,
    xticklabel style={
        /pgf/number format/fixed
    },
    scaled x ticks=false,
    xtick={-0.06, -0.03, 0.0, 0.03, 0.06}
]
\nextgroupplot[ymin=0, ymax=27, ylabel={$j_{shk}(\phibold)$}]
\addplot [black, solid, thick, mark=*, mark size=2, mark options={solid}, mark repeat={8}]  table[x index=0, y index=1, select coords between index={0}{80}] {dat/burg1d_xshk0_midnode_nel0040_p1.sweep.dat}; \label{line:burg1d_midsweep:p1}
\addplot [blue, solid, thick, mark=square*, mark size=2, mark options={solid}, mark repeat={8}]  table[x index=0, y index=1, select coords between index={0}{80}] {dat/burg1d_xshk0_midnode_nel0040_p2.sweep.dat}; \label{line:burg1d_midsweep:p2}
\addplot [red, solid, thick, mark=triangle*, mark size=2, mark options={solid}, mark repeat={8}]  table[x index=0, y index=1, select coords between index={0}{80}] {dat/burg1d_xshk0_midnode_nel0040_p3.sweep.dat}; \label{line:burg1d_midsweep:p3}
\addplot [draw=none, fill=red, fill opacity=0.1]
coordinates {
(  -0.05,  0)
(  -0.05,  27)
(   0.05,  27)
(   0.05,  0)
(  -0.05,  0)}; \label{line:burg1d_midsweep:h_over_2}
\addplot [black, dashed, thin]
coordinates {
(  0.0,  -2e-3)
(  0.0,  27)}; \label{line:burg1d_midsweep:xshk}

\nextgroupplot[ymin=0, ymax=80, ylabel={$j_{h,p}^{h,p+1}(\phibold)$}]
\addplot [black, solid, thick, mark=*, mark size=2, mark options={solid}, mark repeat={8}]  table[x index=0, y index=2, select coords between index={0}{80}] {dat/burg1d_xshk0_midnode_nel0040_p1.sweep.dat};
\addplot [blue, solid, thick, mark=square*, mark size=2, mark options={solid}, mark repeat={8}]  table[x index=0, y index=2, select coords between index={0}{80}] {dat/burg1d_xshk0_midnode_nel0040_p2.sweep.dat};
\addplot [red, solid, thick, mark=triangle*, mark size=2, mark options={solid}, mark repeat={8}]  table[x index=0, y index=2, select coords between index={0}{80}] {dat/burg1d_xshk0_midnode_nel0040_p3.sweep.dat};
\addplot [draw=none, fill=red, fill opacity=0.1]
coordinates {
(  -0.05,  0)
(  -0.05,  80)
(   0.05,  80)
(   0.05,  0)
(  -0.05,  0)};
\addplot [black, dashed, thin]
coordinates {
(  0.0,  -2e-3)
(  0.0,  80)};

\nextgroupplot[ymin=0, ymax=40,
               xlabel={$\phibold$ (position of node closest to shock)},
               ylabel={$j_{h,p}^{h/2,p}(\phibold)$}]
\addplot [black, solid, thick, mark=*, mark size=2, mark options={solid}, mark repeat={8}]  table[x index=0, y index=3, select coords between index={0}{80}] {dat/burg1d_xshk0_midnode_nel0040_p1.sweep.dat};
\addplot [blue, solid, thick, mark=square*, mark size=2, mark options={solid}, mark repeat={8}]  table[x index=0, y index=3, select coords between index={0}{80}] {dat/burg1d_xshk0_midnode_nel0040_p2.sweep.dat};
\addplot [red, solid, thick, mark=triangle*, mark size=2, mark options={solid}, mark repeat={8}]  table[x index=0, y index=3, select coords between index={0}{80}] {dat/burg1d_xshk0_midnode_nel0040_p3.sweep.dat};
\addplot [draw=none, fill=red, fill opacity=0.1]
coordinates {
(  -0.05,  0)
(  -0.05,  40)
(   0.05,  40)
(   0.05,  0)
(  -0.05,  0)};
\addplot [black, dashed, thin]
coordinates {
(  0.0,  -2e-3)
(  0.0,  40)};

\nextgroupplot[ymin=0, ymax=6.25,
               xlabel={$\phibold$ (position of node closest to shock)},
               ylabel={$j_{rh}(\phibold)$}]
\addplot [black, solid, thick, mark=*, mark size=2, mark options={solid}, mark repeat={8}]  table[x index=0, y index=5, select coords between index={0}{80}] {dat/burg1d_xshk0_midnode_nel0040_p1.sweep.dat};
\addplot [blue, solid, thick, mark=square*, mark size=2, mark options={solid}, mark repeat={8}]  table[x index=0, y index=5, select coords between index={0}{80}] {dat/burg1d_xshk0_midnode_nel0040_p2.sweep.dat};
\addplot [red, solid, thick, mark=triangle*, mark size=2, mark options={solid}, mark repeat={8}]  table[x index=0, y index=5, select coords between index={0}{80}] {dat/burg1d_xshk0_midnode_nel0040_p3.sweep.dat};
\addplot [draw=none, fill=red, fill opacity=0.1]
coordinates {
(  -0.05,  0)
(  -0.05,  6.25)
(   0.05,  6.25)
(   0.05,  0)
(  -0.05,  0)};
\addplot [black, dashed, thin]
coordinates {
(  0.0,  -2e-3)
(  0.0,  6.25)};

\end{groupplot}
\end{tikzpicture}
 \caption{A side-by-side comparison of the behavior of the proposed objective
          function with standard physics- and error-based objective functions
          (\emph{top left}: deviation from mean (\ref{eqn:obj-shk-redsp}),
           \emph{top right}: residual-based $p$-refinement
                             (\ref{eqn:obj-resp-redsp}),
           \emph{bottom left}: residual-based $h$-refinement
                               (\ref{eqn:obj-resh-redsp}),
           \emph{bottom right}: Rankine-Hugoniot conditions
                                (\ref{eqn:obj-rh-redsp}))
          corresponding to the solution of the inviscid Burgers' equation
          with a discontinuous source term in (\ref{eqn:burg1d-eqn}) using a
          piecewise polynomial basis of degree
          $p = 1$ (\ref{line:burg1d_midsweep:p1}),
          $p = 2$ (\ref{line:burg1d_midsweep:p2}),
          $p = 3$ (\ref{line:burg1d_midsweep:p3})
          with $|\Ecal_{h,p}| = 40$ elements without mesh regularization
          $\alpha = 0$. The single degree mesh parameterization from
          (\ref{eqn:mshparam1d}) and Figure~\ref{fig:mshparam-1dof-1d} is
          used and the parameter $\phibold$ is swept over the range
          $[x^*-0.75h,\,x^*+0.75h]$, where $x^*$ is the location
          of the discontinuity (\ref{line:burg1d_midsweep:xshk}).
          The red shaded region identifies the neighborhood of radius
          $h/2$ about $x^*$.
          Even with the added complexity of the conservation law, the
          proposed objective function attains a minimum when the mesh
          is aligned with the discontinuity and approaches this monotonically
          in a neighborhood of radius $h/2$. The other discontinuity
          indicators possess the correct minimum, but clearly do not satisfy
          monotonicity.}
 \label{fig:burg1d_midsweep}
\end{figure}
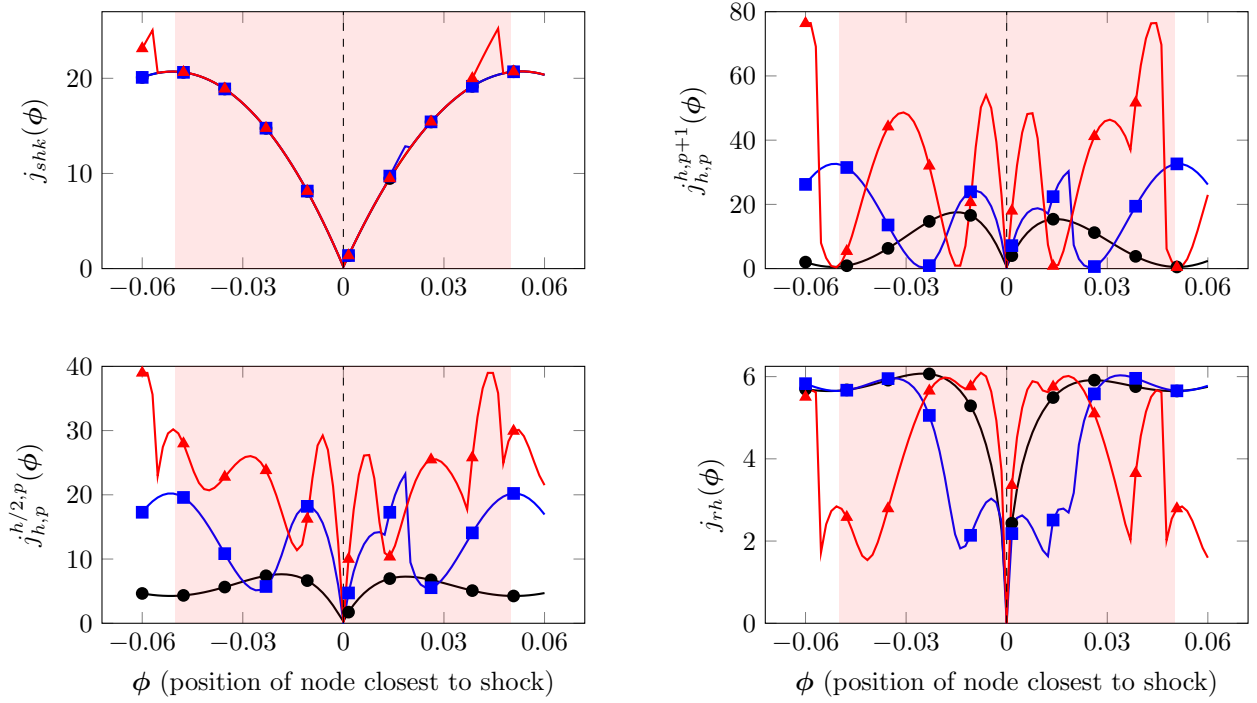
}

\subsection{Full space optimization solver}
\label{sec:track-solver}
The final ingredient to complete the description of the proposed discontinuity
tracking method is the numerical solver for the PDE-constrained optimization
problem in (\ref{eqn:claw-disc-opt1}) with the objective function in
(\ref{eqn:obj}) and constraints in (\ref{eqn:claw-disc}).
The first, and most common, PDE-constrained optimization solver is known as
the reduced space approach \cite{hinze2008optimization} where the PDE
constraint is explicitly enforced by defining the implicit function
$\ubm(\xbm)$ as the solution of $\rbm(\ubm,\,\xbm) = \zerobold$. The
PDE-constrained optimization problem reduces to the unconstrained optimization
problem
\begin{equation} \label{eqn:claw-disc-optred}
 \optunc{\phibold\in\Rbb^{N_\phibold}}{f(\ubm(\phibold),\,\xbm(\phibold))}
\end{equation}
where the state variable $\ubm$ is no longer an optimization variable as it can
be determined given $\phibold$. Standard unconstrained optimization solvers
such as quasi-Newton methods \cite{nocedal2006numerical} can be applied to solve
(\ref{eqn:claw-disc-optred}). The reduced space approach is popular because the
number of parameters is greatly reduced compared to the full space formulation
in (\ref{eqn:claw-disc-opt1}) due to the elimination of $\ubm$ and
sophisticated PDE solvers can be employed as a black box to query the
implicit function $\ubm(\xbm)$ and its gradient using the sensitivity or
adjoint method \cite{hinze2008optimization}.

Despite these advantages, the reduced space approach cannot, in general, be
used to solve the discontinuity-tracking PDE-constrained optimization problem
in (\ref{eqn:claw-disc-opt1}). The reduced space formulation implicitly assumes
the PDE constraint (\ref{eqn:claw-disc}) can be solved given $\xbm(\phibold)$
\emph{for any $\phibold$ in the feasible set} (taken as
$\Rbb^{N_\phibold}$). Recall from the discussion in
Section~\ref{sec:intro} that the PDE constraint cannot be solved unless
$\xbm = \xbm^*$, where $\xbm^*$ is the discrete representation of any
discontinuity-aligned mesh, or enough artificial viscosity is added to
stabilize the solution procedure. Since the mesh will not align with
discontinuities until convergence, the reduced space approach requires
repeated evaluation of the PDE solution at meshes not aligned with the
discontinuity and can only be used if sufficient viscosity, or
Laplacian-based diffusion, is added during intermediate iterations.
Such an approach would lead to a heuristic homotopy method and certainly
impact, if not prevent, convergence of the optimization procedure.

To overcome the numerous issues associated with the reduced space approach
to PDE-constrained optimization for the discontinuity-tracking problem in
(\ref{eqn:claw-disc-opt1}), we propose using the \emph{full space approach}
\cite{gunzburger2002perspectives, biros2005parallel}. In this setting, the
PDE state vector $\ubm$ and mesh parameters $\phibold$ are taken as
\emph{independent} optimization variables, as suggested in the formulation in
(\ref{eqn:claw-disc-opt1}), and converged to their optimal values
simultaneously. This implies the solution of the PDE is
\emph{never required away from a discontinuity-aligned mesh $\xbm^*$} and
overcomes the fundamental difficulty with the reduced space formulation.
Since $\ubm$ and $\phibold$ are varying
independently, the PDE constraint cannot be eliminated and a nonlinearly
constrained optimization solver is required. In this work, we use the
linesearch-based Sequential Quadratic Programming (SQP) method
\cite{boggs2000sequential} with augmented Lagrangian merit function
implemented in the well-known optimization software SNOPT
\cite{gill2002snopt, perez2012pyopt}
as our constrained optimization solver. Any SQP method will reduce the
nonlinear optimization problem in (\ref{eqn:claw-disc-opt1}) to a sequence of
quadratic programs of the form
\begin{equation} \label{eqn:sqp-quadprog}
 \optcona{\Delta\ubm\in\Rbb^{N_\ubm},\,\Delta\phibold\in\Rbb^{N_\phibold}}
     {q(\Delta\ubm,\,\Delta\phibold)}
     {\rbm(\ubm_k,\,\xbm(\phibold_k)) +
     \pder{\rbm}{\ubm}(\ubm_k,\,\xbm(\phibold_k))\Delta\ubm +
     \pder{\rbm}{\phibold}(\ubm_k,\,\xbm(\phibold_k))\Delta\phibold = \zerobold}
\end{equation}
where $\ubm_k$ and $\phibold_k$ correspond to the $k$th iterate of the
optimization variables and
$\func{q}{\Rbb^{N_\ubm}\times\Rbb^{N_\phibold}}{\Rbb}$ is the quadratic
Taylor approximation of the objective function where the Hessian is replaced
by the Hessian of the Lagrangian (or a symmetric positive definite
approximation)
\begin{equation*}
 \Lcal(\ubm,\,\phibold,\,\lambdabold) = f(\ubm,\,\phibold) -
                                        \lambdabold^T\rbm(\ubm,\,\xbm(\phibold))
\end{equation*}
with Lagrange multiplier $\lambdabold \in \Rbb^{N_\ubm}$. The solution of this
quadratic program provides search directions that are used to update the
optimization variables as
\begin{equation} \label{eqn:sqp-update}
 \begin{aligned}
  \ubm_{k+1} &= \ubm_k + \tau\Delta\ubm \\
  \phibold_{k+1} &= \phibold_k + \tau\Delta\phibold,
 \end{aligned}
\end{equation}
where $\tau$ is determined via a linesearch\footnote{Trust-region SQP
methods could be used in place of a linesearch-based SQP method. In this
case, the trust-region radius is determined at the beginning of the iteration
and the quadratic program is restricted to the trust region.} on a
merit function \cite{nocedal2006numerical}. It can be shown that the KKT
conditions of (\ref{eqn:sqp-quadprog}) are equivalent to the linearized KKT
conditions of (\ref{eqn:claw-disc-opt1})
\cite{boggs2000sequential,nocedal2006numerical}
and therefore SQP methods are equivalent to applying damped quasi-Newton
methods to the KKT system of the nonlinear optimization problem
(\ref{eqn:claw-disc-opt1}).

An optimization setting such as the linesearch-based SQP method in
(\ref{eqn:sqp-quadprog})-(\ref{eqn:sqp-update}) demonstrates the superiority
of the full space approach compared to the reduced space approach \emph{for
the discontinuity-tracking problem}. In the
reduced space formulation, the complete nonlinear solution of
$\rbm(\ubm,\,\xbm(\phibold)) = \zerobold$ for a given $\phibold$ is required,
which causes severe robustness issues if $\xbm(\phibold)$ is not aligned with
discontinuities and will likely not be possible since oscillations will
cause the solver to terminate at non-physical values. However, the full space
approach only requires the solution of the quadratic program in
(\ref{eqn:sqp-quadprog}) that only involves the \emph{linearized governing
equation}
\begin{equation*}
 \rbm(\ubm_k,\,\xbm(\phibold_k)) +
     \pder{\rbm}{\ubm}(\ubm_k,\,\xbm(\phibold_k))\Delta\ubm +
     \pder{\rbm}{\phibold}(\ubm_k,\,\xbm(\phibold_k))\Delta\phibold = \zerobold.
\end{equation*}
The solution $\ubm_k$ and the mesh $\xbm(\phibold_k)$ are updated
simultaneously based on the result of the linesearch. Therefore, at every
optimization iteration the mesh is updated alongside the solution and the
linesearch (or trust-region) prevents steps that cause non-physical values
as they would adversely affect the merit function.

While we use the general nonlinear optimization software SNOPT in this
work to demonstrate the merit of the proposed discontinuity-tracking
framework, there is a considerable opportunity for efficiency to develop a
specialized constrained optimization solver that takes advantage of the
considerable structure present in (\ref{eqn:claw-disc-opt1}), namely,
nonlinear equality constraints with an invertible partial Jacobian
$\displaystyle{\pder{\rbm}{\ubm}}$. Furthermore, a specialized optimizer
will be able to re-use data structures and operations from the PDE
implementation and incorporate homotopy measures (viscosity or temporal
relaxation) into the discrete PDE for particularly difficult problems. Such a
solver is the subject of ongoing work.


\subsection{Practical considerations: implementation, initialization,
            robustness, and efficiency}
\label{sec:track-practical}
In this section we discuss a number of practical considerations important in
the application of the proposed method to interesting problems. These include
relevant implementation details, initialization of the optimization problem in
(\ref{eqn:claw-disc-opt1}), and the expected efficiency compared to standard
steady state solvers or a reduced space optimization formulation.

\subsubsection*{Implementation}
The implementation of the proposed method consists of three main components:
the residual $\rbm(\ubm,\,\xbm)$ and objective $f(\ubm,\,\xbm)$ functions,
the definition of the mapping from optimization parameters to nodal
coordinates of the high-order mesh $\xbm = \Acal(\phibold)$, and
the nonlinear optimization solver for (\ref{eqn:claw-disc-opt1}),

The high-dimensional nature of the optimization problem in
(\ref{eqn:claw-disc-opt1}) with $N_\ubm + N_\phibold$ optimization variables,
requires a gradient-based optimization solver as derivative-free methods are
impractical in many dimensions. Any such solver will require
objective function and gradient evaluations as well as the constraint
function and gradient evaluations, i.e.,
\begin{equation} \label{eqn:fullsp-terms0}
 \begin{aligned}
  &f(\ubm,\,\xbm(\phibold)),\quad&&\pder{f}{\ubm}(\ubm,\,\xbm(\phibold)),\quad
  &&\pder{f}{\phibold}(\ubm,\,\xbm(\phibold)), \\
  &\rbm(\ubm,\,\xbm(\phibold)), &&\pder{\rbm}{\ubm}(\ubm,\,\xbm(\phibold)),
  &&\pder{\rbm}{\phibold}(\ubm,\,\xbm(\phibold)).
 \end{aligned}
\end{equation}
These are the same terms that would be required if a
reduced space optimization approach was used; however, in that setting,
the derivative terms would have to be combined using the sensitivity or
adjoint method to compute the \emph{total derivative} of the objective function
with respect to $\phibold$. The partial derivatives with respect to $\phibold$
are computed via a simple application of the chain rule
\begin{equation} \label{eqn:fullsp-terms-chain}
 \begin{aligned}
  \pder{f}{\phibold}(\ubm,\,\xbm(\phibold)) &=
  \pder{f}{\xbm}(\ubm,\,\xbm(\phibold)) \pder{\xbm}{\phibold}(\phibold) \\
  \pder{\rbm}{\phibold}(\ubm,\,\xbm(\phibold)) &=
  \pder{\rbm}{\xbm}(\ubm,\,\xbm(\phibold)) \pder{\xbm}{\phibold}(\phibold)
 \end{aligned}
\end{equation}
where $\displaystyle{\pder{\xbm}{\phibold} = \pder{\Acal}{\phibold}}$.

The terms involving the objective and residual function, i.e.,
\begin{equation} \label{eqn:fullsp-terms1}
 \begin{aligned}
  &f(\ubm,\,\xbm),\quad&&\pder{f}{\ubm}(\ubm,\,\xbm),\quad
  &&\pder{f}{\xbm}(\ubm,\,\xbm), \\
  &\rbm(\ubm,\,\xbm), &&\pder{\rbm}{\ubm}(\ubm,\,\xbm),
  &&\pder{\rbm}{\xbm}(\ubm,\,\xbm)
 \end{aligned}
\end{equation}
are solely dependent on the discretization of the conservation law and
can be completely separated from the mesh parametrization and optimization
solver. The implementation of the objective function and its derivatives are
straightforward due to the element-wise integral form of (\ref{eqn:obj-shk})
that can be expressed in terms of element mass matrices. The residual function
and its derivative with respect to the state variable $\ubm$ are standard terms
required by an implicit steady state solver. The derivative of the residual
function with respect to the nodal coordinates $\xbm$ is not standard and
usually only required for the purposes of gradient-based shape optimization
\cite{zahr2016dgopt}. However, the formulation of the conservation law on a
\emph{reference domain} in (\ref{eqn:claw-weak-dg-numflux}) implies the
elements over which the integrals are performed $K \in \Ecal_{h,p}$ are
independent of $\xbm$ and the dependence on $\xbm$ is restricted to the
\emph{pointwise} definition of the transformed fluxes. This implies the
implementation and data structures used to compute
$\displaystyle{\pder{\rbm}{\xbm}}$ will mimic those for
$\displaystyle{\pder{\rbm}{\ubm}}$.

The mesh parametrization can be treated independently of the PDE
discretization and must implement the mapping from $\phibold$ to the
mesh coordinates $\xbm = \Acal(\phibold)$ as well as the Jacobian of this
mapping, i.e.,
\begin{equation} \label{eqn:mshparam-and-deriv}
 \begin{aligned}
  \xbm(\phibold) &= \Acal(\phibold) \\
  \pder{\xbm}{\phibold}(\phibold) &= \pder{\Acal}{\phibold}(\phibold).
 \end{aligned}
\end{equation}
The most convenient and general mesh parametrization uses the
entire mesh as the optimization variables, i.e.,
\begin{equation}
 \begin{aligned}
  \phibold &= \xbm \\
  \Acal(\phibold) &= \phibold \\
  \pder{\Acal}{\phibold}(\phibold) &= \Ibm.
 \end{aligned}
\end{equation}
However, this burdens the optimization solver due to the large number of
optimization variables and likely introduces local minima that do not
correspond to discontinuity-aligned meshes. A general single-node mesh
parametrization for one-dimensional problems is presented in
Section~\ref{sec:track-obj-behavior} and used in
Sections~\ref{sec:track-obj-behavior}-\ref{sec:track-obj-compare},~
\ref{sec:num-exp:l2proj}-\ref{sec:num-exp:nozzle}. A general approach in
arbitrary dimensions that effectively reduces the number of optimization
variables selects $\phibold$ to represent coordinates of nodes in the
vicinity of discontinuities and smoothing based on linear elasticity is
used to determine the positions of the remaining nodes. This assumes
\emph{a-priori} knowledge of the shock region, which can easily be found
by introducing sufficient viscosity and finding the corresponding steady
state using traditional methods.
Regardless of the choice of mesh parametrization, the terms in
(\ref{eqn:fullsp-terms1}) and (\ref{eqn:mshparam-and-deriv}) can easily be
combined to yield the required derivative information for gradient-based
optimization solvers in (\ref{eqn:fullsp-terms0}).

The full space optimization solver is the final component required to
realize the proposed discontinuity-tracking method in practice, which
was discussed in Section~\ref{sec:track-solver}.

\subsubsection*{Initialization}
The initial guess for the optimization variables $\ubm$, $\phibold$ is a
particularly important consideration in gradient-based optimization with
nonlinear constraints due to the presence of local minima at meshes not
aligned with discontinuities. The objective function in (\ref{eqn:obj}) was
required to mostly eliminate local minima due to the \emph{monotonicity}
constraint, but cannot be guaranteed. In this work we use an initial guess
for the mesh parameters $\phibold_0$ such that the resulting mesh is
identical to the reference mesh $\Ecal_{h,p}$. This simple mesh initialization
strategy will be shown to be quite robust in the numerical examples in
Section~\ref{sec:num-exp}. Unlike the reduced space approach to optimization
that only requires an initial guess for $\phibold$, the full space approach
requires initialization of $\ubm$. Naturally, the solution of
$\rbm(\ubm,\,\xbm(\phibold_0)) = \zerobold$ is a desirable initial guess
since $\ubm_0$, $\phibold_0$ would satisfy the PDE constraint.
However, this is not possible since $\xbm(\phibold_0)$ corresponds to a
mesh that is not aligned with discontinuities. Instead, define
\begin{equation} \label{eqn:vclaw-disc}
 \rbm_\nu(\ubm,\,\xbm) = \zerobold
\end{equation}
as a discontinuous Galerkin discretization of a \emph{viscous} conservation
law with viscosity parameter $\nu$
\begin{equation} \label{eqn:vclaw}
 \nabla \cdot \Fcal(U) + \nabla \cdot \Fcal_\nu(U,\,\nabla U) = 0
 \quad\text{ in } \Gcal(\Omega_0,\,\mu),
\end{equation}
where $\Fcal_\nu$ becomes the zero map as $\nu \rightarrow 0$.
The reference \cite{arnold2002unified} provides an overview of various methods available
for the treatment of viscous terms in a DG setting. Since we are only looking
for an initial guess for the optimizer, the viscous terms in (\ref{eqn:vclaw})
can either correspond to physical viscosity, e.g., the viscous terms in the
Navier-Stokes equations, or non-physical Laplacian smoothing. Let $\bar\nu$
be a value of the viscosity parameter that is as close to zero as possible
but sufficiently large that the solution of
\begin{equation} \label{eqn:vclaw-disc0}
 \rbm_{\bar\nu}(\ubm,\,\xbm(\phibold_0)) = \zerobold
\end{equation}
can be found reliably. Then the solution of (\ref{eqn:vclaw-disc0}) is taken
as the initial guess for the state variable. 

These heuristics provide high quality initial guesses, particularly at
low polynomial orders such as $p=1$, but issues arise for high-order
elements due to the introduction of local minima not necessarily sensed
by the proposed indicator (\ref{eqn:obj-shk}). Therefore, we propose a
homotopy method on the polynomial order since even the $p=1$ solution finds a
quality approximation of the shock location and solution.
The transfer of the PDE solution $\ubm$ and mesh $\xbm$ to a higher order
mesh is a trivial injection operation and the corresponding optimization
parameters $\phibold$ are determined to exactly, if possible, reconstruct the
high-order mesh. Collectively, these initialization and homotopy
strategies and our particular choice of objective function seem to, in
practice, eliminate potential local minima that do not align with
discontinuities. Expections arise only in rather pedantic cases such as
a discontinuity positioned exactly in the center of an element in
one-dimension.

\subsubsection*{Robustness}
Recall from Section~\ref{sec:intro} and \ref{sec:track-obj} that robustness
is a key consideration in the design of the proposed method. In particular,
the objective function aligns the faces of the mesh with the discontinuities
to avoid severe oscillations in the solution and a full space optimization
approach is used to ensure the PDE solution is never required at non-aligned
meshes. These measures effectively solve the robustness issue of traditional
steady-state solvers provided the optimization problem is properly initialized
as discussed in the previous section.

\subsubsection*{Efficiency}
A natural question regarding efficiency of the proposed method arises given
the shift from a traditional nonlinear system of equations in $N_\ubm$
variables to a nonlinearly constrained optimization problem in
$N_\ubm+N_\phibold$ variables and $N_\ubm$ constraints. In fact, the proposed
discontinuity-tracking method is expected to perform similarly to a
traditional steady state solver \emph{on a given mesh}. This comes
from the observation that \emph{per iteration} an implicit steady state
solver will require evaluation of the residual and its Jacobian with
respect to $\ubm$ and a linear solve with this Jacobian. A full space
solver will also require these operations as well as the additional
terms in (\ref{eqn:fullsp-terms0}), i.e., the objective function and its
derivatives and the Jacobian of the residual with respect to $\xbm$. The
evaluation of the objective function and its derivatives are inexpensive
operations: they require a single pass over the mesh, have obvious data
structures, and a parallel implementation is straightforward. The evaluation
of the Jacobian of the residual with respect to $\xbm$ will incur a cost
similar to the evaluation of the Jacobian with respect to $\ubm$.

Despite a similar per iteration performance on a given mesh, the proposed
method that tracks discontinuities in the solution is expected to achieve
\emph{orders of magnitude} reduction in degrees of freedom compared to
traditional shock capturing methods such as limiting and artificial viscosity,
even when combined with adaptive mesh refinement. This reduction is expected because the discontinuity
is perfectly represented by the solution basis given there is sufficient
resolution for the isoparametric elements to track the \emph{geometry} of
the shock. Thus mesh resolution is driven by the geometry of the shock and
the smooth features in the solution away from the shock, which allows for
very coarse discretizations on high-order meshes as will be seen in
Section~\ref{sec:num-exp}.


\subsection{Summary of proposed discontinuity-tracking method}
Sections~\ref{sec:track-optform}-\ref{sec:track-practical} provide a complete
description of the proposed optimization-based discontinuity-tracking method.
The key ingredients of the proposed method are:
\begin{itemize}
 \item a high-order discretization of steady conservation laws in
       (\ref{eqn:claw-phys}) that supports discontinuities along element faces
       and is parametrized by the positions of the high-order mesh nodes;
       in this work we use a high-order discontinuous Galerkin method with
       isoparametric elements
 \item the PDE-constrained optimization formulation in (\ref{eqn:claw-disc-opt})
 \item a parametrization of the mesh as $\xbm = \Acal(\phibold)$ where
       $\phibold$ is a reduced set of optimization parameters and $\Acal$
       maps them to the nodes of the high-order mesh, possibly including
       mesh smoothing; the most general choice is the identity map, with
       other options presented in
       Section~\ref{sec:track-obj-behavior}~and~\ref{sec:num-exp:invisc-2d}
 \item the objective function in (\ref{eqn:obj}) that penalizes spurious
       oscillations and mesh distortion; numerical experiments on
       one-dimensional test cases in Section~\ref{sec:track-obj-behavior}-
       \ref{sec:track-obj-compare} show this indicator
       has a local minima at discontinuity-aligned mesh and decreases
       monotonically to this minima in a neighborhood of radius approximately
       $h/2$
 \item the full space PDE-constrained optimization solver that
       simultaneously converges the solution and mesh to their optimal values
       without ever requiring the PDE solution on a non-aligned mesh and
 \item the practical issue of the initialization of the non-convex optimization
       problem in (\ref{eqn:claw-disc-opt}) with \emph{viscosity} solutions of
       the conservation law or the result of the discontinuity tracking method
       at a lower polynomial order.
\end{itemize}
       
\section{Numerical experiments} \label{sec:num-exp}

\subsection{$L^2$ projection of discontinuous function onto piecewise polynomial
            basis}
\label{sec:num-exp:l2proj}
In this section, the proposed framework is investigated without the
complication of a conservation law. Instead, an $L^2$ projection of a
given discontinuous function onto a piecewise polynomial space
$\hat\Vcal_{h,p}(\xbm)$ induced by the discretization $\Ecal_{h,p}$ of the
reference domain $\Omega_0$ is considered. The purpose of this section is
three-fold:
\begin{inparaenum}[1)]
 \item demonstrate the merit of the proposed framework across a range of
       mesh sizes and polynomial orders,
 \item study the convergence behavior of the full space approach, and
 \item verify optimal $p+1$ convergence rates for the proposed
       discontinuity-tracking framework and compare to alternative methods.
\end{inparaenum}

The $L^2$ projection of a function $\eta \in L^2(\Omega)$ onto
$\hat\Vcal_{h,p}(\xbm)$ is defined as the function
$u \in \hat\Vcal_{h,p}(\xbm)$ such that
\begin{equation}
 \int_\Omega \psi u \,dV = \int_\Omega \psi \eta \, dV , \qquad
 \forall \psi \in \hat\Vcal_{h,p}(\xbm).
\end{equation}
Expansion of the solution and test function in the basis $\{\varphi_i\}$ and
restriction to element $K\in\Ecal_{h,p}$ leads to the discrete problem
\begin{equation} \label{eqn:l2proj-findim}
 \left(\int_{\Gcal(K,\,\xbm)} \varphi_i^K \varphi_j^K \, dV\right)\ubm_i^K =
 \int_{\Gcal(K,\,\xbm)} \varphi_i^K \eta \, dV, \qquad i \in \{1,\,\dots,\,N_p\}
\end{equation}
that can be written compactly as
\begin{equation} \label{eqn:l2proj-elem}
 \Mbm^K(\xbm) \cdot \ubm^K = \etabold^K(\xbm),
\end{equation}
where the element mass matrix and contribution of the function $\eta$ are
\begin{equation} \label{eqn:l2proj-terms}
 \Mbm_{ij}^K(\xbm) = \int_{\Gcal(K,\,\xbm)} \varphi_i^K\varphi_j^K \, dV, \qquad
 \etabold_i^K(\xbm) = \int_{\Gcal(K,\,\xbm)} \varphi_i^K \eta \, dV, \qquad
 i,\,j \in \{1,\,\dots,\,N_p\}.
\end{equation}
Therefore, the discrete residual for the case of an $L^2$ projection takes
the form
\begin{equation} \label{eqn:l2proj-res-disc}
 \rbm(\ubm,\,\xbm) = \Mbm(\xbm) \cdot \ubm - \etabold(\xbm) = \zerobold.
\end{equation}
In general, discontinuities present in $\eta$ will not align with element
faces and care must be taken to ensure the right-hand side of
(\ref{eqn:l2proj-findim}) is integrated accurately using quadrature. In this
work, for any element $K\in\Ecal_{h,p}$ such that $\Gcal(K,\,\xbm)$ contains
discontinuities of $\eta$, we compute special quadrature rules that do not
traverse these discontinuities. This special case of a $L^2$ projection
does not require the full space solver advocated for in
Section~\ref{sec:track-solver} since the governing equation in
(\ref{eqn:l2proj-res-disc}) can be solved given any $\xbm$ that does not
invert the mesh, i.e., the implicit function $\ubm(\xbm)$ is well-defined,
and  the reduced space approach is a valid option. However, we use the full
space approach to test the proposed method in this simple scenario. The mesh
$\xbm$ is initialized as a uniform mesh and $\ubm$ is initialized as the
solution of (\ref{eqn:l2proj-res-disc}) on the uniform mesh, regardless of
polynomial order.

We begin by considering the discontinuous, piecewise sinusoidal function
in one spatial dimension:
\begin{equation} \label{eqn:l2proj1d}
 \eta_k(x) =
 \begin{cases}
  0.25\sin(2\pi k x) - 0.5,& \qquad x < 0 \\
  0.25\sin(2\pi k x) + 0.5,& \qquad x > 0
 \end{cases}
\end{equation}
for all $x \in \Omega = (-1,\,1)$. The one-dimensional mesh $\Ecal_{h,p}$
is parametrized with one degree of freedom corresponding to each
\emph{non-boundary} node in the \emph{continuous, high-order} mesh.
The results of the discontinuity-tracking method, applied to the $L^2$
projection of (\ref{eqn:l2proj1d}) with $k = 2$, are shown in
Figures~\ref{fig:l2proj1d-trkplt}-\ref{fig:l2proj1d-trkplt-fail}
for a range of mesh sizes and polynomial orders. For an appropriate choice
of $\alpha$ ($\alpha = 1$ in this case), the discontinuity is tracked and
the solution is well-resolved with as few as five elements of degree
$p=3$. If $\alpha$ is too small, mesh regularization is ignored and the element
containing the discontinuity is crushed. If $\alpha$ is too large,
discontinuity-tracking is ignored and the mesh remains uniform.

\ifbool{fastcompile}{}{
\begin{figure}
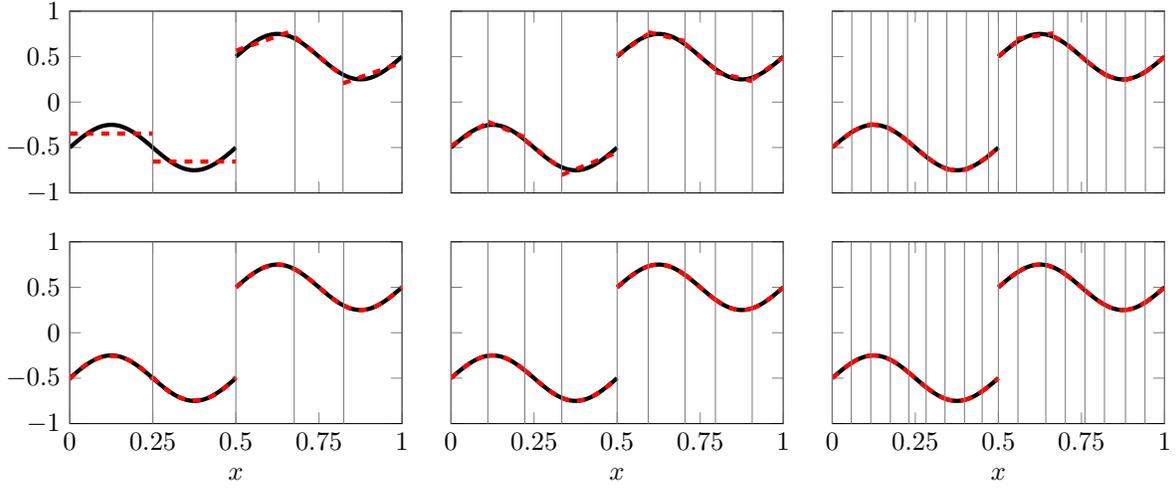

 \centering
 \begin{tikzpicture}
\begin{groupplot}[
    group style={
        group size=3 by 2,
        horizontal sep=0.65cm,
        vertical sep=0.65cm
    },
    xmin=0, xmax=1.0, ymin=-1.0, ymax=1.0,
    xtick={0,0.25,0.5,0.75,1},
    ytick={-1,-0.5,0,0.5,1},
    width=6.0cm,
    height=4.0cm
]

\nextgroupplot[xticklabels={,,}]
\input{dat/l2proj1d_exact.tikz} \label{line:l2proj1d_trkplt:exsln}
\input{dat/l2proj1d_xskel_fullsp_nel0005_p1.dev_from_mean.dev_from_unifmsh_1.trkplt.tikz}
\addplot [ultra thick, dashed, red]
coordinates {
(  -0.1,  0)
(  -0.05, 0)};  \label{line:trkplt:soln}
\addplot [solid, gray]
coordinates {
(  -0.1,  0)
(  -0.05, 0)};  \label{line:trkplt:mesh}

\nextgroupplot[xticklabels={,,}, yticklabels={,,}]
\input{dat/l2proj1d_exact.tikz}
\input{dat/l2proj1d_xskel_fullsp_nel0009_p1.dev_from_mean.dev_from_unifmsh_1.trkplt.tikz}

\nextgroupplot[xticklabels={,,}, yticklabels={,,}]
\input{dat/l2proj1d_exact.tikz}
\input{dat/l2proj1d_xskel_fullsp_nel0017_p1.dev_from_mean.dev_from_unifmsh_1.trkplt.tikz}

\nextgroupplot[xlabel={$x$}]
\input{dat/l2proj1d_exact.tikz}
\input{dat/l2proj1d_xskel_fullsp_nel0005_p3.dev_from_mean.dev_from_unifmsh_1.trkplt.tikz}

\nextgroupplot[xlabel={$x$}, yticklabels={,,}]
\input{dat/l2proj1d_exact.tikz}
\input{dat/l2proj1d_xskel_fullsp_nel0009_p3.dev_from_mean.dev_from_unifmsh_1.trkplt.tikz}

\nextgroupplot[xlabel={$x$}, yticklabels={,,}]
\input{dat/l2proj1d_exact.tikz}
\input{dat/l2proj1d_xskel_fullsp_nel0017_p3.dev_from_mean.dev_from_unifmsh_1.trkplt.tikz}

\end{groupplot}
\end{tikzpicture}
 \caption{Results of discontinuity-tracking framework using the objective
          function defined in (\ref{eqn:obj}) with $\alpha = 1$ applied to
          the $L^2$ projection of the discontinuous function in
          (\ref{eqn:l2proj1d}) with $k = 2$ onto discretizations with
          $|\Ecal_{h,p}| = 5$ (\emph{left column}),
          $|\Ecal_{h,p}| = 9$ (\emph{middle column}), $|\Ecal_{h,p}| = 17$
          (\emph{right column}) elements and polynomial order $p = 1$
          (\emph{top row}) and $p = 3$ (\emph{bottom row}). The $p=2$
          results are excluded for brevity as they are visually identical to
          the $p=3$ results. Legend: piecewise sinusoidal function in
          (\ref{eqn:l2proj1d}) with $k = 2$ (\ref{line:l2proj1d_trkplt:exsln}),
          mesh (\ref{line:trkplt:mesh}) and solution (\ref{line:trkplt:soln})
          obtained from discontinuity-tracking framework.}
 \label{fig:l2proj1d-trkplt}
\end{figure}
}

\ifbool{fastcompile}{}{
\begin{figure}
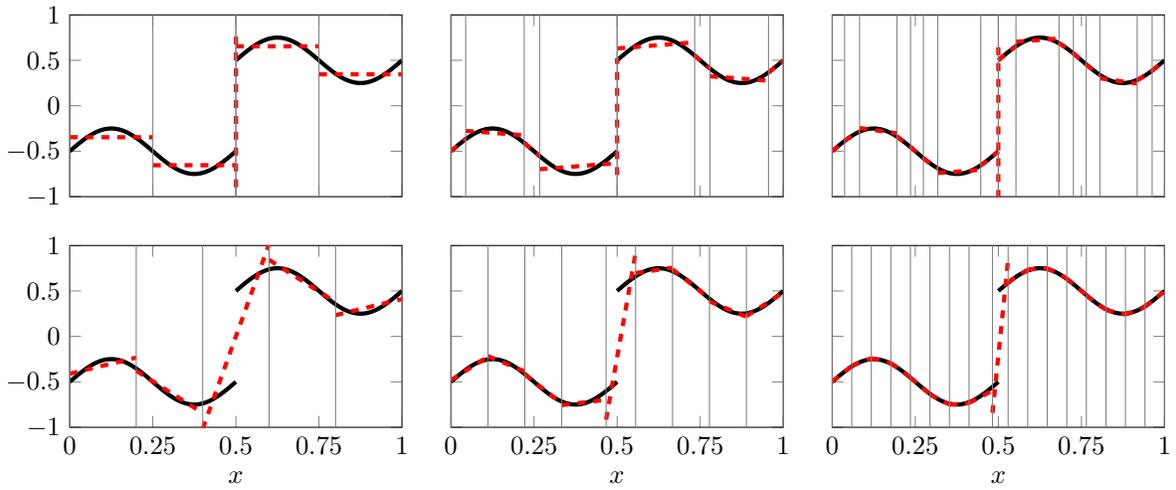

 \centering
 \begin{tikzpicture}
\begin{groupplot}[
    group style={
        group size=3 by 2,
        horizontal sep=0.65cm,
        vertical sep=0.65cm
    },
    xmin=0, xmax=1.0, ymin=-1.0, ymax=1.0,
    xtick={0,0.25,0.5,0.75,1},
    ytick={-1,-0.5,0,0.5,1},
    width=6.0cm,
    height=4.0cm
]

\nextgroupplot[xticklabels={,,}]
\input{dat/l2proj1d_exact.tikz} \label{line:l2proj1d_trkplt:exsln}
\input{dat/l2proj1d_xskel_fullsp_nel0005_p1.dev_from_mean.dev_from_unifmsh_0.trkplt.tikz}

\nextgroupplot[xticklabels={,,}, yticklabels={,,}]
\input{dat/l2proj1d_exact.tikz}
\input{dat/l2proj1d_xskel_fullsp_nel0009_p1.dev_from_mean.dev_from_unifmsh_0.trkplt.tikz}

\nextgroupplot[xticklabels={,,}, yticklabels={,,}]
\input{dat/l2proj1d_exact.tikz}
\input{dat/l2proj1d_xskel_fullsp_nel0017_p1.dev_from_mean.dev_from_unifmsh_0.trkplt.tikz}

\nextgroupplot[xlabel={$x$}]
\input{dat/l2proj1d_exact.tikz}
\input{dat/l2proj1d_xskel_fullsp_nel0005_p1.dev_from_mean.dev_from_unifmsh_10.trkplt.tikz}

\nextgroupplot[xlabel={$x$}, yticklabels={,,}]
\input{dat/l2proj1d_exact.tikz}
\input{dat/l2proj1d_xskel_fullsp_nel0009_p1.dev_from_mean.dev_from_unifmsh_10.trkplt.tikz}

\nextgroupplot[xlabel={$x$}, yticklabels={,,}]
\input{dat/l2proj1d_exact.tikz}
\input{dat/l2proj1d_xskel_fullsp_nel0017_p1.dev_from_mean.dev_from_unifmsh_10.trkplt.tikz}

\end{groupplot}
\end{tikzpicture}
 \caption{Results of discontinuity-tracking framework using the objective
          function defined in (\ref{eqn:obj}) with $\alpha = 0$
          (\emph{top row}) and $\alpha = 10$ (\emph{bottom row}) applied to
          the $L^2$ projection of the discontinuous function in
          (\ref{eqn:l2proj1d}) with $k = 2$ onto discretizations with
          $|\Ecal_{h,p}| = 5$ (\emph{left column}),
          $|\Ecal_{h,p}| = 9$ (\emph{middle column}), $|\Ecal_{h,p}| = 17$
          (\emph{right column}) elements and polynomial order $p = 1$.
          Legend: piecewise sinuosity function in (\ref{eqn:l2proj1d})
          with $k = 2$ (\ref{line:l2proj1d_trkplt:exsln}),
          mesh (\ref{line:trkplt:mesh}) and solution (\ref{line:trkplt:soln})
          obtained from discontinuity-tracking framework.}
 \label{fig:l2proj1d-trkplt-fail}
\end{figure}
}

The convergence history of the full space optimization solver corresponding
to the mesh with $17$ $p = 1$ elements for $\alpha = 1$ is provided in
Figure~\ref{fig:l2proj1d_trkcnvg}. It shows the solution of the discretized
PDE is only required \emph{at convergence}, i.e., once the objective is
minimized and the grid is aligned with the discontinuities. In fact, the
PDE constraint violation is relatively large until the objective function has
been reduced to nearly its final value and the discontinuity is closely
tracked. This is a feature of the full space approach and provides robustness
not possible with the reduced space approach since a PDE solution is never
required on non-aligned meshes.

\ifbool{fastcompile}{}{
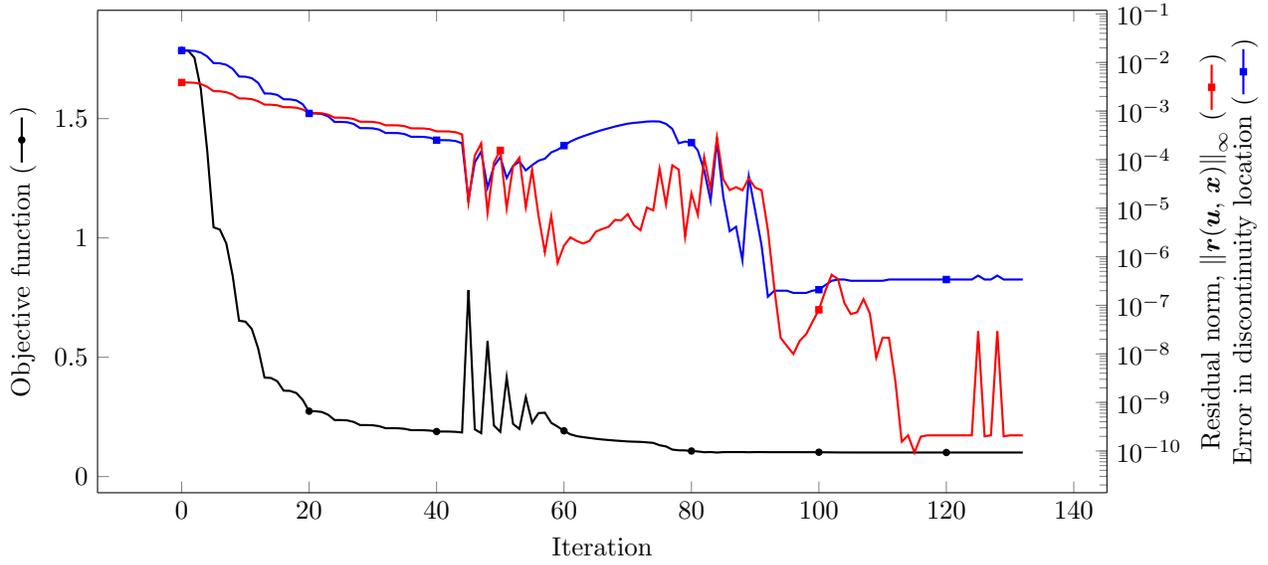
\begin{figure}
 \centering
 \begin{tikzpicture}

\begin{axis}[
  width=15cm, height=8cm,
  axis y line*=left,
  xlabel={Iteration},
  ylabel={Objective function (\ref{line:l2proj1d_trkcnvg:obj})}
]
\addplot [black, solid, thick, mark=*, mark size=1, mark options={solid}, mark repeat={20}]  table[x index=0, y index=1, select coords between index={0}{300}] {dat/l2proj1d_xskel_fullsp_nel0017_p1.dev_from_mean.dev_from_unifmsh_0p01.trkcnvg.dat}; \label{line:l2proj1d_trkcnvg:obj}
\end{axis}

\begin{axis}[
  width=15cm, height=8cm,
  ymode=log,
  axis y line*=right,
  axis x line=none,
  ylabel style={align=center},
  ylabel={Residual norm, $\norm{\rbm(\ubm,\,\xbm)}_\infty$
          (\ref{line:l2proj1d_trkcnvg:rinf}) \\
          Error in discontinuity location
          (\ref{line:l2proj1d_trkcnvg:shkloc})}
]
\addplot [blue, solid, thick, mark=square*, mark size=1, mark options={solid}, mark repeat={20}]  table[x index=0, y expr=abs(\thisrowno{2}-0.5), select coords between index={0}{300}] {dat/l2proj1d_xskel_fullsp_nel0017_p1.dev_from_mean.dev_from_unifmsh_0p01.trkcnvg.dat}; \label{line:l2proj1d_trkcnvg:shkloc}
\addplot [red, solid, thick, mark=square*, mark size=1, mark options={solid}, mark repeat={50}]  table[x index=0, y index=3, select coords between index={0}{300}] {dat/l2proj1d_xskel_fullsp_nel0017_p1.dev_from_mean.dev_from_unifmsh_0p01.trkcnvg.dat}; \label{line:l2proj1d_trkcnvg:rinf}
\end{axis}

\end{tikzpicture}
 \caption{Convergence history of full space optimization solver
          (SNOPT \cite{gill2002snopt}) applied to discontinuity-tracking
          optimization problem in (\ref{eqn:claw-disc-opt1}) with objective
          function in (\ref{eqn:obj}) ($\alpha = 1$) corresponding to the
          $L^2$ projection of the discontinuous function in (\ref{eqn:l2proj1d})
          with $k = 2$ onto a mesh with $|\Ecal_{h,p}| = 17$, $p = 1$ elements.
          Legend: Value of the objective function
          (\ref{line:l2proj1d_trkcnvg:obj}, \emph{left axis}),
          the residual norm $\norm{\rbm(\ubm,\,\xbm)}_\infty$
          (\ref{line:l2proj1d_trkcnvg:rinf}, \emph{right axis}),
          and error in the location of the discontinuity
          (\ref{line:l2proj1d_trkcnvg:shkloc}, \emph{right axis}).}
 \label{fig:l2proj1d_trkcnvg}
\end{figure}
}

The convergence of the $L^2$ projection of (\ref{eqn:l2proj1d}) with $k = 2$
onto a sequence of increasingly refined meshes with polynomial orders
$p = 1,\,2,\,3,\,4,\,5,\,6$ to the exact function is provided in
Figure~\ref{fig:l2proj1d_cnvg}. To ensure the mesh regularization
term does not impede convergence, the regularization parameter is set to
zero ($\alpha = 0$) and the single degree of freedom parametrization of
(\ref{fig:mshparam-1dof-1d}) is used. The $L^1$ error
\begin{equation} \label{eqn:l1err}
 e_{L^1}(u,\,u_{h,p}) = \int_\Omega |u - u_{h,p}|\,dV
\end{equation}
is used to measure convergence of the finite-dimensional solution $u_{h,p}$ to
the analytical solution $u$ given by (\ref{eqn:l2proj1d}) with $k = 2$.
Despite the $L^1$ norm being the weakest norm, it is appropriate to study
problems with discontinuities as it encapsulates the error in the smooth
region of the solution and discontinuity location and is expected to converge
at the appropriate rate of $\Ocal(h^{p+1})$ if the discontinuity location
converges at this rate. Figure~\ref{fig:l2proj1d_cnvg} shows the expected
convergence rate for polynomial orders up to $p = 6$ and a highly accurate
solution is possible with as few as four $p=6$ elements.

\ifbool{fastcompile}{}{
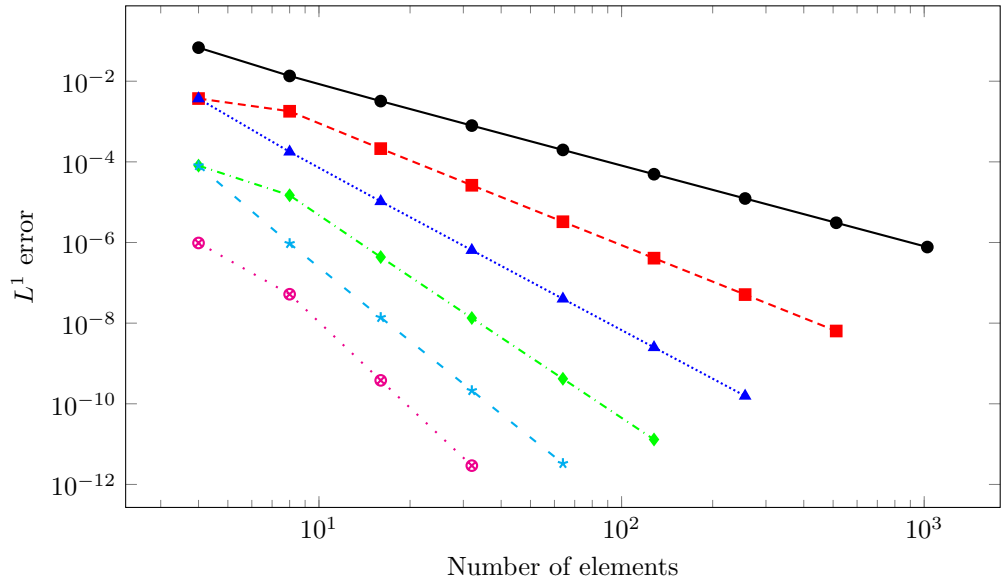
\begin{figure}
 \centering
 \begin{tikzpicture}
\begin{axis}[
    xmode=log,
    ymode=log,
    width=0.8\textwidth,
    height=0.5\textwidth,
    xlabel={Number of elements},
    ylabel={$L^1$ error}]

\addplot [black, solid, thick, mark=*, mark size=2, mark options={solid}]  table[x index=0, y index=2] {dat/l2proj1d_midnode_p1.dev_from_mean.cnvg.dat}; \label{line:l2proj1d_cnvg:p1}
\addplot [red, densely dashed, thick, mark=square*, mark size=2, mark options={solid}]  table[x index=0, y index=2] {dat/l2proj1d_midnode_p2.dev_from_mean.cnvg.dat}; \label{line:l2proj1d_cnvg:p2}
\addplot [blue, densely dotted, thick, mark=triangle*, mark size=2, mark options={solid}]  table[x index=0, y index=2] {dat/l2proj1d_midnode_p3.dev_from_mean.cnvg.dat}; \label{line:l2proj1d_cnvg:p3}
\addplot [green, dashdotted, thick, mark=diamond*, mark size=2, mark options={solid}]  table[x index=0, y index=2] {dat/l2proj1d_midnode_p4.dev_from_mean.cnvg.dat}; \label{line:l2proj1d_cnvg:p4}
\addplot [cyan, loosely dashed, thick, mark=star, mark size=2, mark options={solid}]  table[x index=0, y index=2] {dat/l2proj1d_midnode_p5.dev_from_mean.cnvg.dat}; \label{line:l2proj1d_cnvg:p5}
\addplot [magenta, loosely dotted, thick, mark=otimes, mark size=2, mark options={solid}]  table[x index=0, y index=2] {dat/l2proj1d_midnode_p6.dev_from_mean.cnvg.dat}; \label{line:l2proj1d_cnvg:p6}

\end{axis}

\end{tikzpicture}
 \caption{Convergence of discontinuity-tracking method using the single
          degree of freedom parametrization in (\ref{eqn:mshparam1d-disc})
          and no mesh regularization ($\alpha = 0$) applied to the $L^2$
          projection of the discontinuous function in (\ref{eqn:l2proj1d})
          with $k = 2$ for polynomial orders
          $p = 1$ (\ref{line:l2proj1d_cnvg:p1}),
          $p = 2$ (\ref{line:l2proj1d_cnvg:p2}),
          $p = 3$ (\ref{line:l2proj1d_cnvg:p3}),
          $p = 4$ (\ref{line:l2proj1d_cnvg:p4}),
          $p = 5$ (\ref{line:l2proj1d_cnvg:p5}),
          $p = 6$ (\ref{line:l2proj1d_cnvg:p6}).
          The expected convergence rates of $p+1$ are obtained in all cases.}
 \label{fig:l2proj1d_cnvg}
\end{figure}
}

With the convergence rates of the proposed method established, we compare
its performance to alternative methods:
\begin{inparaenum}[1)]
 \item uniform refinement and
 \item adaptive mesh refinement.
\end{inparaenum}
The mesh adaptation algorithm uses a straight-forward refinement strategy
where elements
with large values of the residual in an enriched ($h$-refined) test space
are split. The comparison is provided in Figure~\ref{fig:l2proj1d_cnvg_cmpr}
for $p = 1,\,3,\,5$. As expected, uniform refinement is only first-order
$\Ocal(h)$ accurate and the adaptive mesh refinement achieves the optimal
order of accuracy $\Ocal(h^{p+1})$ eventually. Our discontinuity-tracking
provides better accuracy for a given number of elements across a range of 
mesh sizes and polynomial orders since the discontinuity is accurately
tracked and the basis functions are optimally utilized in regions of the
domain where the solution is smooth.

\ifbool{fastcompile}{}{
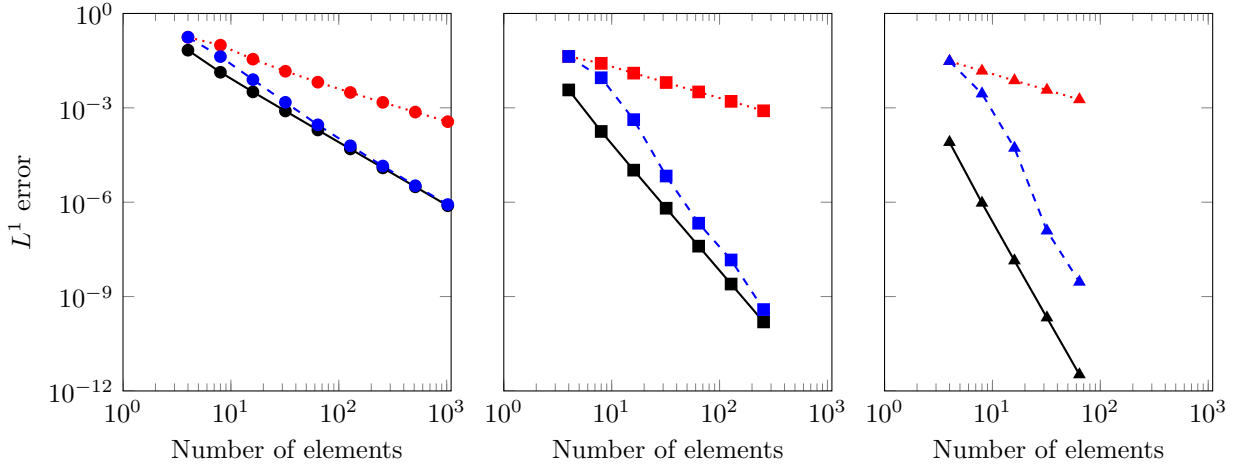
\begin{figure}
 \centering
 \begin{tikzpicture}
\begin{groupplot}[
    group style={
        group size=3 by 2,
        horizontal sep=0.7cm,
        vertical sep=0.9cm
    },
    xmin=1, xmax=1100,
    ymin=1e-12, ymax=1,
    xmode=log,
    ymode=log,
    width=0.36\textwidth,
    height=0.4\textwidth,
    xlabel={Number of elements}
]

\nextgroupplot[ylabel={$L^1$ error}]
\addplot [black, solid, thick, mark=*, mark size=2, mark options={solid}]  table[x index=0, y index=2] {dat/l2proj1d_midnode_p1.dev_from_mean.cnvg.dat}; \label{line:l2proj1d_cnvg_cmpr:p1}
\addplot [red, dotted, thick, mark=*, mark size=2, mark options={solid}]  table[x index=0, y index=2] {dat/l2proj1d_p1.dev_from_mean.cnvg_fixmsh.dat}; \label{line:l2proj1d_cnvg_cmpr:uref:p1}
\addplot [blue, dashed, thick, mark=*, mark size=2, mark options={solid}]  table[x index=0, y index=2] {dat/l2proj1d_p1.dev_from_mean.cnvg_amrmsh.dat}; \label{line:l2proj1d_cnvg_cmpr:amr:p1}

\nextgroupplot[yticklabels={,,}]
\addplot [black, solid, thick, mark=square*, mark size=2, mark options={solid}]  table[x index=0, y index=2] {dat/l2proj1d_midnode_p3.dev_from_mean.cnvg.dat}; \label{line:l2proj1d_cnvg_cmpr:p3}
\addplot [red, dotted, thick, mark=square*, mark size=2, mark options={solid}]  table[x index=0, y index=2] {dat/l2proj1d_p3.dev_from_mean.cnvg_fixmsh.dat}; \label{line:l2proj1d_cnvg_cmpr:uref:p3}
\addplot [blue, dashed, thick, mark=square*, mark size=2, mark options={solid}]  table[x index=0, y index=2] {dat/l2proj1d_p3.dev_from_mean.cnvg_amrmsh.dat}; \label{line:l2proj1d_cnvg_cmpr:amr:p3}

\nextgroupplot[yticklabels={,,}]
\addplot [black, solid, thick, mark=triangle*, mark size=2, mark options={solid}]  table[x index=0, y index=2] {dat/l2proj1d_midnode_p5.dev_from_mean.cnvg.dat}; \label{line:l2proj1d_cnvg_cmpr:p5}
\addplot [red, dotted, thick, mark=triangle*, mark size=2, mark options={solid}]  table[x index=0, y index=2] {dat/l2proj1d_p5.dev_from_mean.cnvg_fixmsh.dat}; \label{line:l2proj1d_cnvg_cmpr:uref:p5}
\addplot [blue, dashed, thick, mark=triangle*, mark size=2, mark options={solid}]  table[x index=0, y index=2] {dat/l2proj1d_p5.dev_from_mean.cnvg_amrmsh.dat}; \label{line:l2proj1d_cnvg_cmpr:amr:p5}

\end{groupplot}
\end{tikzpicture}
 \caption{Comparison of convergence of discontinuity-tracking method
          (single degree of freedom parametrization in
          (\ref{eqn:mshparam1d-disc}) and no mesh regularization, $\alpha = 0$)
          to uniform and adaptive mesh refinement at fixed polynomial orders
          for the $L^2$ projection of the discontinuous function in
          (\ref{eqn:l2proj1d}) with $k = 2$. The discontinuity-tracking method
          and adaptive mesh refinement achieve the expected $p+1$ convergence
          rate in the asymptotic regime while uniform refinement is limited to
          first order due to poor shock resolution. Legend:
          discontinuity-tracking with
          $p=1$ (\ref{line:l2proj1d_cnvg_cmpr:p1}),
          $p=3$ (\ref{line:l2proj1d_cnvg_cmpr:p3}),
          $p=5$ (\ref{line:l2proj1d_cnvg_cmpr:p5});
          uniform refinement with
          $p=1$ (\ref{line:l2proj1d_cnvg_cmpr:uref:p1}),
          $p=3$ (\ref{line:l2proj1d_cnvg_cmpr:uref:p3}),
          $p=5$ (\ref{line:l2proj1d_cnvg_cmpr:uref:p5});
          and adaptive mesh refinement with
          $p=1$ (\ref{line:l2proj1d_cnvg_cmpr:amr:p1}),
          $p=3$ (\ref{line:l2proj1d_cnvg_cmpr:amr:p3}),
          $p=5$ (\ref{line:l2proj1d_cnvg_cmpr:amr:p5}).}
 \label{fig:l2proj1d_cnvg_cmpr}
\end{figure}
}

\begin{figure}
 \centering
 \begin{tikzpicture}

\begin{groupplot}[
    group style={ 
        group size=2 by 2,
        horizontal sep=1.3cm
    },
    enlargelimits=false,
    axis on top, axis equal image,
    xmin=0, xmax=1, ymin=0, ymax=1,
    width=0.55\textwidth,
]

\nextgroupplot[ylabel={$y$}]
\addplot graphics [xmin=0, xmax=1, ymin=0, ymax=1] {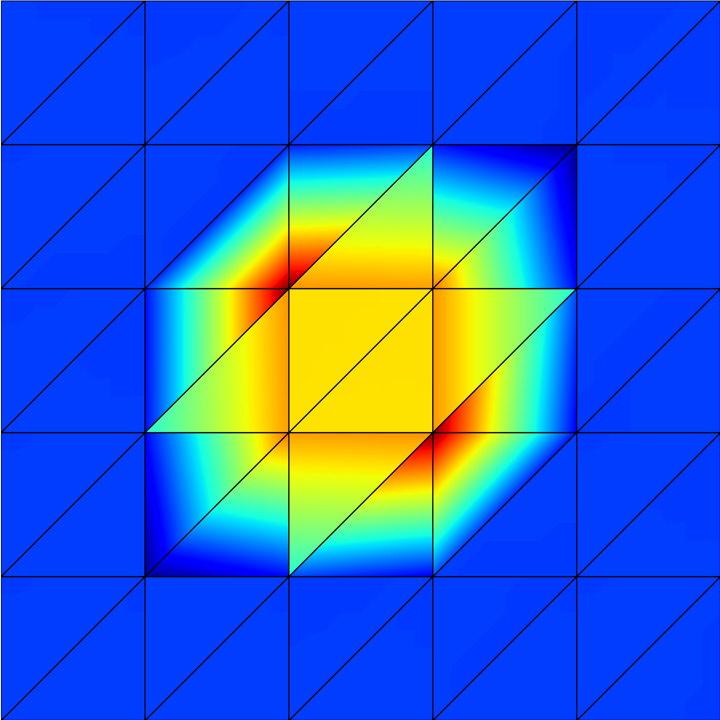};
\draw[black, dashed] (axis cs:0.5,0.5) circle [radius=0.25];

\nextgroupplot
\addplot graphics [xmin=0, xmax=1, ymin=0, ymax=1] {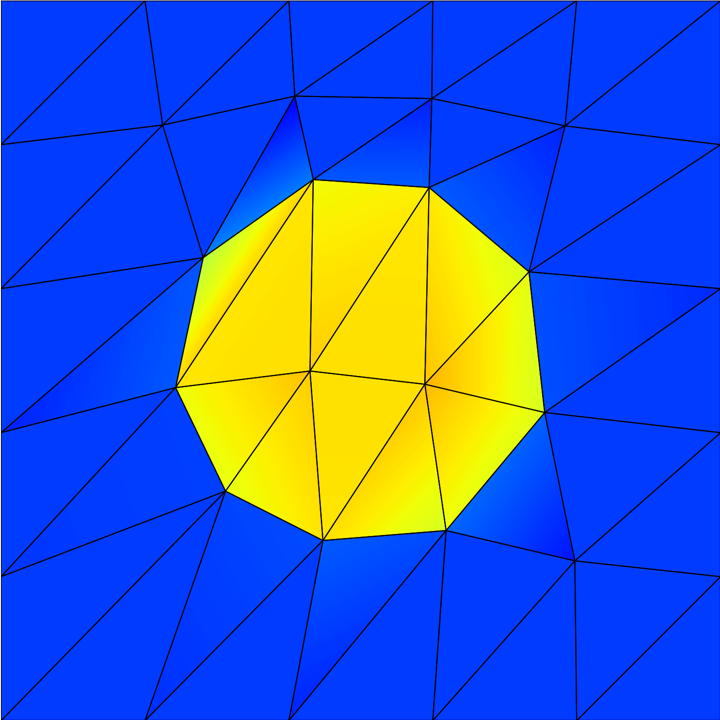};
\draw[black, dashed] (axis cs:0.5,0.5) circle [radius=0.25];

\nextgroupplot[xlabel={$x$}, ylabel={$y$}]
\addplot graphics [xmin=0, xmax=1, ymin=0, ymax=1] {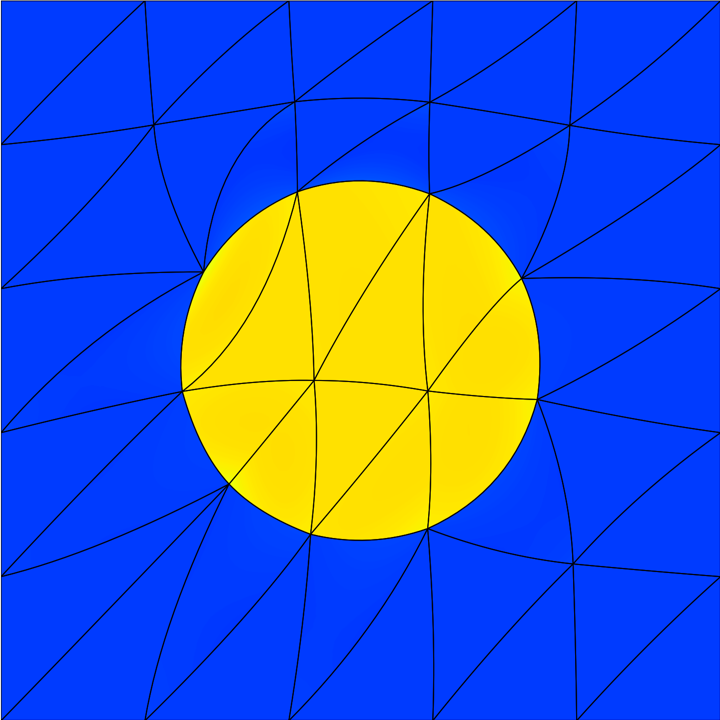};
\draw[black, dashed] (axis cs:0.5,0.5) circle [radius=0.25];

\nextgroupplot[xlabel={$x$}]
\addplot graphics [xmin=0, xmax=1, ymin=0, ymax=1] {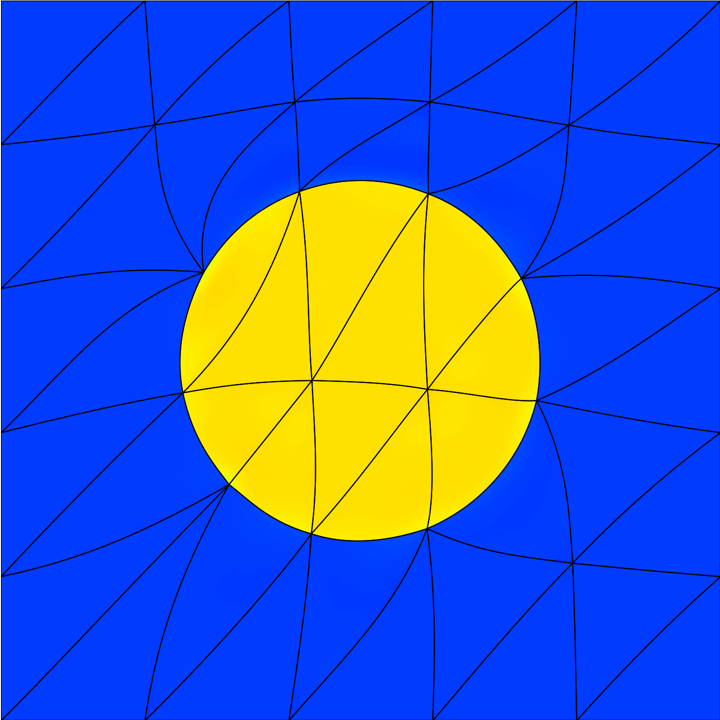};
\draw[black, dashed] (axis cs:0.5,0.5) circle [radius=0.25];

\end{groupplot}

\end{tikzpicture}
 \caption{Results of discontinuity-tracking framework using the objective
          function defined in (\ref{eqn:obj}) with $\alpha = 1$ applied to
          the $L^2$ projection of the function in (\ref{eqn:l2proj2d-smooth})
          with steep gradients onto discretizations with $|\Ecal_{h,p}| = 50$
          and polynomial order $p=1$ (\emph{top row}), $p = 2$
          (\emph{bottom left}), and $p=3$ (\emph{bottom right}).
          \emph{Top left}: $L^2$ projection of (\ref{eqn:l2proj2d-smooth})
          onto uniform non-aligned mesh with $50$ $p=1$ elements.
          \emph{Top right}: result of the discontinuity-tracking framework with
          $p=1$ elements using the non-aligned mesh and solution in the
          \emph{top left} as the initial guess. \emph{Bottom left}: result of
          the discontinuity-tracking framework with $p=2$ elements using the
          $p=1$ aligned mesh and solution in the \emph{top right} as the
          initial guess. \emph{Bottom right}: result of discontinuity-tracking
          framework with $p=3$ elements using the $p=2$ aligned mesh and
          solution in the \emph{bottom left} as the initial guess.}
 \label{fig:l2proj2d-trkplt}
\end{figure}
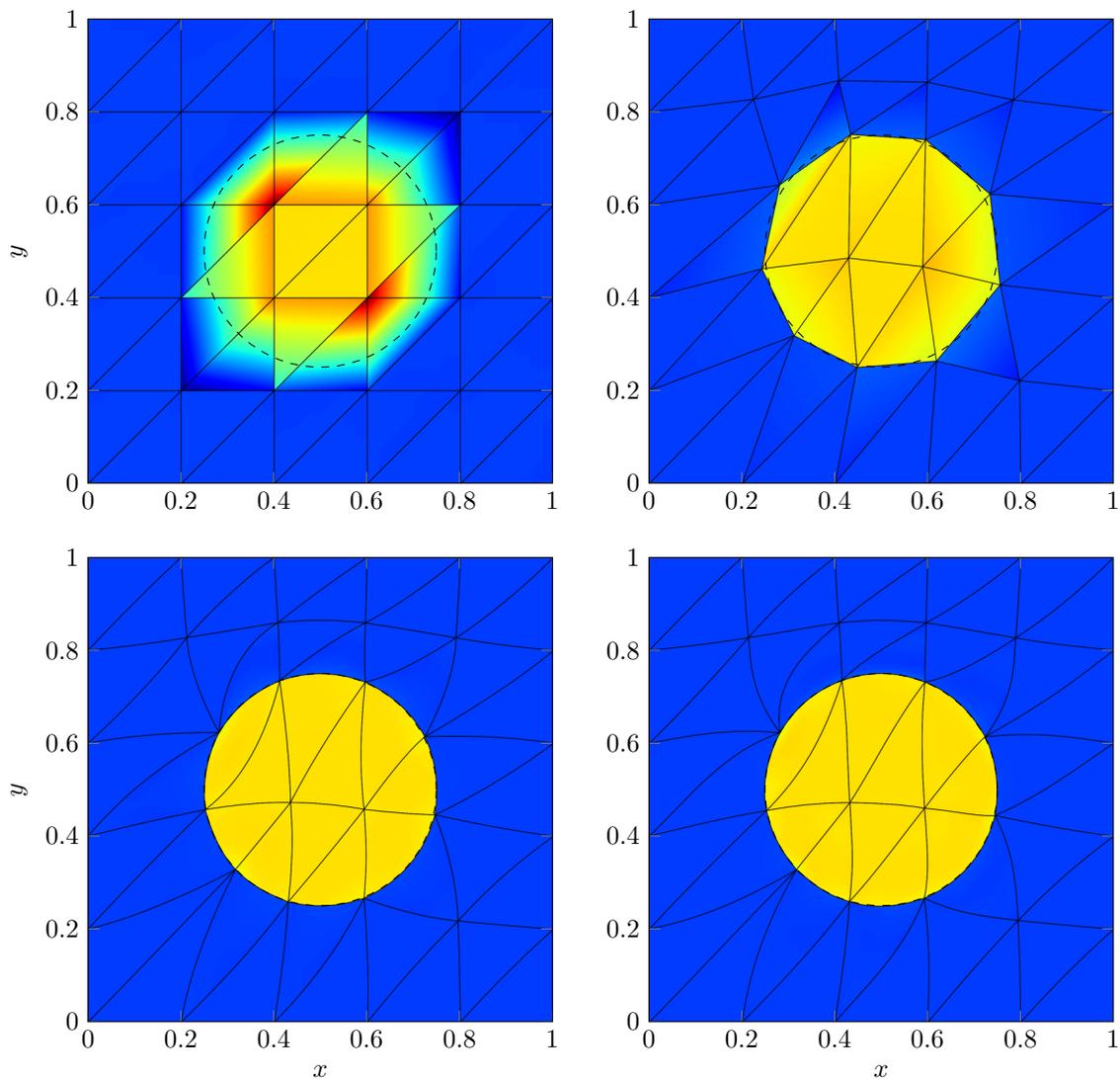

To close this section, we consider the $L^2$ projection of the discontinuous,
piecewise constant function in two spatial dimensions where the discontinuity
surface is a circle
\begin{equation} \label{eqn:l2proj2d}
 \eta(x) =
 \begin{cases}
  2, \qquad x^2 + y^2 < r^2 \\
  1, \qquad x^2 + y^2 > r^2
 \end{cases}
\end{equation}
for all $(x,\,y) \in \Omega = (0,\,1) \times (0,\,1)$. Due to the difficulty
associated with producing quadrature rules that do not cross discontinuities
in two-dimensions on curved meshes, the discontinuous function in
(\ref{eqn:l2proj2d}) is replaced with a smooth approximation
\begin{equation} \label{eqn:l2proj2d-smooth}
 \eta_\nu(x) = 1 + \frac12\left[1 - \tanh\left(\nu^{-1}
                       \left(\sqrt{(x-c_x)^2 + (y-c_y)^2}-r\right)\right)\right]
\end{equation}
where $\nu$ is a viscosity-like parameter that controls the steepness of the
smoothed discontinuity. In this work, $\nu = 10^{-3}$ is used. To ensure the
term on the right-hand-side of (\ref{eqn:l2proj-findim}) is integrated
accurately, a Gaussian quadrature rule with $20$ quadrature points is used.
The results of the discontinuity-tracking framework applied to the
$L^2$ projection of (\ref{eqn:l2proj2d-smooth}) onto a coarse mesh of only
$50$ elements with polynomial orders $p=1,\,2,\,3$ is shown in
Figure~\ref{fig:l2proj2d-trkplt}. The mesh is parametrized using all
continuous, high-order, non-boundary mesh nodes for the most general
parametrization. The $L^2$ projection onto a uniform,
non-aligned mesh is clearly a poor approximation to the piecewise constant
function where the discontinuity surface is a circle. Despite being a poor
approximation to the function in (\ref{eqn:l2proj2d}), this mesh and solution
are used as the initial guess for the discontinuity-tracking method with the
$p=1$ mesh. Subsequently, the mesh and solution produced by our
discontinuity-tracking method using polynomial order $p=\bar{p}$ is used to
initialize the method with polynomial order $p=\bar{p}+1$. The coarse mesh
with $p=2$ and $p=3$ track the discontinuity very well, sufficiently
regularize the mesh, and accurately approximate the piecewise constant
function in (\ref{eqn:l2proj2d}).

\subsection{Inviscid Burgers' equation with discontinuous source term}
\label{sec:num-exp:burg1d}
In this section, the proposed discontinuity-tracking framework is applied to
the inviscid, modified one-dimensional Burgers' equation with a discontinuous
source term from \cite{barter08}
\begin{equation} \label{eqn:burg1d-eqn}
 \pder{}{x}\left(\frac{1}{2}u^2\right) = \beta u + f(x), \quad
 \text{for}~x\in\Omega\subset\Rbb,
\end{equation}
where $\beta = -0.1$ and
\begin{equation} \label{eqn:burg1d-src}
 f(x) =
 \begin{cases}
  (2+\sin(\frac{\pi x}{2}))(\frac{\pi}{2}\cos(\frac{\pi x}{2}) - \beta), \qquad x < 0 \\
  (2+\sin(\frac{\pi x}{2}))(\frac{\pi}{2}\cos(\frac{\pi x}{2}) + \beta), \qquad x > 0.
 \end{cases}
\end{equation}
It can be shown that the solution of (\ref{eqn:burg1d-eqn}) is
\begin{equation} \label{eqn:burg1d-soln}
 u(x) =
 \begin{cases}
  \phantom{+}2 + \sin(\frac{\pi x}{2}), \qquad x < 0 \\
  -2 - \sin(\frac{\pi x}{2}), \qquad x > 0.
 \end{cases}
\end{equation}
The domain is taken as $\Omega = (-2,\,2)$ with Dirichlet boundary conditions
specified at both ends corresponding to the exact solution, i.e.,
$u(-2) = 2$ and $u(2) = -2$. The discretization of (\ref{eqn:burg1d-eqn})
proceeds according to the formulation in Section~\ref{sec:disc} and the
numerical flux is taken as the solution of the exact Riemann problem
corresponding to the time-dependent version of (\ref{eqn:burg1d-eqn}) without
source terms ($\beta = 0$ and $f(x) = 0$).

The one-dimensional mesh $\Ecal_{h,p}$ is parametrized with one degree of
freedom corresponding to each \emph{non-boundary} node in the
\emph{continuous, high-order} mesh. The results of the discontinuity-tracking
method, applied to the solution of the modified, inviscid Burgers' equation
are shown in Figure~\ref{fig:burg1d-trkplt} for a range of mesh sizes and
polynomial orders. For an appropriate choice
of $\alpha$ ($\alpha = 0.1$ in this case), the discontinuity is tracked and
the solution is well-resolved with as few as five $p=3$ elements.

Convergence of the proposed discontinuity-tracking method to the exact
solution in (\ref{eqn:burg1d-soln}) using a sequence of increasingly
refined meshes with polynomial orders $p = 1,\,2,\,3,\,4,\,5,\,6$ is provided
in Figure~\ref{fig:burg1d-cnvg}. Similarly to the previous section, the $L^1$
error is used to quantify the error in the solution and the single degree
of freedom mesh parametrization (\ref{fig:mshparam-1dof-1d}) without
regularization ($\alpha = 0$) is used to ensure mesh regularization does not
impede convergence. Figure~\ref{fig:burg1d-cnvg} shows the expected
convergence rate for most polynomial orders up to $p = 6$ and a highly
accurate solution is possible with as few as four $p=6$ elements.
Figure~\ref{fig:burg1d-cnvg-cmpr} compares the performance of the
discontinuity-tracking method to uniform and adaptive mesh refinement.
As expected, uniform refinement is only first-order $\Ocal(h)$ accurate and
adaptive mesh refinement nearly achieves the optimal order of accuracy
$\Ocal(h^{p+1})$. For this test case involving a conservation law, the
improved accuracy provided by tracking the discontinuity is more pronounced
than the previous case of an $L^2$ projection of an analytical function,
e.g., for the $p=3$ case, an error of $10^{-4}$ is possible with only
$10$ elements when the discontinuity is tracked while adaptive mesh
refinement requires over $700$ elements.

\ifbool{fastcompile}{}{
\begin{figure}
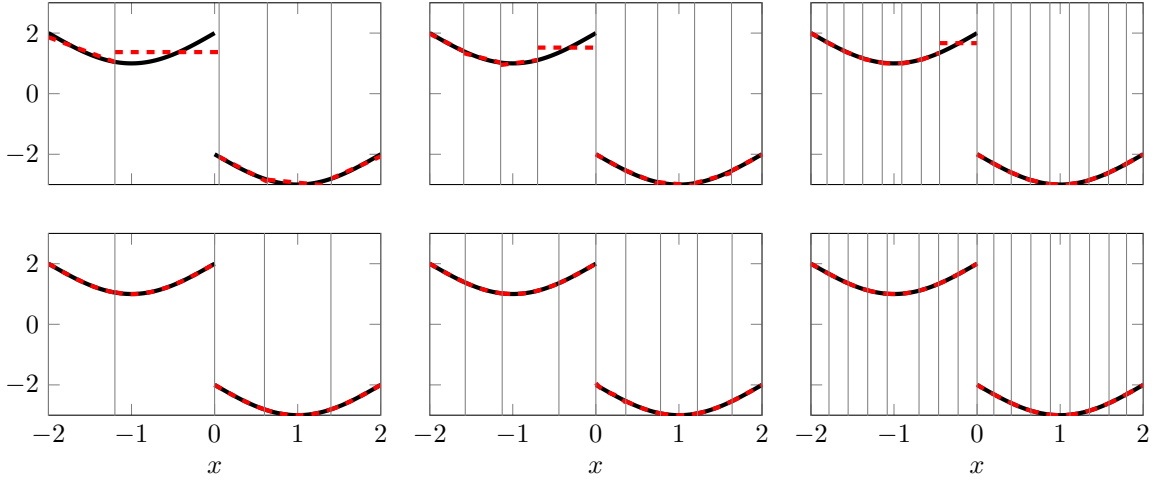

 \centering
 \begin{tikzpicture}
\begin{groupplot}[
    group style={
        group size=3 by 2,
        horizontal sep=0.65cm,
        vertical sep=0.65cm
    },
    xmin=-2.0, xmax=2.0, ymin=-3.0, ymax=3.0,
    width=6.0cm,
    height=4.0cm
]

\nextgroupplot[xticklabels={,,}]
\input{dat/burg1d_exact.tikz} \label{line:burg1d_trkplt:exsln}
\input{dat/burg1d_xshk0_xskel_nel0005_p1.dev_from_mean.dev_from_unifmsh_0p1.trkplt.tikz}

\nextgroupplot[xticklabels={,,}, yticklabels={,,}]
\input{dat/burg1d_exact.tikz}
\input{dat/burg1d_xshk0_xskel_nel0009_p1.dev_from_mean.dev_from_unifmsh_0p1.trkplt.tikz}

\nextgroupplot[xticklabels={,,}, yticklabels={,,}]
\input{dat/burg1d_exact.tikz}
\input{dat/burg1d_xshk0_xskel_nel0017_p1.dev_from_mean.dev_from_unifmsh_0p1.trkplt.tikz}

\nextgroupplot[xlabel={$x$}]
\input{dat/burg1d_exact.tikz}
\input{dat/burg1d_xshk0_xskel_nel0005_p3.dev_from_mean.dev_from_unifmsh_0p1.trkplt.tikz}

\nextgroupplot[xlabel={$x$}, yticklabels={,,}]
\input{dat/burg1d_exact.tikz}
\input{dat/burg1d_xshk0_xskel_nel0009_p3.dev_from_mean.dev_from_unifmsh_0p1.trkplt.tikz}

\nextgroupplot[xlabel={$x$}, yticklabels={,,}]
\input{dat/burg1d_exact.tikz}
\input{dat/burg1d_xshk0_xskel_nel0017_p3.dev_from_mean.dev_from_unifmsh_0p1.trkplt.tikz}

\end{groupplot}
\end{tikzpicture}
 \caption{Results of discontinuity-tracking framework (objective function
          defined in (\ref{eqn:obj}) with $\alpha = 0.1$) applied to the
          solution of the inviscid Burgers' equation (\ref{eqn:burg1d-eqn}) onto
          discretizations with $|\Ecal_{h,p}| = 5$ (\emph{left column}),
          $|\Ecal_{h,p}| = 9$ (\emph{middle column}), $|\Ecal_{h,p}| = 17$
          (\emph{right column}) elements and polynomial order $p = 1$
          (\emph{top row}) and $p = 3$ (\emph{bottom row}). The $p=2$
          results are excluded for brevity as they are visually identical to
          the $p=3$ results. Legend: exact solution
          (\ref{line:burg1d_trkplt:exsln}) of the modified inviscid Burgers'
          equation in (\ref{eqn:burg1d-eqn}),
          mesh (\ref{line:trkplt:mesh}) and solution (\ref{line:trkplt:soln})
          obtained from discontinuity-tracking framework.}
 \label{fig:burg1d-trkplt}
\end{figure}

\begin{figure}
 \centering
 \begin{tikzpicture}
\begin{axis}[
    xmin=1, xmax=1000,
    xmode=log,
    ymode=log,
    width=0.8\textwidth,
    height=0.5\textwidth,
    xlabel={Number of elements},
    ylabel={$L^1$ error}]

\addplot [black, solid, thick, mark=*, mark size=2, only marks, mark options={solid}]  table[x index=0, y index=5, select coords between index={0}{10}] {dat/burg1d_xshk0_midnode_p1.dev_from_mean.cnvg.dat}; \label{line:burg1d_cnvg:p1}
\addplot [black, solid, thin]
coordinates {
(4.0, 0.606865790997)
(200.0, 0.000287356203088)}; \label{line:burg1d_cnvg:slp2}

\addplot [red, solid, thick, mark=square*, mark size=2, only marks, mark options={solid}]  table[x index=0, y index=5, select coords between index={0}{10}] {dat/burg1d_xshk0_midnode_p2.dev_from_mean.cnvg.dat}; \label{line:burg1d_cnvg:p2}
\addplot [red, densely dashed, thin]
coordinates {
(2.000000e+00, 4.239738e-01)
(1.000000e+02, 2.020198e-06)}; \label{line:burg1d_cnvg:slp3}

\addplot [blue, solid, thick, mark=triangle*, mark size=2, only marks, mark options={solid}]  table[x index=0, y index=5, select coords between index={0}{10}] {dat/burg1d_xshk0_midnode_p3.dev_from_mean.cnvg.dat}; \label{line:burg1d_cnvg:p3}
\addplot [blue, densely dotted, thin]
coordinates {
(2.000000e+00, 5.660268e-02)
(5.000000e+01, 2.351972e-07)}; \label{line:burg1d_cnvg:slp4}

\addplot [green, solid, thick, mark=diamond*, mark size=2, only marks, mark options={solid}]  table[x index=0, y index=5, select coords between index={0}{10}] {dat/burg1d_xshk0_midnode_p4.dev_from_mean.cnvg.dat}; \label{line:burg1d_cnvg:p4}
\addplot [green, dashdotted, thin]
coordinates {
(4.000000e+00, 5.473396e-03)
(3.000000e+01, 8.822887e-08)}; \label{line:burg1d_cnvg:slp5}

\addplot [cyan, solid, thick, mark=star, mark size=2, only marks, mark options={solid}]  table[x index=0, y index=5, select coords between index={0}{10}] {dat/burg1d_xshk0_midnode_p5.dev_from_mean.cnvg.dat}; \label{line:burg1d_cnvg:p5}
\addplot [cyan, loosely dashed, thin]
coordinates {
(6.000000e+00, 8.081133e-06)
(2.000000e+01, 4.240036e-08)}; \label{line:burg1d_cnvg:slp6}

\addplot [magenta, solid, thick, mark=otimes, mark size=2, only marks, mark options={solid}]  table[x index=0, y index=5, select coords between index={0}{3}] {dat/burg1d_xshk0_midnode_p6.dev_from_mean.cnvg.dat}; \label{line:burg1d_cnvg:p6}
\addplot [magenta, loosely dotted, thin]
coordinates {
(4.000000e+00, 2.130902e-05)
(8.000000e+00, 5.229575e-08)}; \label{line:burg1d_cnvg:slp7}

\end{axis}

\end{tikzpicture}
 \caption{Convergence of discontinuity-tracking method using the single
          degree of freedom parametrization in (\ref{eqn:mshparam1d-disc})
          and no mesh regularization ($\alpha = 0$) applied to the modified
          inviscid Burgers' equation in (\ref{eqn:burg1d-eqn}) for
          polynomial orders
          $p = 1$ (\ref{line:burg1d_cnvg:p1}),
          $p = 2$ (\ref{line:burg1d_cnvg:p2}),
          $p = 3$ (\ref{line:burg1d_cnvg:p3}),
          $p = 4$ (\ref{line:burg1d_cnvg:p4}),
          $p = 5$ (\ref{line:burg1d_cnvg:p5}),
          $p = 6$ (\ref{line:burg1d_cnvg:p6}).
          The expected convergence rates of $p+1$ are obtained in most cases.
          The slopes of the best-fit lines to the data points in the asymptotic
          regime are: $\angle\,-1.95$ (\ref{line:burg1d_cnvg:slp2}),
                      $\angle\,-3.13$ (\ref{line:burg1d_cnvg:slp3}),
                      $\angle\,-3.85$ (\ref{line:burg1d_cnvg:slp4}),
                      $\angle\,-5.47$ (\ref{line:burg1d_cnvg:slp5}),
                      $\angle\,-4.36$ (\ref{line:burg1d_cnvg:slp6}),
                      $\angle\,-8.67$ (\ref{line:burg1d_cnvg:slp7}).}
 \label{fig:burg1d-cnvg}
\end{figure}
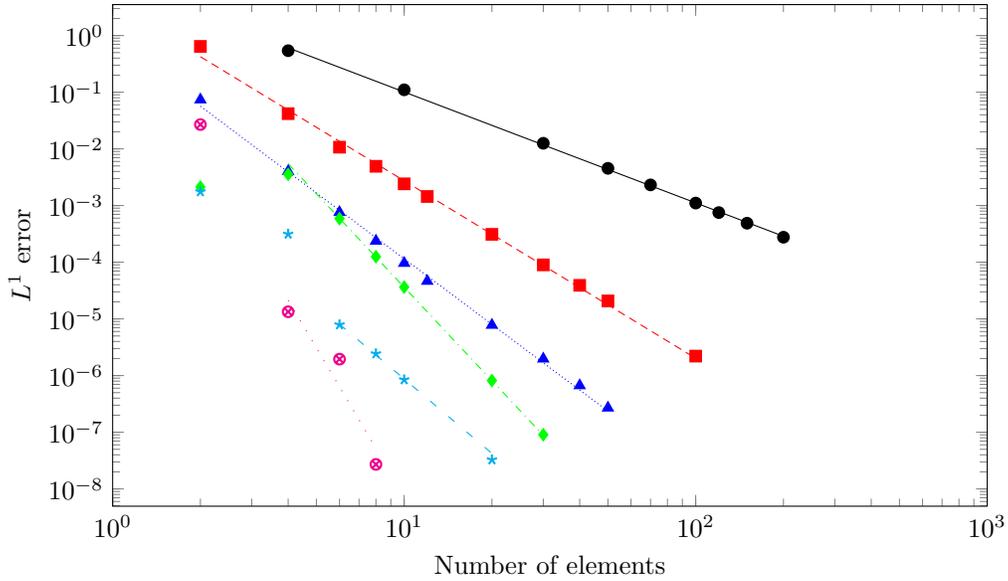

\begin{figure}
 \centering
 \begin{tikzpicture}
\begin{groupplot}[
    group style={
        group size=3 by 2,
        horizontal sep=0.7cm,
        vertical sep=0.9cm
    },
    xmin=1, xmax=1000,
    ymin=1e-7, ymax=1,
    xmode=log,
    ymode=log,
    width=0.36\textwidth,
    height=0.3\textwidth,
    xlabel={Number of elements}
]

\nextgroupplot[ylabel={$L^1$ error}]
\addplot [black, solid, thick, mark=*, mark size=2, mark options={solid}]  table[x index=0, y index=5, select coords between index={0}{10}] {dat/burg1d_xshk0_midnode_p1.dev_from_mean.cnvg.dat}; \label{line:burg1d_cnvg_cmpr:p1}
\addplot [black, dashed, thick, mark=*, mark size=2, mark options={solid}]  table[x index=0, y index=5, select coords between index={0}{100}] {dat/burg1d_xshk0_midnode_p1.cnvg_uref.dat}; \label{line:burg1d_cnvg_cmpr:uref:p1}
\addplot [black, dotted, thick, mark=*, mark size=2, mark options={solid}]  table[x index=0, y index=5, select coords between index={0}{100}] {dat/burg1d_xshk0_midnode_p1.cnvg_amr_res_h.dat}; \label{line:burg1d_cnvg_cmpr:amr:p1}

\nextgroupplot[yticklabels={,,}]
\addplot [red, solid, thick, mark=square*, mark size=2,  mark options={solid}]  table[x index=0, y index=5, select coords between index={0}{10}] {dat/burg1d_xshk0_midnode_p2.dev_from_mean.cnvg.dat}; \label{line:burg1d_cnvg_cmpr:p2}
\addplot [red, dashed, thick, mark=square*, mark size=2, mark options={solid}]  table[x index=0, y index=5, select coords between index={0}{100}] {dat/burg1d_xshk0_midnode_p2.cnvg_uref.dat}; \label{line:burg1d_cnvg_cmpr:uref:p2}
\addplot [red, dotted, thick, mark=square*, mark size=2, mark options={solid}]  table[x index=0, y index=5, select coords between index={0}{100}] {dat/burg1d_xshk0_midnode_p2.cnvg_amr_res_h.dat}; \label{line:burg1d_cnvg_cmpr:amr:p2}

\nextgroupplot[yticklabels={,,}]
\addplot [blue, solid, thick, mark=triangle*, mark size=2,  mark options={solid}]  table[x index=0, y index=5, select coords between index={0}{10}] {dat/burg1d_xshk0_midnode_p3.dev_from_mean.cnvg.dat}; \label{line:burg1d_cnvg_cmpr:p3}
\addplot [blue, dashed, thick, mark=triangle*, mark size=2, mark options={solid}]  table[x index=0, y index=5, select coords between index={0}{100}] {dat/burg1d_xshk0_midnode_p3.cnvg_uref.dat}; \label{line:burg1d_cnvg_cmpr:uref:p3}
\addplot [blue, dotted, thick, mark=triangle*, mark size=2, mark options={solid}]  table[x index=0, y index=5, select coords between index={0}{100}] {dat/burg1d_xshk0_midnode_p3.cnvg_amr_res_h.dat}; \label{line:burg1d_cnvg_cmpr:amr:p3}

\end{groupplot}
\end{tikzpicture}
 \caption{Comparison of convergence of discontinuity-tracking method
          (single degree of freedom parametrization in
          (\ref{eqn:mshparam1d-disc}) and no mesh regularization, $\alpha = 0$)
          to uniform and adaptive mesh refinement at fixed polynomial orders
          for modified inviscid Burgers' equation. The discontinuity-%
          tracking method achieves the expected $p+1$ convergence rate in the
          asymptotic regime while uniform refinement is limited to first order
          and adaptive mesh refinement falls short of optimal convergence
          rates due to poor shock resolution.
          Legend: discontinuity-tracking with
          $p=1$ (\ref{line:burg1d_cnvg_cmpr:p1}),
          $p=2$ (\ref{line:burg1d_cnvg_cmpr:p2}),
          $p=3$ (\ref{line:burg1d_cnvg_cmpr:p3});
          uniform refinement with
          $p=1$ (\ref{line:burg1d_cnvg_cmpr:uref:p1}),
          $p=2$ (\ref{line:burg1d_cnvg_cmpr:uref:p2}),
          $p=3$ (\ref{line:burg1d_cnvg_cmpr:uref:p3});
          and adaptive mesh refinement with
          $p=1$ (\ref{line:burg1d_cnvg_cmpr:amr:p1}),
          $p=2$ (\ref{line:burg1d_cnvg_cmpr:amr:p2}),
          $p=3$ (\ref{line:burg1d_cnvg_cmpr:amr:p3}).}
 \label{fig:burg1d-cnvg-cmpr}
\end{figure}
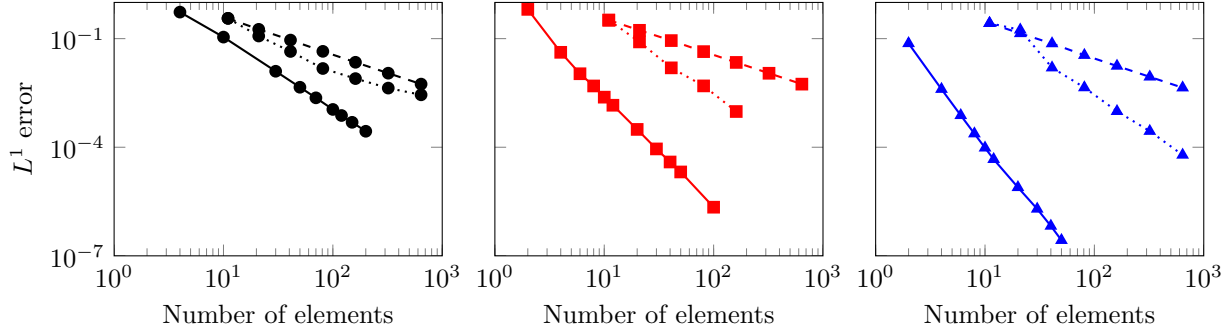
}

\subsection{Transonic, inviscid flow through nozzle}
\label{sec:num-exp:nozzle}
Our final one-dimensional example considers the relevant situation of transonic,
inviscid flow through a converging-diverging nozzle. The quasi-one-dimensional
Euler equations are used to model the inviscid, compressible flow in a
variable-area stream tube $A(x)$
\begin{equation} \label{eqn:nozzle-eqn}
 \begin{aligned}
  \pder{}{x}(A\rho u) &= 0, \\
  \pder{}{x}(A[\rho u^2 + p]) &= \frac{p}{A}\pder{A}{x}, \\
  \pder{}{x}(A[\rho E+p]u) &= 0
 \end{aligned}
\end{equation}
where $\rho$ is the fluid density, $u$ is the fluid velocity, $p$ is the
thermodynamic pressure, and
\begin{equation}
 \rho E = \rho e + \frac{1}{2}\rho u^2
\end{equation}
is the total energy. The pressure is related to $\rho E$ by the
equation of state
\begin{equation}
 p = (\gamma-1)\left(\rho E - \frac{1}{2}\rho u^2\right)
\end{equation}
for a perfect gas with ratio of specific heats $\gamma = 1.4$. The domain
is taken as $\Omega = (0,\,1)$ and the nozzle profile takes the form
\begin{equation}
 A(x) =
 \begin{cases}
  1 - (1-T)\cos(\pi(x-0.5)/0.8)^2 \quad &x \in [0.1,\,0.9] \\
  1 \quad &\text{otherwise}
 \end{cases}
\end{equation}
where $T = 0.8$ is the height of the nozzle throat. The boundary conditions
weakly impose the farfield conditions
$\rho_i = 1.0$, $u_i = 1.0$, $M_i = 0.40$ at the inflow and
$\rho_o = 1.0$, $u_o = 1.0$, $M_o = 0.45$ at the outflow using Roe's
approximate Riemann solver. A schematic of the domain, including boundary
conditions, is provided in Figure~\ref{fig:nozzle-geom}. The discretization
of the governing equations proceed according to the  formulation in
Section~\ref{sec:disc} and Roe's approximate Riemann solver corresponding
to the time-dependent version of (\ref{eqn:nozzle-eqn}) is used for
the numerical flux.

\ifbool{fastcompile}{}{
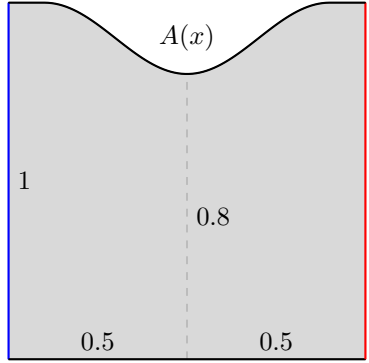
\begin{figure}
\centering
\begin{tikzpicture}
\begin{axis}[
axis equal image,
axis lines=none,
xmin=-0.1,
ymin=-0.1,
ymax=1.1,
xmax=1.1]
\addplot [fill opacity=0.5, fill=black!30!white, forget plot]
coordinates {
(  0.00000000,   0.00000000)
(  1.00000000,   0.00000000)
(  1.00000000,   1.00000000)
(  0.98989899,   1.00000000)
(  0.97979798,   1.00000000)
(  0.96969697,   1.00000000)
(  0.95959596,   1.00000000)
(  0.94949495,   1.00000000)
(  0.93939394,   1.00000000)
(  0.92929293,   1.00000000)
(  0.91919192,   1.00000000)
(  0.90909091,   1.00000000)
(  0.89898990,   0.99999685)
(  0.88888889,   0.99961947)
(  0.87878788,   0.99861544)
(  0.86868687,   0.99699107)
(  0.85858586,   0.99475658)
(  0.84848485,   0.99192603)
(  0.83838384,   0.98851723)
(  0.82828283,   0.98455161)
(  0.81818182,   0.98005412)
(  0.80808081,   0.97505306)
(  0.79797980,   0.96957988)
(  0.78787879,   0.96366902)
(  0.77777778,   0.95735764)
(  0.76767677,   0.95068546)
(  0.75757576,   0.94369445)
(  0.74747475,   0.93642858)
(  0.73737374,   0.92893356)
(  0.72727273,   0.92125653)
(  0.71717172,   0.91344579)
(  0.70707071,   0.90555047)
(  0.69696970,   0.89762023)
(  0.68686869,   0.88970496)
(  0.67676768,   0.88185446)
(  0.66666667,   0.87411810)
(  0.65656566,   0.86654454)
(  0.64646465,   0.85918144)
(  0.63636364,   0.85207510)
(  0.62626263,   0.84527023)
(  0.61616162,   0.83880964)
(  0.60606061,   0.83273396)
(  0.59595960,   0.82708142)
(  0.58585859,   0.82188756)
(  0.57575758,   0.81718507)
(  0.56565657,   0.81300352)
(  0.55555556,   0.80936922)
(  0.54545455,   0.80630503)
(  0.53535354,   0.80383022)
(  0.52525253,   0.80196036)
(  0.51515152,   0.80070721)
(  0.50505051,   0.80007866)
(  0.49494949,   0.80007866)
(  0.48484848,   0.80070721)
(  0.47474747,   0.80196036)
(  0.46464646,   0.80383022)
(  0.45454545,   0.80630503)
(  0.44444444,   0.80936922)
(  0.43434343,   0.81300352)
(  0.42424242,   0.81718507)
(  0.41414141,   0.82188756)
(  0.40404040,   0.82708142)
(  0.39393939,   0.83273396)
(  0.38383838,   0.83880964)
(  0.37373737,   0.84527023)
(  0.36363636,   0.85207510)
(  0.35353535,   0.85918144)
(  0.34343434,   0.86654454)
(  0.33333333,   0.87411810)
(  0.32323232,   0.88185446)
(  0.31313131,   0.88970496)
(  0.30303030,   0.89762023)
(  0.29292929,   0.90555047)
(  0.28282828,   0.91344579)
(  0.27272727,   0.92125653)
(  0.26262626,   0.92893356)
(  0.25252525,   0.93642858)
(  0.24242424,   0.94369445)
(  0.23232323,   0.95068546)
(  0.22222222,   0.95735764)
(  0.21212121,   0.96366902)
(  0.20202020,   0.96957988)
(  0.19191919,   0.97505306)
(  0.18181818,   0.98005412)
(  0.17171717,   0.98455161)
(  0.16161616,   0.98851723)
(  0.15151515,   0.99192603)
(  0.14141414,   0.99475658)
(  0.13131313,   0.99699107)
(  0.12121212,   0.99861544)
(  0.11111111,   0.99961947)
(  0.10101010,   0.99999685)
(  0.09090909,   1.00000000)
(  0.08080808,   1.00000000)
(  0.07070707,   1.00000000)
(  0.06060606,   1.00000000)
(  0.05050505,   1.00000000)
(  0.04040404,   1.00000000)
(  0.03030303,   1.00000000)
(  0.02020202,   1.00000000)
(  0.01010101,   1.00000000)
(  0.00000000,   1.00000000)
(  0.00000000,   0.00000000)};

\addplot [thin, black, dashed, opacity=0.25, forget plot]
coordinates {
(  0.50000000,   0.00000000)
(  0.50000000,   0.80000000)};

\addplot [thick, color=blue]
coordinates {
(  0.00000000,   1.00000000)
(  0.00000000,   0.00000000)};\label{line:nozzle0:inflow}

\addplot [thick, color=black]
coordinates {
(  0.00000000,   0.00000000)
(  1.00000000,   0.00000000)};\label{line:nozzle0:wall}

\addplot [thick, color=red]
coordinates {
(  1.00000000,   0.00000000)
(  1.00000000,   1.00000000)};\label{line:nozzle0:outflow}

\addplot [thick, color=black]
coordinates {
(  1.00000000,   1.00000000)
(  0.98989899,   1.00000000)
(  0.97979798,   1.00000000)
(  0.96969697,   1.00000000)
(  0.95959596,   1.00000000)
(  0.94949495,   1.00000000)
(  0.93939394,   1.00000000)
(  0.92929293,   1.00000000)
(  0.91919192,   1.00000000)
(  0.90909091,   1.00000000)
(  0.89898990,   0.99999685)
(  0.88888889,   0.99961947)
(  0.87878788,   0.99861544)
(  0.86868687,   0.99699107)
(  0.85858586,   0.99475658)
(  0.84848485,   0.99192603)
(  0.83838384,   0.98851723)
(  0.82828283,   0.98455161)
(  0.81818182,   0.98005412)
(  0.80808081,   0.97505306)
(  0.79797980,   0.96957988)
(  0.78787879,   0.96366902)
(  0.77777778,   0.95735764)
(  0.76767677,   0.95068546)
(  0.75757576,   0.94369445)
(  0.74747475,   0.93642858)
(  0.73737374,   0.92893356)
(  0.72727273,   0.92125653)
(  0.71717172,   0.91344579)
(  0.70707071,   0.90555047)
(  0.69696970,   0.89762023)
(  0.68686869,   0.88970496)
(  0.67676768,   0.88185446)
(  0.66666667,   0.87411810)
(  0.65656566,   0.86654454)
(  0.64646465,   0.85918144)
(  0.63636364,   0.85207510)
(  0.62626263,   0.84527023)
(  0.61616162,   0.83880964)
(  0.60606061,   0.83273396)
(  0.59595960,   0.82708142)
(  0.58585859,   0.82188756)
(  0.57575758,   0.81718507)
(  0.56565657,   0.81300352)
(  0.55555556,   0.80936922)
(  0.54545455,   0.80630503)
(  0.53535354,   0.80383022)
(  0.52525253,   0.80196036)
(  0.51515152,   0.80070721)
(  0.50505051,   0.80007866)
(  0.49494949,   0.80007866)
(  0.48484848,   0.80070721)
(  0.47474747,   0.80196036)
(  0.46464646,   0.80383022)
(  0.45454545,   0.80630503)
(  0.44444444,   0.80936922)
(  0.43434343,   0.81300352)
(  0.42424242,   0.81718507)
(  0.41414141,   0.82188756)
(  0.40404040,   0.82708142)
(  0.39393939,   0.83273396)
(  0.38383838,   0.83880964)
(  0.37373737,   0.84527023)
(  0.36363636,   0.85207510)
(  0.35353535,   0.85918144)
(  0.34343434,   0.86654454)
(  0.33333333,   0.87411810)
(  0.32323232,   0.88185446)
(  0.31313131,   0.88970496)
(  0.30303030,   0.89762023)
(  0.29292929,   0.90555047)
(  0.28282828,   0.91344579)
(  0.27272727,   0.92125653)
(  0.26262626,   0.92893356)
(  0.25252525,   0.93642858)
(  0.24242424,   0.94369445)
(  0.23232323,   0.95068546)
(  0.22222222,   0.95735764)
(  0.21212121,   0.96366902)
(  0.20202020,   0.96957988)
(  0.19191919,   0.97505306)
(  0.18181818,   0.98005412)
(  0.17171717,   0.98455161)
(  0.16161616,   0.98851723)
(  0.15151515,   0.99192603)
(  0.14141414,   0.99475658)
(  0.13131313,   0.99699107)
(  0.12121212,   0.99861544)
(  0.11111111,   0.99961947)
(  0.10101010,   0.99999685)
(  0.09090909,   1.00000000)
(  0.08080808,   1.00000000)
(  0.07070707,   1.00000000)
(  0.06060606,   1.00000000)
(  0.05050505,   1.00000000)
(  0.04040404,   1.00000000)
(  0.03030303,   1.00000000)
(  0.02020202,   1.00000000)
(  0.01010101,   1.00000000)
(  0.00000000,   1.00000000)};\label{line:nozzle0:wall}

\node[above]    at    (axis cs:0.25, 0.0) {$0.5$};
\node[above]    at    (axis cs:0.75, 0.0) {$0.5$};
\node[right]    at    (axis cs:0.0, 0.5) {$1$};
\node[right]    at    (axis cs:0.5, 0.4) {$0.8$};
\node[above]    at    (axis cs:0.5, 0.84) {$A(x)$};
\end{axis}
\end{tikzpicture} \\
\caption{Geometry and boundary conditions for nozzle flow. Boundary conditions:
         inviscid wall (\ref{line:nozzle0:wall}),
         inflow (\ref{line:nozzle0:inflow}),
         outflow (\ref{line:nozzle0:outflow}).}
\label{fig:nozzle-geom}
\end{figure}
}

The one-dimensional mesh $\Ecal_{h,p}$ is parametrized with one degree of
freedom corresponding to each \emph{non-boundary} node in the
\emph{continuous, high-order} mesh. The results of the discontinuity-tracking
method, applied to the solution of the quasi-1d Euler equations
are shown in Figure~\ref{fig:nozzle-trkplt} for a range of
mesh sizes and polynomial orders. For an appropriate choice
of $\alpha$ ($\alpha = 0.1$ in this case), the discontinuity is tracked and
the solution is well-resolved with as few as eight $p=2$ elements.

\ifbool{fastcompile}{}{
\begin{figure}
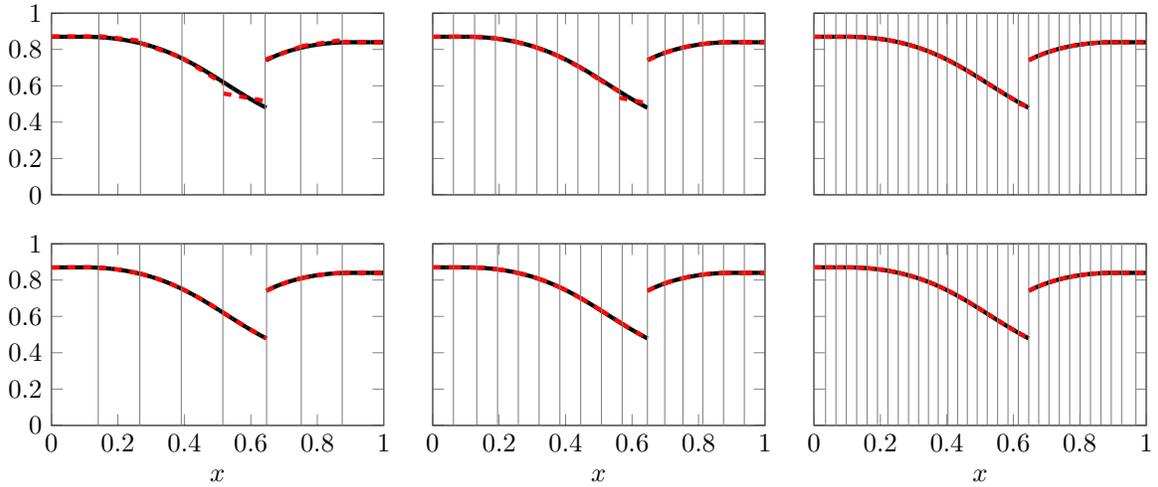

 \centering
 \begin{tikzpicture}
\begin{groupplot}[
    group style={
        group size=3 by 2,
        horizontal sep=0.65cm,
        vertical sep=0.65cm
    },
    xmin=0.0, xmax=1.0, ymin=0.0, ymax=1.0,
    width=6.0cm,
    height=4.0cm
]

\nextgroupplot[xticklabels={,,}]
\input{dat/nozzle_exact.tikz} \label{line:nozzle_trkplt:exsln}
\input{dat/nozzle_Mi0p4_Mo0p45_T0p8_epslm1em07_xskel_nel0008_p1.dev_from_mean.dev_from_unifmsh_0p1.trkplt.tikz}

\nextgroupplot[xticklabels={,,}, yticklabels={,,}]
\input{dat/nozzle_exact.tikz}
\input{dat/nozzle_Mi0p4_Mo0p45_T0p8_epslm1em07_xskel_nel0016_p1.dev_from_mean.dev_from_unifmsh_0p1.trkplt.tikz}

\nextgroupplot[xticklabels={,,}, yticklabels={,,}]
\input{dat/nozzle_exact.tikz}
\input{dat/nozzle_Mi0p4_Mo0p45_T0p8_epslm1em07_xskel_nel0032_p1.dev_from_mean.dev_from_unifmsh_0p1.trkplt.tikz}

\nextgroupplot[xlabel={$x$}]
\input{dat/nozzle_exact.tikz}
\input{dat/nozzle_Mi0p4_Mo0p45_T0p8_epslm1em07_xskel_nel0008_p2.dev_from_mean.dev_from_unifmsh_0p1.trkplt.tikz}

\nextgroupplot[xlabel={$x$}, yticklabels={,,}]
\input{dat/nozzle_exact.tikz}
\input{dat/nozzle_Mi0p4_Mo0p45_T0p8_epslm1em07_xskel_nel0016_p2.dev_from_mean.dev_from_unifmsh_0p1.trkplt.tikz}

\nextgroupplot[xlabel={$x$}, yticklabels={,,}]
\input{dat/nozzle_exact.tikz}
\input{dat/nozzle_Mi0p4_Mo0p45_T0p8_epslm1em07_xskel_nel0032_p2.dev_from_mean.dev_from_unifmsh_0p1.trkplt.tikz}

\end{groupplot}
\end{tikzpicture}
 \caption{Results of discontinuity-tracking framework (objective function
          defined in (\ref{eqn:obj}) with $\alpha = 0.1$) applied to the
          quasi-1d Euler equations in (\ref{eqn:nozzle-eqn}) with
          $|\Ecal_{h,p}| = 8$ (\emph{left column}), $|\Ecal_{h,p}| = 16$
          (\emph{middle column}), $|\Ecal_{h,p}| = 32$ (\emph{right column})
          elements and polynomial order $p = 1$ (\emph{top row}) and
          $p = 2$ (\emph{bottom row}).
          Legend: reference solution to the quasi-1d Euler equations determined
          from aggressive adaptive mesh refinement with over $80000$ $p=1$
          elements (\ref{line:nozzle_trkplt:exsln}),
          mesh (\ref{line:trkplt:mesh}) and solution (\ref{line:trkplt:soln})
          obtained from discontinuity-tracking framework applied to quasi-1d
          Euler equations.}
 \label{fig:nozzle-trkplt}
\end{figure}

Convergence of the proposed discontinuity-tracking method to a reference
solution using a sequence of increasingly refined meshes with polynomial
orders $p = 1,\,2$ is provided in Figure~\ref{fig:nozzle-cnvg}. The reference
solution is computed using DG on a highly adapted mesh near the discontinuity
consisting of over $80000$ $p=1$ elements. Similar to the previous sections,
the $L^1$ error is used to quantify the error in the solution. To ensure the
mesh regularization term does not impede convergence, the regularization
parameter is set to zero ($\alpha = 0$) and the single degree of freedom
parametrization of (\ref{fig:mshparam-1dof-1d}) is used.
Figure~\ref{fig:nozzle-cnvg} shows the expected convergence rate for
polynomial orders $p = 1,\,2$ and a highly accurate solution is
possible with as few as $20$ $p=2$ elements. The discontinuity-tracking method
is shown to outperform the adaptive mesh refinement method with $p=1$, at
least in terms of number of elements for a required accuracy, which is limited
to second-order accuracy.

\begin{figure}
 \centering
 \begin{tikzpicture}
\begin{axis}[
    xmin=1, xmax=1000,
    xmode=log,
    ymode=log,
    width=0.8\textwidth,
    height=0.5\textwidth,
    xlabel={Number of elements},
    ylabel={$L^1$ error}]

\addplot [black, solid, thick, mark=*, mark size=2, only marks, mark options={solid}]  table[x index=0, y index=6, select coords between index={0}{8}] {dat/nozzle_Mi0p4_Mo0p45_T0p8_epslm1em07_midnode_p1.dev_from_mean.cnvg.dat}; \label{line:nozzle_cnvg:p1}
\addplot [black, solid, thin]
coordinates {
(4.000000e+00, 1.927106e-01)
(4.000000e+02, 2.223102e-05)}; \label{line:nozzle_cnvg:slp2}

\addplot [red, solid, thick, mark=square*, mark size=2, only marks, mark options={solid}]  table[x index=0, y index=6, select coords between index={0}{8}] {dat/nozzle_Mi0p4_Mo0p45_T0p8_epslm1em07_midnode_p2.dev_from_mean.cnvg.dat}; \label{line:nozzle_cnvg:p2}
\addplot [red, densely dashed, thin]
coordinates {
(4.000000e+00, 3.121176e-02)
(7.000000e+01, 1.384733e-05)}; \label{line:nozzle_cnvg:slp3}


\addplot [blue, solid, thick, mark=x, mark size=2, only marks, mark options={solid}]  table[x index=0, y index=6, select coords between index={0}{5}] {dat/nozzle_Mi0p4_Mo0p45_T0p8_epslm1em07_fvm_amr_cnvg.dat}; \label{line:nozzle_cnvg:fvm}

\end{axis}

\end{tikzpicture}
 \caption{Convergence of discontinuity-tracking method using the single
          degree of freedom parametrization in (\ref{eqn:mshparam1d-disc})
          and no mesh regularization ($\alpha = 0$) applied to the quasi-1d
          Euler equations (\ref{eqn:nozzle-eqn}) for polynomial orders $p = 1$
          (\ref{line:nozzle_cnvg:p1}) and $p = 2$ (\ref{line:nozzle_cnvg:p2}).
          The expected convergence rates of $p+1$ are obtained in both cases.
          The slopes of the best-fit lines to the data points in
          the asymptotic regime are:
          $\angle\,-1.97$ (\ref{line:nozzle_cnvg:slp2}) and
          $\angle\,-2.70$ (\ref{line:nozzle_cnvg:slp3}).
          A reference second-order method that uses adaptive mesh refinement
          with $p=1$ elements (\ref{line:nozzle_cnvg:fvm}).}
 \label{fig:nozzle-cnvg}
\end{figure}
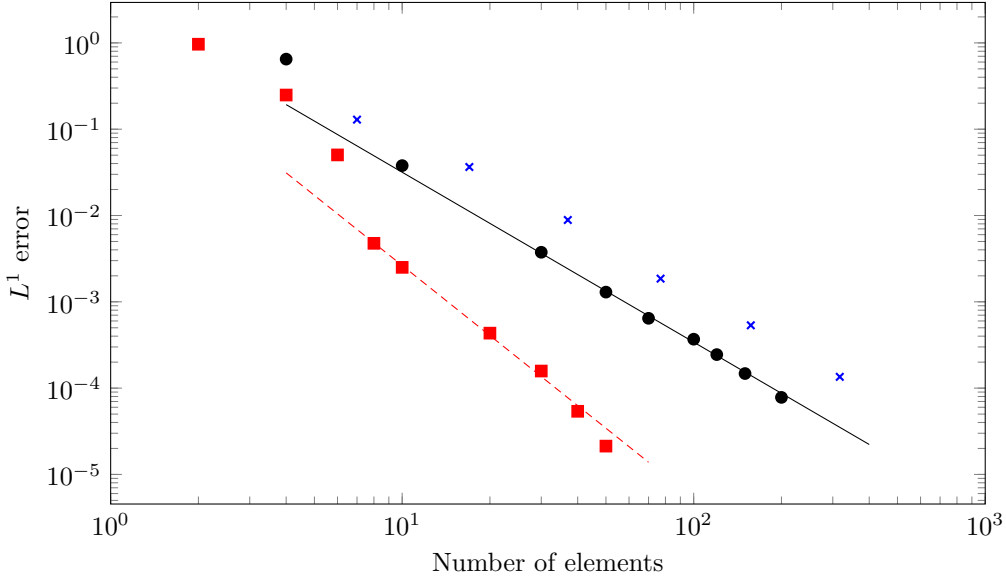
}

\subsection{Supersonic, inviscid flow around cylinder}
\label{sec:num-exp:invisc-2d}
With the merit of the proposed discontinuity tracking framework established
in one spatial dimension, we turn our focus to its performance on
two-dimensional conservation laws, specifically the two-dimensional steady
Euler equations
\begin{equation} \label{eqn:euler2d}
 \begin{aligned}
  \frac{\partial}{\partial x_j}(\rho u_j) &= 0, \\
  \frac{\partial}{\partial x_j} (\rho u_i u_j+ p)  &= 0
  \quad\text{for }i=1,2, \\
  \frac{\partial}{\partial x_j} \left(u_j(\rho E+p)\right) &= 0,
 \end{aligned}
\end{equation}
where $\rho$ is the fluid density, $u_1,\,u_2$ are the velocity
components, and $E$ is the total energy. For an ideal gas,
the pressure $p$ has the form
\begin{equation}
 p = (\gamma-1)\rho \left( E - \frac12 u_k u_k\right),
\end{equation}
where $\gamma$ is the adiabatic gas constant. The discretization
of the governing equations proceed according to the  formulation in
Section~\ref{sec:disc} and Roe's approximate Riemann solver corresponding
to the time-dependent version of (\ref{eqn:euler2d}) is used for
the numerical flux. All farfield and inviscid wall boundary conditions are
weakly imposed using the Roe solver.

The two-dimensional problem considered is supersonic flow around a circle at
Mach $2$. Due to symmetry, only one quarter of the domain is modeled; see
Figure~\ref{fig:cylshk1-geom} for a schematic of the domain and
boundary conditions.
\ifbool{fastcompile}{}{
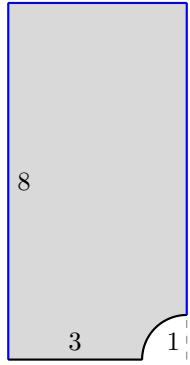
\begin{figure}
 \centering
 \begin{tikzpicture}
\begin{axis}[
axis equal image,
axis lines=none,
xmin=-0.8,
ymin=-0.8,
ymax=8.8,
xmax=4.8]
\addplot [fill opacity=0.5, fill=black!30!white, forget plot]
coordinates {
(  0.00000000,   0.00000000)
(  3.00000000,   0.00000000)
(  3.00012587,   0.01586596)
(  3.00050346,   0.03172793)
(  3.00113266,   0.04758192)
(  3.00201332,   0.06342392)
(  3.00314522,   0.07924996)
(  3.00452808,   0.09505604)
(  3.00616154,   0.11083820)
(  3.00804519,   0.12659245)
(  3.01017856,   0.14231484)
(  3.01256111,   0.15800140)
(  3.01519225,   0.17364818)
(  3.01807130,   0.18925124)
(  3.02119755,   0.20480667)
(  3.02457021,   0.22031053)
(  3.02818843,   0.23575894)
(  3.03205130,   0.25114799)
(  3.03615784,   0.26647381)
(  3.04050703,   0.28173256)
(  3.04509776,   0.29692038)
(  3.04992888,   0.31203345)
(  3.05499918,   0.32706796)
(  3.06030738,   0.34202014)
(  3.06585214,   0.35688622)
(  3.07163207,   0.37166246)
(  3.07764571,   0.38634513)
(  3.08389154,   0.40093054)
(  3.09036800,   0.41541501)
(  3.09707346,   0.42979491)
(  3.10400623,   0.44406661)
(  3.11116455,   0.45822652)
(  3.11854664,   0.47227107)
(  3.12615062,   0.48619674)
(  3.13397460,   0.50000000)
(  3.14201659,   0.51367739)
(  3.15027457,   0.52722547)
(  3.15874647,   0.54064082)
(  3.16743015,   0.55392006)
(  3.17632342,   0.56705986)
(  3.18542405,   0.58005691)
(  3.19472974,   0.59290793)
(  3.20423816,   0.60560969)
(  3.21394691,   0.61815899)
(  3.22385354,   0.63055267)
(  3.23395556,   0.64278761)
(  3.24425043,   0.65486073)
(  3.25473555,   0.66676900)
(  3.26540829,   0.67850941)
(  3.27626596,   0.69007901)
(  3.28730583,   0.70147489)
(  3.29852511,   0.71269417)
(  3.30992099,   0.72373404)
(  3.32149059,   0.73459171)
(  3.33323100,   0.74526445)
(  3.34513927,   0.75574957)
(  3.35721239,   0.76604444)
(  3.36944733,   0.77614646)
(  3.38184101,   0.78605309)
(  3.39439031,   0.79576184)
(  3.40709207,   0.80527026)
(  3.41994309,   0.81457595)
(  3.43294014,   0.82367658)
(  3.44607994,   0.83256985)
(  3.45935918,   0.84125353)
(  3.47277453,   0.84972543)
(  3.48632261,   0.85798341)
(  3.50000000,   0.86602540)
(  3.51380326,   0.87384938)
(  3.52772893,   0.88145336)
(  3.54177348,   0.88883545)
(  3.55593339,   0.89599377)
(  3.57020509,   0.90292654)
(  3.58458499,   0.90963200)
(  3.59906946,   0.91610846)
(  3.61365487,   0.92235429)
(  3.62833754,   0.92836793)
(  3.64311378,   0.93414786)
(  3.65797986,   0.93969262)
(  3.67293204,   0.94500082)
(  3.68796655,   0.95007112)
(  3.70307962,   0.95490224)
(  3.71826744,   0.95949297)
(  3.73352619,   0.96384216)
(  3.74885201,   0.96794870)
(  3.76424106,   0.97181157)
(  3.77968947,   0.97542979)
(  3.79519333,   0.97880245)
(  3.81074876,   0.98192870)
(  3.82635182,   0.98480775)
(  3.84199860,   0.98743889)
(  3.85768516,   0.98982144)
(  3.87340755,   0.99195481)
(  3.88916180,   0.99383846)
(  3.90494396,   0.99547192)
(  3.92075004,   0.99685478)
(  3.93657608,   0.99798668)
(  3.95241808,   0.99886734)
(  3.96827207,   0.99949654)
(  3.98413404,   0.99987413)
(  4.00000000,   1.00000000)
(  4.00000000,   8.00000000)
(  0.00000000,   8.00000000)
(  0.00000000,   0.00000000)};

\addplot [gray, dashed, forget plot]
coordinates {
(  4.00000000,   0.00000000)
(  4.00000000,   1.00000000)};

\addplot [thick, color=blue]
coordinates {
(  0.00000000,   8.00000000)
(  0.00000000,   0.00000000)};\label{line:cylshk1:farfield}

\addplot [thick, color=black]
coordinates {
(  0.00000000,   0.00000000)
(  3.00000000,   0.00000000)};\label{line:cylshk1:wall}

\addplot [thick, color=black]
coordinates {
(  3.00000000,   0.00000000)
(  3.00012587,   0.01586596)
(  3.00050346,   0.03172793)
(  3.00113266,   0.04758192)
(  3.00201332,   0.06342392)
(  3.00314522,   0.07924996)
(  3.00452808,   0.09505604)
(  3.00616154,   0.11083820)
(  3.00804519,   0.12659245)
(  3.01017856,   0.14231484)
(  3.01256111,   0.15800140)
(  3.01519225,   0.17364818)
(  3.01807130,   0.18925124)
(  3.02119755,   0.20480667)
(  3.02457021,   0.22031053)
(  3.02818843,   0.23575894)
(  3.03205130,   0.25114799)
(  3.03615784,   0.26647381)
(  3.04050703,   0.28173256)
(  3.04509776,   0.29692038)
(  3.04992888,   0.31203345)
(  3.05499918,   0.32706796)
(  3.06030738,   0.34202014)
(  3.06585214,   0.35688622)
(  3.07163207,   0.37166246)
(  3.07764571,   0.38634513)
(  3.08389154,   0.40093054)
(  3.09036800,   0.41541501)
(  3.09707346,   0.42979491)
(  3.10400623,   0.44406661)
(  3.11116455,   0.45822652)
(  3.11854664,   0.47227107)
(  3.12615062,   0.48619674)
(  3.13397460,   0.50000000)
(  3.14201659,   0.51367739)
(  3.15027457,   0.52722547)
(  3.15874647,   0.54064082)
(  3.16743015,   0.55392006)
(  3.17632342,   0.56705986)
(  3.18542405,   0.58005691)
(  3.19472974,   0.59290793)
(  3.20423816,   0.60560969)
(  3.21394691,   0.61815899)
(  3.22385354,   0.63055267)
(  3.23395556,   0.64278761)
(  3.24425043,   0.65486073)
(  3.25473555,   0.66676900)
(  3.26540829,   0.67850941)
(  3.27626596,   0.69007901)
(  3.28730583,   0.70147489)
(  3.29852511,   0.71269417)
(  3.30992099,   0.72373404)
(  3.32149059,   0.73459171)
(  3.33323100,   0.74526445)
(  3.34513927,   0.75574957)
(  3.35721239,   0.76604444)
(  3.36944733,   0.77614646)
(  3.38184101,   0.78605309)
(  3.39439031,   0.79576184)
(  3.40709207,   0.80527026)
(  3.41994309,   0.81457595)
(  3.43294014,   0.82367658)
(  3.44607994,   0.83256985)
(  3.45935918,   0.84125353)
(  3.47277453,   0.84972543)
(  3.48632261,   0.85798341)
(  3.50000000,   0.86602540)
(  3.51380326,   0.87384938)
(  3.52772893,   0.88145336)
(  3.54177348,   0.88883545)
(  3.55593339,   0.89599377)
(  3.57020509,   0.90292654)
(  3.58458499,   0.90963200)
(  3.59906946,   0.91610846)
(  3.61365487,   0.92235429)
(  3.62833754,   0.92836793)
(  3.64311378,   0.93414786)
(  3.65797986,   0.93969262)
(  3.67293204,   0.94500082)
(  3.68796655,   0.95007112)
(  3.70307962,   0.95490224)
(  3.71826744,   0.95949297)
(  3.73352619,   0.96384216)
(  3.74885201,   0.96794870)
(  3.76424106,   0.97181157)
(  3.77968947,   0.97542979)
(  3.79519333,   0.97880245)
(  3.81074876,   0.98192870)
(  3.82635182,   0.98480775)
(  3.84199860,   0.98743889)
(  3.85768516,   0.98982144)
(  3.87340755,   0.99195481)
(  3.88916180,   0.99383846)
(  3.90494396,   0.99547192)
(  3.92075004,   0.99685478)
(  3.93657608,   0.99798668)
(  3.95241808,   0.99886734)
(  3.96827207,   0.99949654)
(  3.98413404,   0.99987413)
(  4.00000000,   1.00000000)};\label{line:cylshk1:wall}

\addplot [thick, color=blue]
coordinates {
(  4.00000000,   1.00000000)
(  4.00000000,   8.00000000)};\label{line:cylshk1:farfield}

\addplot [thick, color=blue]
coordinates {
(  4.00000000,   8.00000000)
(  0.00000000,   8.00000000)};\label{line:cylshk1:farfield}

\node[above]    at    (axis cs:1.5, 0.0) {$3$};
\node[right]    at    (axis cs:0.0, 4.0) {$8$};
\node[]    at    (axis cs:3.7, 0.4) {$1$};
\end{axis}
\end{tikzpicture} \\
 \caption{Geometry and boundary conditions for supersonic flow around cylinder.
          Boundary conditions: inviscid wall/symmetry condition
          (\ref{line:cylshk1:wall}) and farfield (\ref{line:cylshk1:farfield}).}
 \label{fig:cylshk1-geom}
\end{figure}
}
For simplicity, we use a simple mesh parametrization that
does not include the position of all the nodes of the continuous high-order
mesh, but rather a well-chosen subset of these nodes. The remainder of the nodes
are determined through a linear operator that incorporates mesh smoothing.
In this work, we use linear elasticity with prescribed displacements at
the parametrized nodes and in the normal direction along domain boundaries.
For this problem, we use a mesh with only $|\Ecal_{h,p}| = 48$ elements and
polynomial degree up to $p = 4$. Figure~\ref{fig:cylshk-msh} shows these meshes
and identifies the parametrized nodes whose displacements compose
$\phibold$ and correspond to $N_\phibold = 12,\,24,\,36,\,48$ for the
$p = 1,\,2,\,3,\,4$ meshes, respectively.
\ifbool{fastcompile}{}{
\begin{figure}
 \centering
 \includegraphics[width=0.2\textwidth]{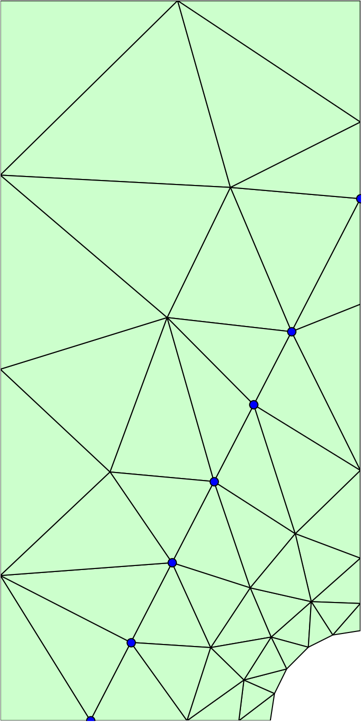} \quad
 \includegraphics[width=0.2\textwidth]{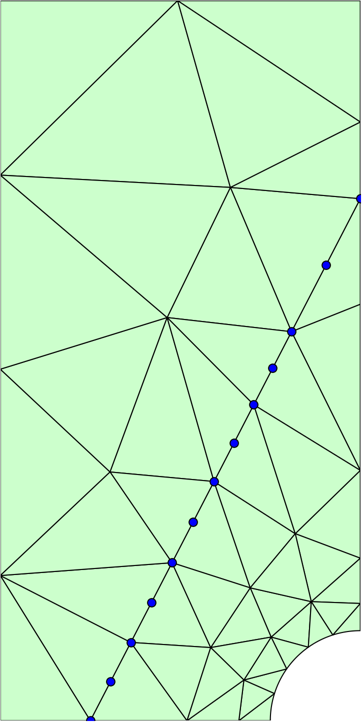} \quad
 \includegraphics[width=0.2\textwidth]{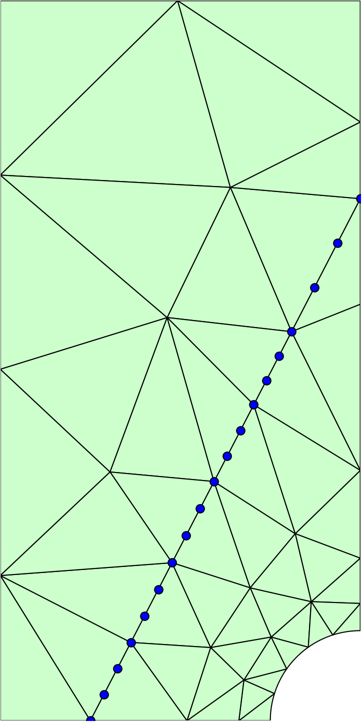} \quad
 \includegraphics[width=0.2\textwidth]{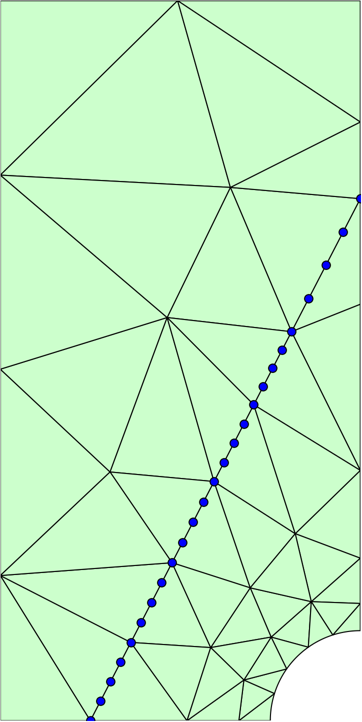}
 \caption{Reference domain for and mesh for supersonic flow around cylinder
          with $|\Ecal_{h,p}| = 48$ elements and polynomial orders $p = 1$
          (\emph{left}), $p = 2$ (\emph{middle left}), $p = 3$
          (\emph{middle right}), and $p = 4$ (\emph{right}). The blue
          circles identify parametrized nodes, i.e., nodal positions that
          whose displacements compose the optimization variables $\phibold$.
          Only the displacement normal to the boundary are taken as
          optimization variables for the two parametrized nodes that lie
          on boundaries.}
\label{fig:cylshk-msh}
\end{figure}
}

As a non-convex optimization problem underlies the discontinuity-tracking
framework, the performance of the full space solver relies on a quality
initial guess. For the case of $p = 1$, the initial guess for the mesh is
taken as the reference mesh, shown in Figure~\ref{fig:cylshk-msh}, while
for the $p > 1$ case, the mesh is initialized with the solution of the
discontinuity-tracking problem on the $p-1$ mesh. The initial state vector,
regardless of polynomial degree, is taken as a viscosity solution of
(\ref{eqn:euler2d}) with a global viscosity parameter $\nu = 0.05$, i.e., a solution of the Navier-Stokes equations, on the initial mesh, as discussed
in Section~\ref{sec:track-practical}.

The result of the discontinuity-tracking method using the above initialization
strategy and the reference domain in Figure~\ref{fig:cylshk-msh} is provided
in Figures~\ref{fig:cylshk-soln-indic}-\ref{fig:cylshk-soln}.
Figure~\ref{fig:cylshk-soln-indic} shows the high-order computational mesh
superimposed on the fluid density, element-wise discontinuity indicator
($f_{shk}$), and element-wise mesh distortion indicator ($f_{msh}$) for
the viscosity solution of the $p=1$, non-aligned mesh (e.g., initial guess
for the $p=1$ discontinuity-tracking problem) and the inviscid solution using
the discontinuity-tracking framework for the $p = 1,\,2,\,3,\,4$ meshes
(mesh regularization parameter is set to $\alpha = 0.1$).
The non-aligned mesh is clearly underresolved and the shock is smeared out.
Several elements contain non-trivial under- and over-shoot, which is
identified by the discontinuity indicator that takes large values in elements
near the discontinuity. The mesh distortion indicator is identically zero since
the reference mesh is used. It is clear that the discontinuity-tracking
framework, even for $p = 1$, considerably reduces the overshoot in the solution
near the shock at the price of mesh distortion. The $p = 1$ solution is
still clearly underresolved and the shock surface is faceted since the
element cannot bend in the piecewise linear setting. However, the
$p = 2,\,3,\,4$ reference meshes deliver quite accurate and smooth solutions
away from the shock and track the smooth, curved shock surface very well.
Since the search space is enriched in the cases relative to the $p = 1$ case,
the discontinuity indicator and mesh distortion indicator also improve in
most elements.

The bottom row of Figure~\ref{fig:cylshk-soln-indic} shows the mesh quality
marginally suffers from aligning the mesh with discontinuities, particularly
in the $p=1$ case. However, as mentioned above, the mesh quality improves as
the polynomial order increases due to the enriched search space and the
presence of the $f_{msh}$ term in the objective function. This term also
prevents intermediate optimization iterations from inverting or severely
distorting the mesh since this would cause a substantial increase in the
objective function and be rejected by the linesearch. This is significant as
an inverted or distorted mesh would cause the condition number of the
Jacobian matrix to degrade, which would in turn lead to a failure of
the linear solver.

\ifbool{fastcompile}{}{
\begin{figure}
 \centering
 \includegraphics[width=0.19\textwidth]{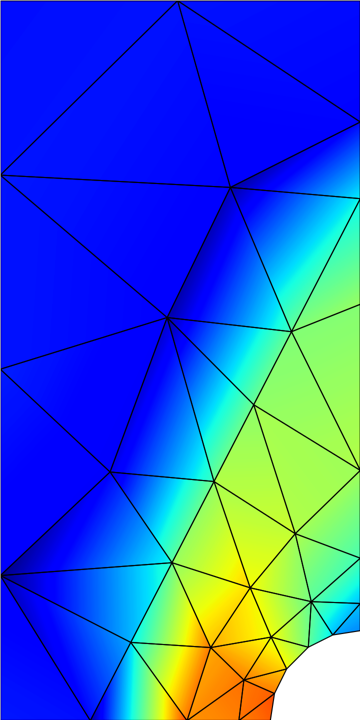} \hfill
 \includegraphics[width=0.19\textwidth]{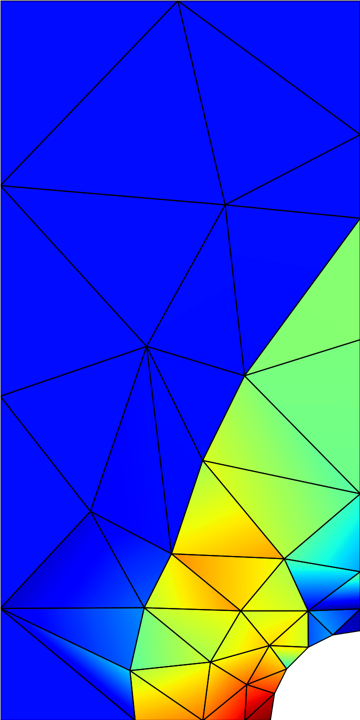} \hfill
 \includegraphics[width=0.19\textwidth]{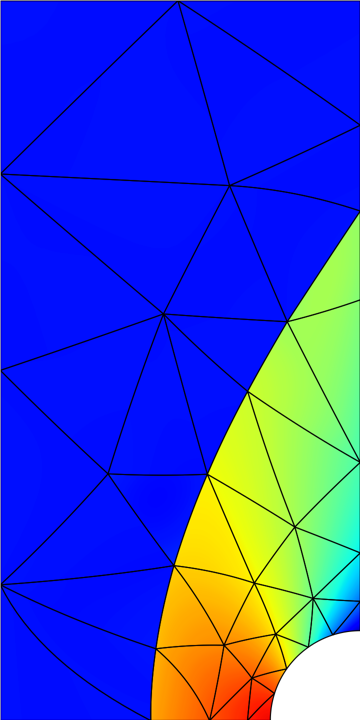} \hfill
 \includegraphics[width=0.19\textwidth]{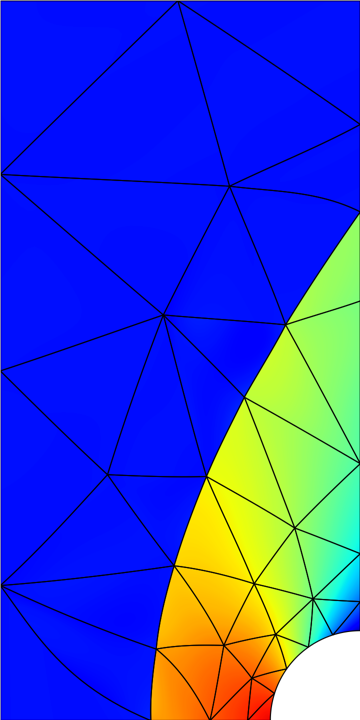} \hfill
 \includegraphics[width=0.19\textwidth]{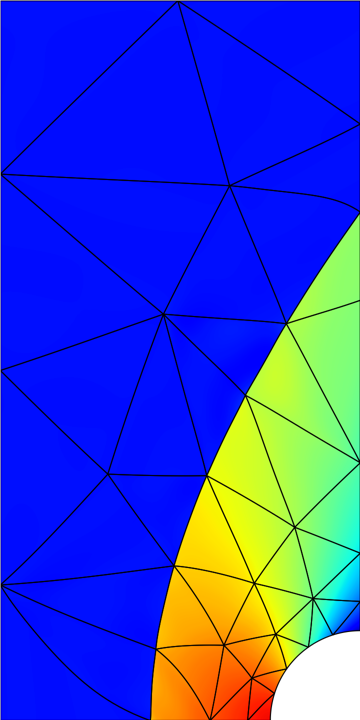} \\
 \includegraphics[width=0.19\textwidth]{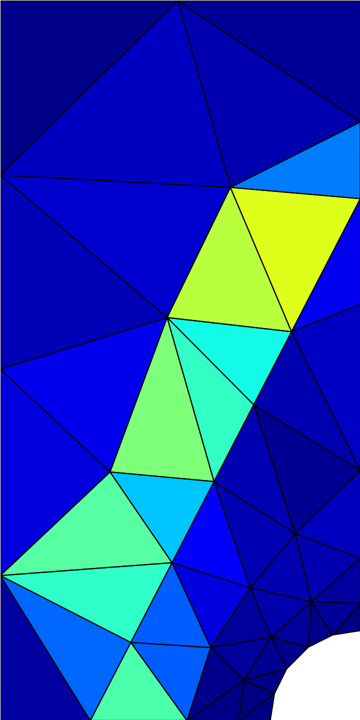} \hfill
 \includegraphics[width=0.19\textwidth]{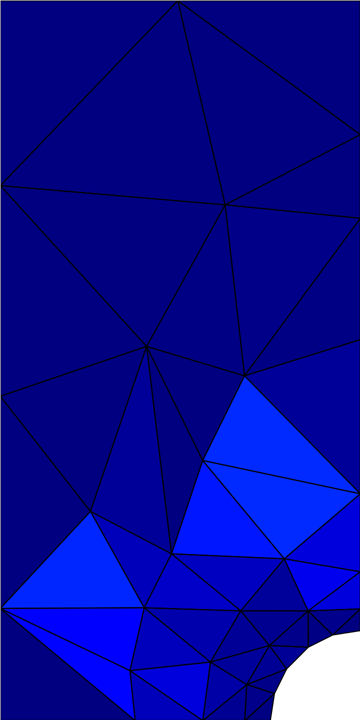} \hfill
 \includegraphics[width=0.19\textwidth]{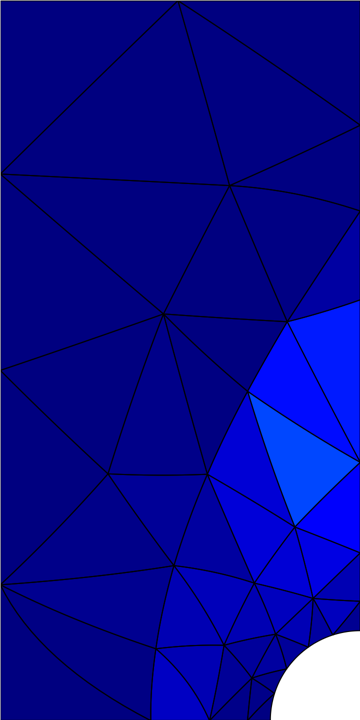} \hfill
 \includegraphics[width=0.19\textwidth]{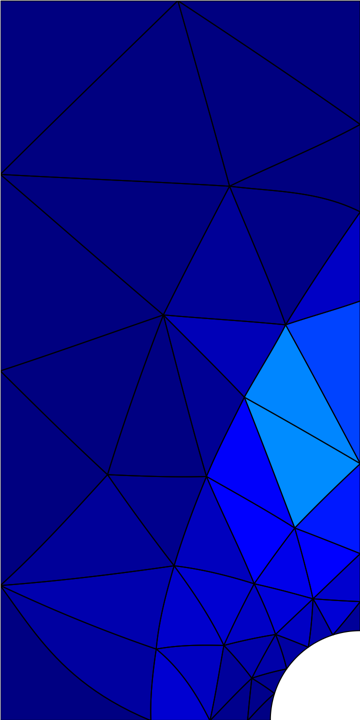} \hfill
 \includegraphics[width=0.19\textwidth]{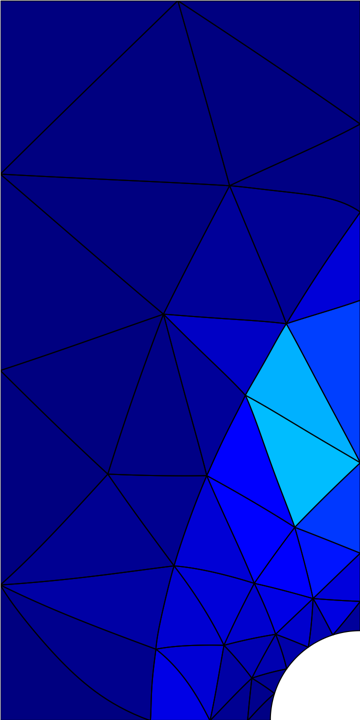} \\
 \includegraphics[width=0.19\textwidth]{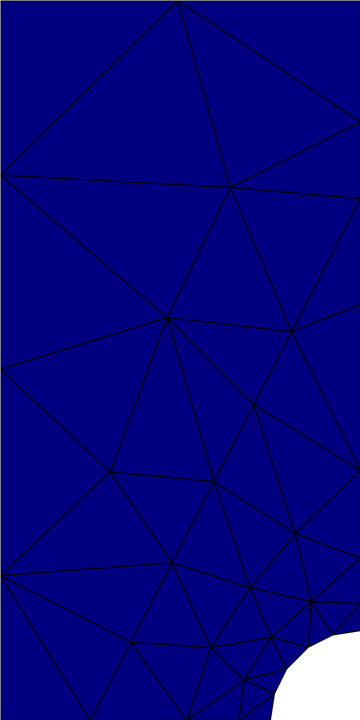} \hfill
 \includegraphics[width=0.19\textwidth]{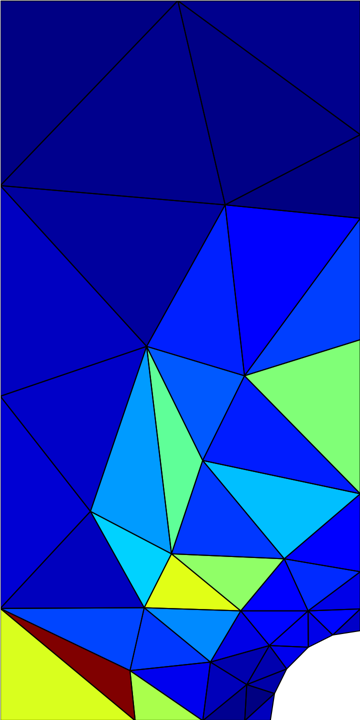} \hfill
 \includegraphics[width=0.19\textwidth]{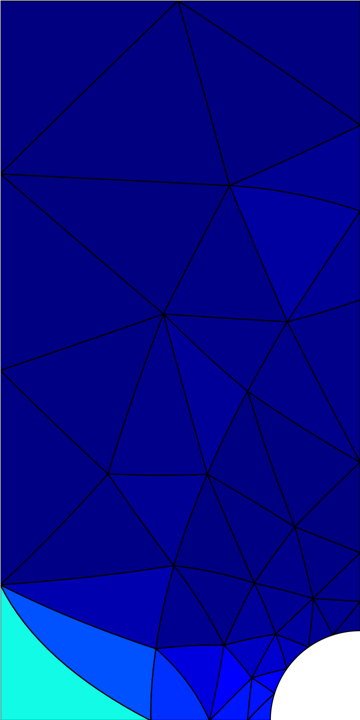} \hfill
 \includegraphics[width=0.19\textwidth]{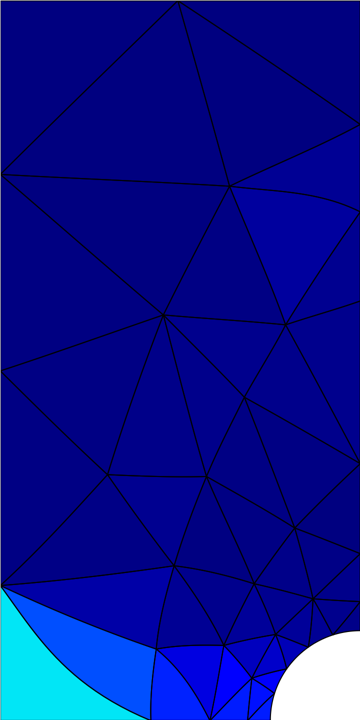} \hfill
 \includegraphics[width=0.19\textwidth]{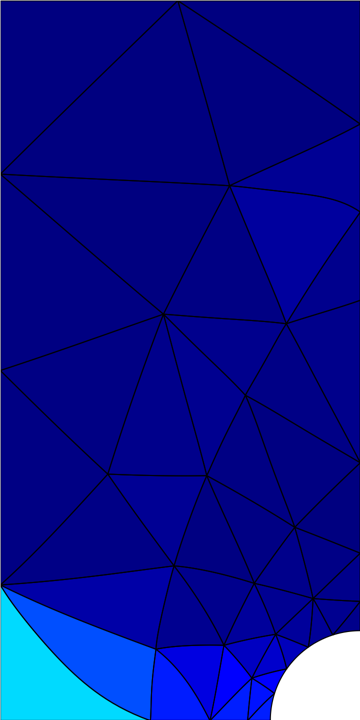}
 \caption{Discontinuity-tracking results corresponding to the reference
          meshes in Figure~\ref{fig:cylshk-msh} with mesh edges superimposed.
          \emph{Rows}: fluid density (\emph{top}), element-wise discontinuity
          indicator (\emph{middle}), element-wise mesh distortion indicator
          (\emph{bottom}). \emph{Columns}: viscosity solution on the
          non-aligned $p=1$ mesh (\emph{left}) and the inviscid,
          discontinuity-tracking solution on the $p = 1$ (\emph{middle left}),
          $p = 2$ (\emph{middle}), $p = 3$ (\emph{middle right}),
          $p = 4$ (\emph{right}) meshes with mesh regularization
          $\alpha = 0.1$.}
 \label{fig:cylshk-soln-indic}
\end{figure}
}
Figure~\ref{fig:cylshk-soln} shows the fluid density, entropy, and enthalpy
corresponding to the non-aligned, viscosity solver and the inviscid, tracking
framework without the mesh edges to emphasize the smooth, high-quality solution
obtained away from the shock.
\ifbool{fastcompile}{}{
\begin{figure}
 \centering
 \includegraphics[width=0.19\textwidth]{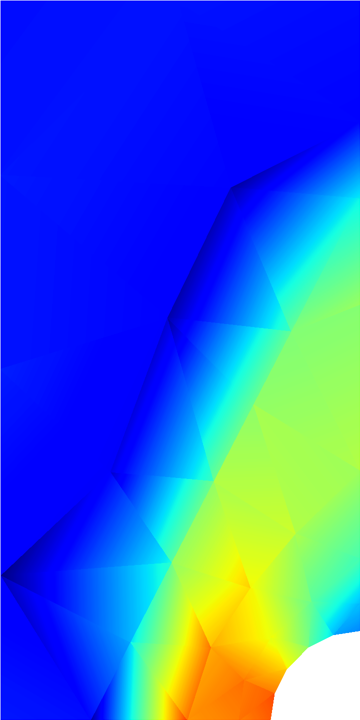} \hfill
 \includegraphics[width=0.19\textwidth]{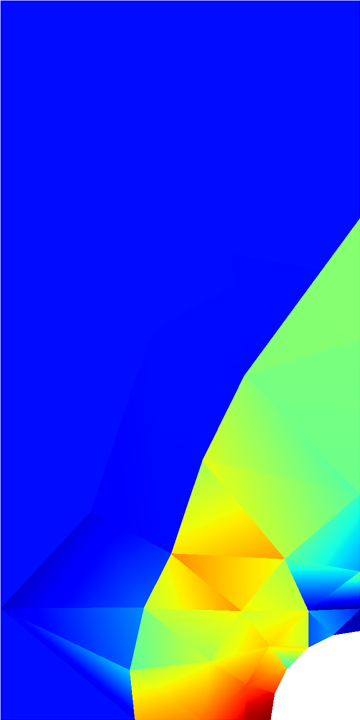} \hfill
 \includegraphics[width=0.19\textwidth]{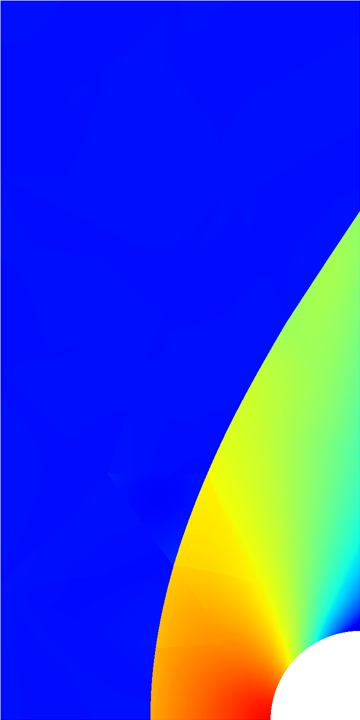} \hfill
 \includegraphics[width=0.19\textwidth]{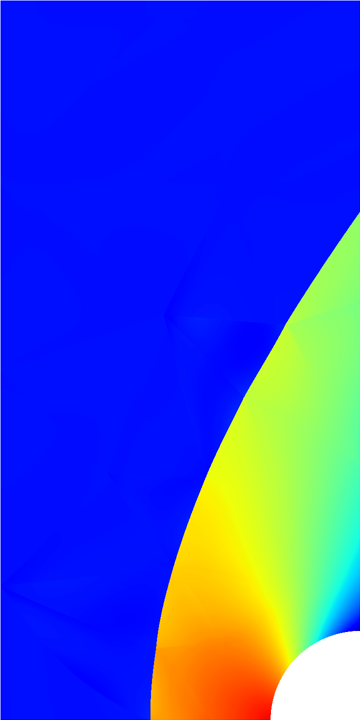} \hfill
 \includegraphics[width=0.19\textwidth]{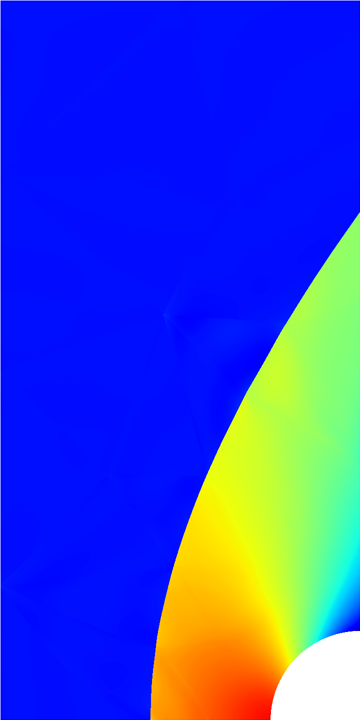} \\
 \includegraphics[width=0.19\textwidth]{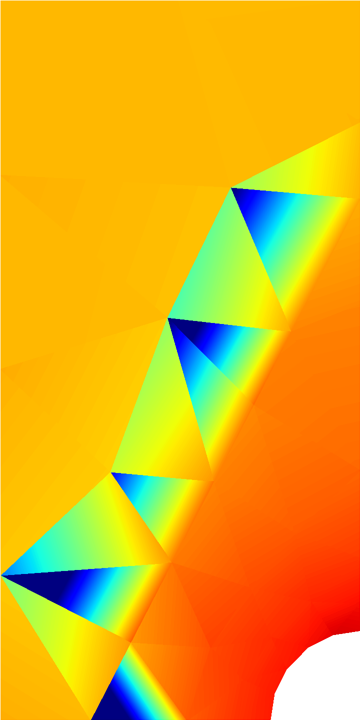} \hfill
 \includegraphics[width=0.19\textwidth]{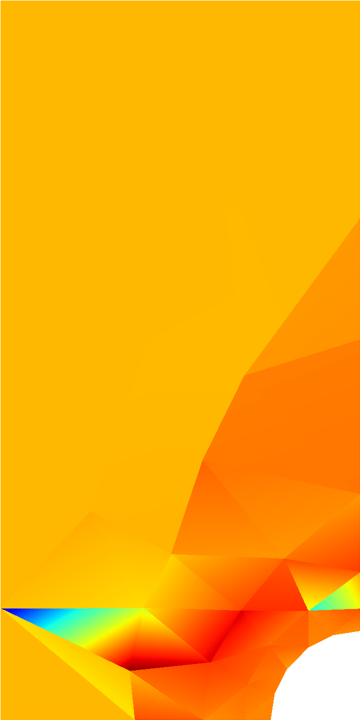} \hfill
 \includegraphics[width=0.19\textwidth]{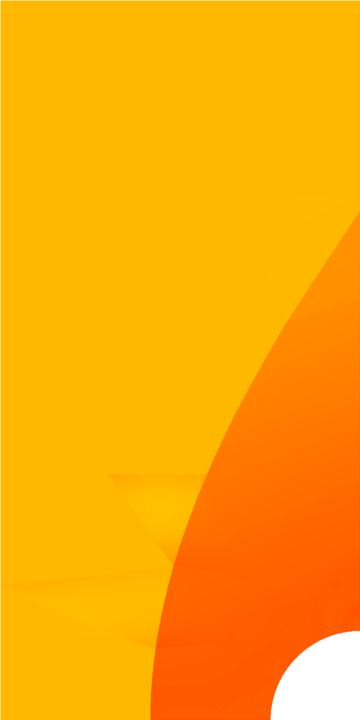} \hfill
 \includegraphics[width=0.19\textwidth]{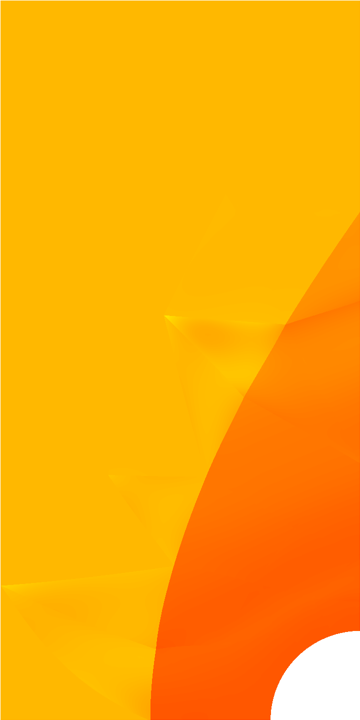} \hfill
 \includegraphics[width=0.19\textwidth]{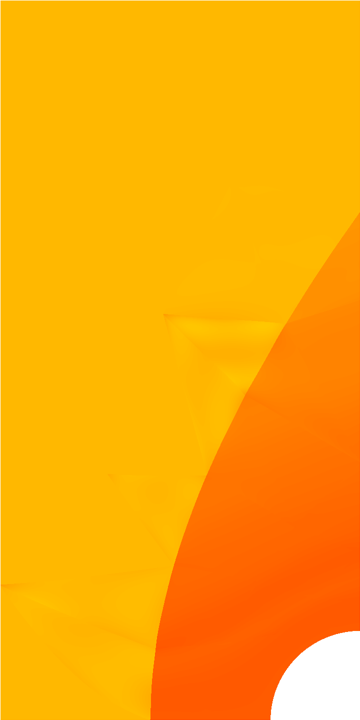} \\
 \includegraphics[width=0.19\textwidth]{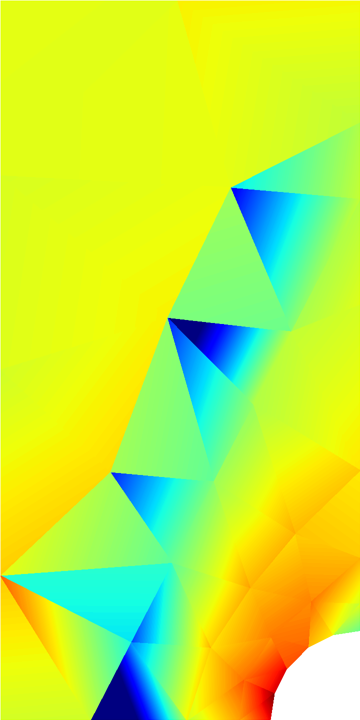} \hfill
 \includegraphics[width=0.19\textwidth]{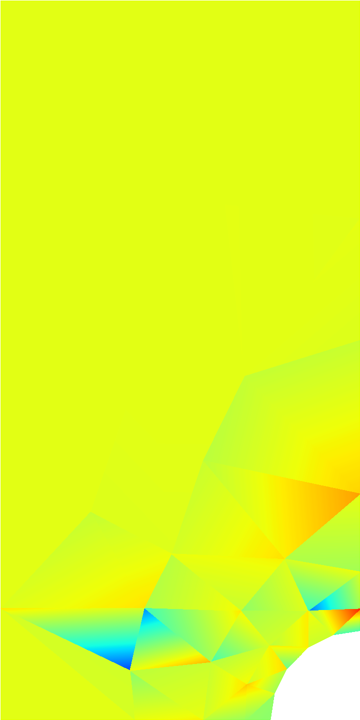} \hfill
 \includegraphics[width=0.19\textwidth]{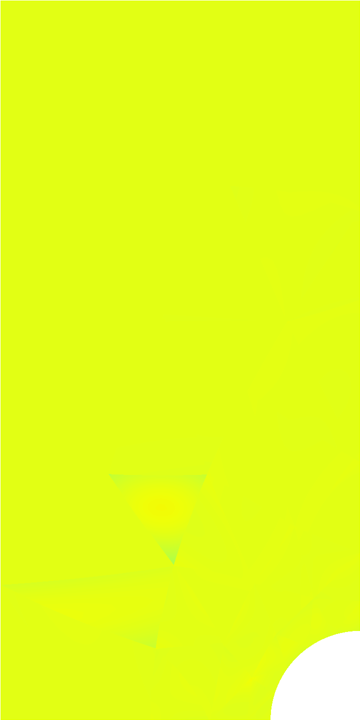} \hfill
 \includegraphics[width=0.19\textwidth]{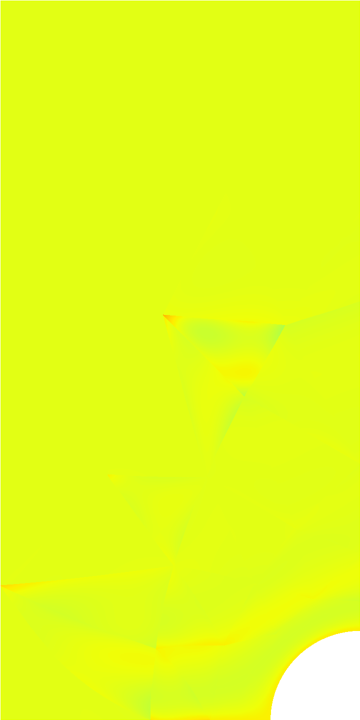} \hfill
 \includegraphics[width=0.19\textwidth]{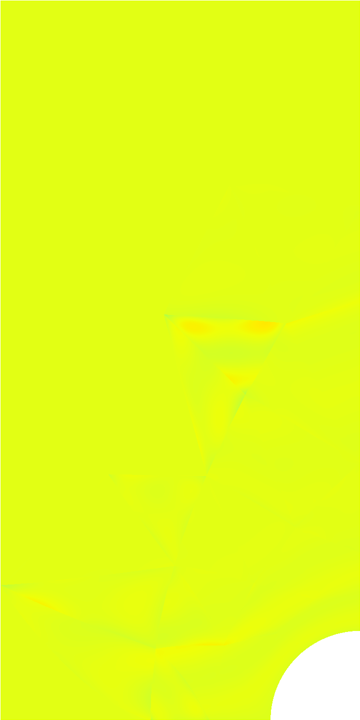}
 \caption{Fluid density (\emph{top}), entropy (\emph{middle}), and enthalpy
          (\emph{bottom}) corresponding to viscosity solution on the non-aligned
          $p = 1$ mesh (\emph{left}) and the inviscid, discontinuity-tracking
          solution on the $p = 1$ (\emph{middle left}), $p = 2$
          (\emph{middle}), $p = 3$ (\emph{middle right}), $p=4$ (\emph{right})
          meshes with mesh regularization $\alpha = 0.1$.}
\label{fig:cylshk-soln}
\end{figure}
}

To quantify the quality of the solution obtained from the discontinuity-tracking
framework, we consider two metrics from inviscid fluid mechanics: the
total enthalpy $H = (\rho E+p)/\rho$, and the stagnation pressure. In a
steady, inviscid flow the total enthalpy is constant and the stagnation
pressure, $p_0$, is given by
\begin{equation} \label{eqn:stag-pres}
 p_0 = p_\infty\frac{1-\gamma+2\gamma M_\infty^2}{\gamma+1}
       \left(\frac{(\gamma+1)^2 M_\infty^2}{4\gamma M_\infty^2-2(\gamma-1)}
       \right)^\frac{\gamma}{\gamma-1},
\end{equation}
where $p_\infty$ and $M_\infty$ are the freestream pressure and Mach number,
respectively. Therefore we use the following error metrics to quantify the
performance of the discontinuity tracking framework without requiring
a reference solution
\begin{equation}
 \begin{aligned}
  e_p(\ubm,\,\xbm) &= |p_s(\ubm,\,\xbm) - p_0| \\
  e_H(\ubm,\,\xbm) &= \sqrt{\frac{\int_\Omega (H(\ubm,\,\xbm)-H_0)^2 \, dv}
                                 {\int_\Omega dv}},
 \end{aligned}
\end{equation}
where $p_s(\ubm,\,\xbm)$ is the pressure at the stagnation point, $p_0$ is
the analytical stagnation pressure (\ref{eqn:stag-pres}), $H(\ubm,\,\xbm)$
is the total enthalpy, and $H_0$ is the farfield total enthalpy. These
quantities are summarized in Table~\ref{tab:cylshk} for the
discontinuity-tracking framework at all polynomial orders considered.
\begin{table}[ht!]
 \centering
 \begin{tabular} {c|cccc}
 Polynomial order ($p$) &  $1$  &  $2$  &  $3$  &  $4$  \\\hline
 Degrees of freedom ($N_\ubm$) &  $576$  &  $1152$  &  $1920$  &  $2880$  \\\hline
 Enthalpy error ($e_H$) &  $0.0106$  &  $0.000462$  &  $0.00151$  &  $0.000885$  \\\hline
 Stagnation pressure error ($e_p$) &  $0.0711$  &  $0.00479$  &  $0.0112$  &  $0.000616$ \\\hline
 Number of optimization iterations & 396 & 283 & 103 & 121
 \end{tabular}
 \caption{Discontinuity-tracking performance summary, including the
          number of degrees of freedom and enthalpy and stagnation pressure
          error for the reference mesh in Figure~\ref{fig:cylshk-msh} with
          $48$ elements and polynomials orders $p=1$, $p=2$, $p=3$, $p=4$
          with mesh regularization $\alpha = 0.1$. The number of
          optimization iterations are reported to quantify the computational
          expense of the full space solver. The solution and mesh are
          considered converged when the first-order optimality conditions
          of (\ref{eqn:claw-disc-opt1}) are satisfied within a tolerance of
          $10^{-6}$.}
 \label{tab:cylshk}
\end{table}
The table shows the discontinuity-tracking framework provides
accurate solutions on discretizations with very few degrees of freedom, e.g.,
$N_\ubm \sim \Ocal(10^3)$ with errors on the order of $e_H \sim \Ocal(10^{-3})$.
However, the table also shows counter-intuitive behavior in that the
enthalpy and stagnation pressure errors do not decrease monotonically
as the polynomial order is increased. The $p = 2$ errors are
exceptionally low and without this column the error decreases
monotonically. While monotonicity is expected, it is not guaranteed
by the finite element convergence theory since the mesh is extremely
coarse and the output quantities used to quantify the error are
complicated nonlinear functionals. Furthermore, the entropy plot in
Figure~\ref{fig:cylshk-soln} shows a minor imperfection along the shock
near the right boundary for the $p = 3$ and $p = 4$ solution, which
suggests the shock is not perfectly aligned with its true, physical location.
This again may be attributed to non-monotonicity of the solution error on
coarse meshes or due to the approximate Riemann solver, Roe's method with
entropy fix, not being consistent with the governing equations.

Table~\ref{tab:cylshk} also reports the number of optimization
iterations required for SNOPT to satsify the first-order optimality
conditions to within a tolerance of $10^{-6}$. Based on the discussion
in Section~\ref{sec:track-practical}, an optimization iteration is
comparable to a standard nonlinear iteration since the main additional
cost comes from the evaluation of the objective function and its
derivatives and $\displaystyle{\pder{\rbm}{\xbm}}$, which is quite small
compared to a linear solve with the Jacboian matrix. Since a $p$-homotopy
strategy is used, the cost of, e.g., the $p = 2$ solution is
$396$ iterations on the $p = 1$ mesh and $283$ iterations on the
$p = 2$ mesh. These are reasonable iteration counts given that
pseudo-transient continutation is the standard nonlinear solver for
supersonic flow problems, which usually require tens to
hundreds of iterations for convergence. In fact, the iteration count
in Table~\ref{tab:cylshk} are competitive given that our method only
requires an extremely coarse mesh with $48$ elements to deliver an
accurate solution.

\section{Conclusion}\label{sec:concl}
This paper introduces a high-order accurate, nonlinearly stable
discontinuity-tracking framework for solving steady conservation laws
with discontinuous solution features such as shock waves. The method
leverages the discontinuities between computational elements present
in the finite-dimensional solution basis in the context of a discontinuous
Galerkin or finite volume discretization to track discontinuities in the
underlying solution. The approximate Riemann solver ensures the numerical
flux is consistent with the governing equations and provides appropriate
stabilization through upwinding. Central to the proposed tracking framework
is a PDE-constrained optimization formulation of the discrete conservation
law whose objective drives the computational mesh to align with
discontinuities and constraints ensure the discrete conservation law is
satisfied. The proposed objective function not only attains its minimum when
the mesh is aligned with discontinuities, but also approaches the minimum
monotonically in a neighborhood of radius $h/2$, which is critical for
a gradient-based optimizer to locate it. A mesh regularization term is
included to avoid poor-quality meshes and is shown to not destroy the
minima or monotonicity of the indicator, for reasonable values of the
mesh quality parameter $\alpha$. The optimization problem is solved
using a full space optimization solver whereby the PDE solution and mesh
simultaneously converge to their optimal values, ensuring the discrete
PDE solution is never required on non-aligned meshes and avoiding
non-robustness issues that arise due to Gibbs' phenomena in these situations.
This simultaneous convergence behavior of the full space approach is
demonstrated in Figure~\ref{fig:l2proj1d_trkcnvg} where the discrete
conservation law residual is quite large until the objective is sufficiently
reduced, i.e., the discontinuity is tracked.

A host of numerical experiments in one and two spatial dimensions are provided
to demonstrate the merit of proposed framework. Simple one-dimensional tests
problems are provided in Section~\ref{sec:track-obj} to show the objective
function in (\ref{eqn:obj-shk}) possesses the appropriate minima and
monotonicity property, whereas other popular shock indicators have the
appropriate minima, but are highly oscillatory even when the mesh is
nearly aligned with the discontinuity. The entire tracking framework is
studied on a slew of one- and two-dimensional test problems, including the
$L^2$ projection of discontinuous functions onto a piecewise polynomial
basis, the modified steady inviscid Burgers' equation with a discontinuous
source term, the quasi-one-dimensional Euler equations (transonic, inviscid
flow through a nozzle), and the supersonic two-dimensional Euler equations.
We showed the framework successfully tracks discontinuities in the solution
even when the nodal positions of the entire mesh are taken as optimization
variables, provided the mesh regularization parameter is set appropriately.
In the two-dimensional test problems with $p > 1$, the high-order mesh
\emph{curves} to provide a high-order alignment with the shock surface. We
also showed optimal $\Ocal(h^{p+1})$ convergence rates in the $L^1$ norm up
to $p = 6$ for the $L^2$ projection problem and inviscid Burgers' equation
and showed the discontinuity tracking framework can significantly outperform
popular alternatives such as uniform or adaptive mesh refinement. Finally,
the discontinuity-tracking framework proved extremely effective in resolving
supersonic (Mach $2$) flow around a cylinder. 

While the proposed discontinuity-tracking framework shows considerable
promise as an efficient, high-order, nonlinearly stable method for resolving
discontinuous solutions of conservation laws, a number of important research
issues must be investigated for this to be a truly viable approach.
First, as discussed in Section~\ref{sec:track-solver}, there is considerable
structure in the optimization problem (\ref{eqn:claw-disc-opt1}), e.g., only
equality constraints, partition of optimization variables into $\ubm$ and
$\phibold$, invertible constraint Jacobian with respect to $\ubm$, efficient
sparse parallel solvers available for operations involving $\rbm$ and its
derivatives, that black-box optimizers cannot leverage. Therefore, there is
significant opportunity for efficiency and robustness by developing a
specialized solver that can leverage this structure and incorporate well-known
homotopy strategies into the constraint since $\rbm(\ubm,\,\xbm) = 0$ is
difficult to solve for high-speed flows that contain shocks.
Additionally, it is often the case that, during the full space iterations,
the discontinuity indicator and mesh quality can be significantly improved
by performing \emph{local topology changes} between elements, for example using the so-called edge flip and face swap operations. Future work
will focus on incorporating these discrete mesh updates into the optimization
solution procedure since it has the potential to drastically improve the
quality of the converged mesh and reduce the number of optimization iterations.
We also intend to extend the proposed framework to the \emph{time-dependent}
setting, where discrete mesh operations will be crucial.
Further research into the choice of numerical flux is also warranted since the
inter-element jumps do not tend to zero under refinement along the shock
surface, an unfamiliar situation for high-order discontinuous Galerkin methods.
While the method was only introduced and tested for
inviscid conservation laws, where perfect discontinuities arise, we expect
the framework to also be useful for viscous problems with smooth, high-gradient
solution features, but further investigation is required. Finally, this paper
only considers relatively simple model problems in one and two dimensions.
Further testing is required for more complex flows in two and three dimensions
and will be the subject of future work.

\section*{Acknowledgments}
This work was supported in part by the Luis W. Alvarez Postdoctoral Fellowship
(MZ) by the Director, Office of Science, Office of Advanced Scientific Computing
Research, U.S. Department of Energy under Contract No. DE-AC02-05CH11231
(MZ, PP) and by the AFOSR Computational Mathematics program under grant number
FA9550-15-1-0010 (MZ, PP). The content of this publication does not
necessarily reflect the position or policy of any of these supporters, and
no official endorsement should be inferred.

\bibliographystyle{plain}
\bibliography{biblio,shock_biblio,track_biblio}
\end{document}